\documentclass[bj]{imsart}

\RequirePackage{amsthm,amsmath,amsfonts,amssymb}
\RequirePackage[numbers]{natbib}

\startlocaldefs

\endlocaldefs

\usepackage[utf8]{inputenc}
\usepackage{float}
\usepackage{latexsym}
\usepackage{color}
\usepackage{graphicx}
\usepackage{hyperref}

\makeatletter
\newcommand*{\hyperlinkcite}[1]{\hyper@link{cite}{cite.#1}}
\makeatother

\newcommand{\rnc}{\renewcommand}
\newcommand{\nc}{\newcommand}
\newcommand{\mrm}{\mathrm}
\renewcommand{\b}{\textbf}
\newcommand{\bs}{\boldsymbol}
\newcommand{\wh}{\widehat}
\nc{\mb}{\mathbb}
\nc{\mc}{\mathcal}
\nc{\E}{\mb{E}}
\nc{\N}{\mb{N}}
\nc{\R}{\mb{R}}
\nc{\Q}{\mb{Q}}
\nc{\Qn}{\mb{Q}_n}
\rnc{\P}{\mb P}
\nc{\Pn}{\P_n}
\nc{\RPn}{\tilde{\P}_n}
\nc{\PnG}{\Pn^{\mc G}}
\nc{\PG}{\P^{\mc G}}
\rnc{\d}{\mrm d}
\nc{\C}{\mc{C}}
\nc{\D}{\mc{D}}
\nc{\B}{\mc{B}}
\nc{\oPo}{\stackrel{\mrm p}{\rightarrow}}
\nc{\oWo}{\stackrel{w}{\rightarrow}}
\nc{\oDo}{\rightsquigarrow}
\nc{\eff}{\|F\|}
\nc{\cred}{\color{red}}
\nc{\blue}{\textcolor{black}}
\nc{\tcb}{\textcolor{black}}
\nc{\tcr}{\textcolor{red}}

\allowdisplaybreaks[3]

\newtheorem{thm}{Theorem}
\newtheorem{lemma}{Lemma}
\newtheorem{rem}{Remark}
\newtheorem{cor}{Corollary}

\begin{document}

\begin{frontmatter}

\title{Randomized Empirical Processes by Algebraic Groups, and Tests for Weak Null Hypotheses}
\runtitle{Randomized Empirical Processes by Algebraic Groups, and Tests for Weak Null Hypotheses}

\begin{aug}
\author[A]{\fnms{Dennis} \snm{Dobler}\ead[label=e1]{d.dobler@vu.nl}}
\address[A]{Department of Mathematics, Faculty of Science, Vrije Universiteit Amsterdam, The Netherlands, \printead{e1}}
\end{aug}

\begin{abstract}
Randomization tests are based on a re-randomization of existing data to gain data-dependent critical values that lead to exact hypothesis tests under special circumstances. However, it is not always possible to re-randomize data in accordance to the physical randomization from which the data has been obained. As a consequence, most statistical tests cannot control the type I error probability. Still, similarly as the bootstrap, data re-randomization can be used to improve the type I error control. However, no general asymptotic theory under weak null hypotheses has been developed for such randomization tests yet. It is the aim of this paper to provide a conveniently applicable theory on the asymptotic validity of randomization tests with asymptotically normal test statistics. Similarly, confidence intervals will be developed.

This will be achieved by creating a link between two well-established fields in mathematical statistics: empirical processes and inference based on randomization via algebraic groups. A broadly applicable conditional weak convergence theorem is developed for empirical processes that are based on randomized observations. Random elements of an algebraic group are applied to the data vectors from which the randomized version of a statistic is derived. Combining a variant of the functional delta-method with a suitable studentization of the statistic, asymptotically exact hypothesis tests is deduced, while the finite sample exactness property under group-invariant sub-hypotheses is preserved. The methodology is exemplified with: the Pearson correlation coefficient, a Mann-Whitney effect based on right-censored paired data, and a competing risks analysis. The practical usefulness of the approaches is assessed through simulation studies and an application to data from patients suffering from diabetic retinopathy.
\end{abstract}


\end{frontmatter}

\section{Introduction}
\label{sec:intro}

Randomization methods are a powerful tool for drawing reliable statistical inferences.
\tcb{A randomization test is clearly motivated from a \emph{(physical) randomization} used in the underlying experiment.
The perhaps most famous example is Fisher's Lady Testing Tea Experiment \cite[Ch. 2]{fisher}: a certain lady should detect for eight cups of tea in which four first the milk has been added, then the tea, and in which four cups it was the other way round. 
As \cite{hemerik21} nicely pointed out, the randomization test is based on reordering all eight cups' labels, while the lady's answers are kept fixed.
This is used to check whether there was a sufficiently strong relation between the lady's answers and the actual content of the cups.
This randomization reflects the assumption that the experimenter has labeled the cups at random in the physical experiment.}

\tcb{When we speak of \emph{randomization} in this paper, we normally mean that already obtained experimental data are in some way (randomly) transformed in order to draw possible inferences on hypotheses to be tested.
The transformations, i.e.\ the randomization procedure, can sometimes be chosen based on the experimental randomization method that had been used to gain the data.
With this choice it is possible to construct exact hypothesis tests for the (sharp) null hypothesis that the sampling distribution of the data is invariant with respect to the transformations; we refer to \cite{hemerik18,hemerik21} and the references cited therein for more details.\footnote{To avoid confusion, it should be pointed out that Hemerik and Goeman used the term \emph{group invariance tests} for what we will call \emph{randomization tests}; all our tests will be based on algebraic groups.}
However, if the data can only be re-randomized in a way that does not reflect the physical randomization or if a weak null hypothesis is to be tested, 
then the randomization test will in general not be exact.
A \emph{weak null hypothesis} is here understood to be a claim about certain aspects of the sampling distribution, e.g.\ the mean or the variance, 
but not about the stricter transformation-invariance of the sampling distribution.
As we will later see, it is possible to construct randomization tests in a way to make them asymptotically exact.
It is the aim of this paper to develop a general randomization testing theory to achieve this.
The resulting tests can thus be considered as competitors to bootstrap or permutation tests.}

A popular example of a randomization method is the  permutation technique applied to the observations of two (or more) samples. 
Randomization methods date back to Fisher (see e.g.\ the discussion of the \tcb{paired} $t$-test in Section~21 of \citep{fisher}) 
and Pitman who discussed randomization by permutation in a series of papers \citep{pitmanI, pitmanII, pitmanIII}.
Permutation is usually carried out \emph{randomly} because it is computationally infeasible to realize all possible permutations;
the growth in their total number as a function of the sample sizes is super-exponentially.
\tcb{Even though random permutation and randomization are sometimes used synonymously, permutation constitutes just one particular randomization possibility.
And, to put things in the right perspective, Hemerik and Goeman pointed out that ``Fisher's famous Lady Tasting Tea experiment, which is commonly referred
to as the first permutation test, is in fact a randomisation test. This distinction is important to avoid
confusion and invalid tests.'' \cite{hemerik21}.}

One of the strong advantages of random permutation is that it results in exact tests if the samples are exchangeable.
\tcb{Often, it can also} be shown that  permutation-based inference methods asymptotically keep the significance level even under non-exchangeability; cf.\ \cite{janssen1997} for a conditional central limit theorem for the permutation version of two-sample $t$-tests.
In the $k \geq 2$ independent samples setup, \cite{ChungRomano2013} analyzed the asymptotics of permutation tests under null hypotheses beyond the case of exchangeability.
Thereby, they assumed the asymptotic linearity of estimators and proved a variant of Slutzky's theorem for randomization procedures.
\cite{chung16b} investigated permuted two-sample $U$-statistics and \cite{chung16} considered permutation tests for multivariate data in multiple samples, with applications to Hotelling's $T^2$ and a maximum statistic.
\cite{diciccio17} constructed tests for the Pearson correlation coefficient and partial correlation coefficients based on random permutation and random coordinate mirrorings.
%
%
%
%
%
%
See also \cite{pauly2015asymptotic} for extensive simulation studies concerning Wald-type permutation tests in general factorial designs,  and \cite{friedrich17} for good results of permutation methods applied to longitudinal data.
\tcb{Recently, \cite{wu20} analyzed permutation tests -- they called them \emph{Fisher randomization tests} -- for weak null hypotheses for the use in various factorial, randomized experimental design settings based on studentized quadratic forms.}
Early applications of random permutation in the
independently right-censored survival analytic context 
were developed in \cite{neuhaus1993} and \cite{janssen2001} and were extended to generalized weighted logrank permutation tests by \cite{Brendel_etal_2014}.
\cite{dobler18test} utilized studentized pooled bootstrapping and permutation techniques for constructing confidence intervals for Mann-Whitney effects in an unpaired, right-censored two-sample problem.
In general, take note of the book \cite{good05} as a source for permutation tests in various fields of application.
\tcb{Apart from a subsection, we will not consider permutation tests in this paper because their theory is sufficiently well developed.}

\tcb{Let us come back to general randomization tests.}
\tcb{Romano} \cite{romano89} made use of the group invariance assumption to construct empirical process-based Kolmogorov-Smirnov-type tests for testing independence,  spherical symmetry,  exchangeability, homogeneity, and change points; also comparisons with the bootstrap approach were made.
\tcb{In \cite{romano90} the same author} analyzed the asymptotic behaviour of randomization-based tests under broader null hypotheses which go beyond the group invariance case. He considered one- and two-sample problems, e.g.\ testing the equality of means or medians.
\cite{janssen1999nonparametric} addressed studentized randomization tests for symmetry functionals of multivariate data.
\tcb{Random coordinate permutations of bivariate data points led to} 
asymptotically exact studentized randomization-based tests for the nonparametric Behrens-Fisher problem in paired data \cite{KonietschkePauly2012}.
A general theorem for the convergence of the conditional distribution of randomized statistics was developed by \cite{hoeffding52}, Theorem~3.2.
This theorem was generalized to the multivariate case by \cite{dumbgen13}; cf. Lemma~4.1 therein.
\tcb{We also refer to} Section~15.2 in \cite{lehmann2005testing} for a collection of general properties and multiple examples and to \cite{janssen07} for a connection to the optimality of tests.
The popularity of randomization techniques has not faltered even though other competitors such as the bootstrap \citep{efron1979bootstrap} and many \tcb{of its} variants have been developed along the way.

\tcb{The bootstrap and random permutation both have} been thoroughly treated in the context of empirical process theory; see Section~3.7 in \cite{vaart96:_weak} for an overview.
Donsker theorems and functional delta-methods for permutation empirical processes provide a modern and powerful technical tool for the development of statistical inference procedures.
Until now, a similar empirical process-type theory has not been available for other randomization methods.
\tcb{One} aim of this paper is to fill this \tcb{methodological} gap.
We will develop a generally applicable randomization empirical process theory 
that allows the construction of asymptotically exact hypothesis tests in multivariate data.
At the same time, finite sample exactness of these tests will be guaranteed  for certain \tcb{sharper} hypotheses under which the data distribution is randomization-invariant.
This will be achieved by randomizing studentized test statistics and combining a conditional central limit theorem for the randomization empirical process with a new functional delta-method for randomization empirical processes.
The just mentioned studentization is in the spirit of \tcb{Janssen} \cite{janssen1997} who considered a permutation version of the two-sample $t$-test by using a suitable studentization.

The article is organized as follows.
Section~\ref{sec:ex} introduces three exemplary testing problems we will later solve with the help of different randomization approaches:
the Pearson sample correlation coefficient for paired data, a Mann-Whitney effect for right-censored paired data, and the relation of cumulative incidence functions in competing risks situations.
Empirical processes and the notion of randomization are introduced in Section~\ref{sec:eprt}.
The main results, i.e.\ a conditional weak convergence theorem and a functional delta-method for the randomization empirical process, are given in Section~\ref{sec:main}.
Also, connections to the classical bootstrap and permutation tests in two-sample problems are made.
Section~\ref{sec:excontd} revisits the previously mentioned examples and particular randomization-based tests are derived.
The practical performance of the randomization procedures is analyzed with the help of a simulation study in Section~\ref{sec:simus} and the randomization test for the Mann-Whitney effect is applied to a real data-set in Section~\ref{sec:data}.
We conclude with a discussion in Section~\ref{sec:disc}.
All proofs are given in the appendix.

\section{Three examples}
\label{sec:ex}

We motivate the use of randomization empirical processes with the help of three particular examples that we are going to revisit multiple times in the upcoming sections: the Pearson correlation coefficient, the Mann-Whitney effect for right-censored paired data, and cumulative incidence functions in competing risks situations.
Throughout the article, we write $(\Omega, \mc A, P)$ for the underlying probability space. 
We denote expectations, variances, and covariances as $E$, $var$, and $cov$, respectively.
Multivariate quantities are printed in bold-type, random quantities and some functions usually get capital letters.

\subsection{Pearson's correlation coefficient}
\label{ex:pearson}

Let $(Y,Z)$ be a bivariate random vector with positive marginal variances and the Pearson correlation coefficient 
 $\rho_{Y,Z} = {cov(Y,Z)}/[var(Y) var(Z)]^{1/2} \ \in [-1,1].$
A well-known test for independence of $Y$ and $Z$ tests whether $\rho_{Y,Z}$ is equal to zero by using a studentized version of the empirical correlation coefficient as a test statistic.
\cite{pitmanII} suggested to randomize the pairings of all $Y$ and $Z$-values, of which there are $n!$ possibilities if the sample consists of $n \in \N$ pairs. 
\cite{omelka12} applied a random permutation approach in the two-sample problem of testing equality of two correlation coefficients.
Even though random permutation is not covered by the theory developed in this article, as we are going to randomize \emph{each} data point separately, Section~\ref{sec:perm} discusses possible connections with permutation tests.

Instead of random permutation, we consider  in Section~\ref{sec:corr} the following randomization approaches to test the null hypothesis $H: \rho_{Y,Z} = 0$ against one- or two-sided alternatives: 
  random rotations around the origin, corresponding to the restricted null hypothesis of rotation invariance of the joint distribution of $(Y,Z)$, and
 random sign flips for each component, 
 corresponding to the sub-null hypothesis of joint distributions of $(Y,Z)$ that are symmetric with respect to the coordinate axes.

\tcb{It is the aim of this paper} to develop asymptotically exact tests that are also exact for finite sample sizes under the above-mentioned \tcb{sharp} sub-hypotheses.
Even though there is \emph{per se} no flaw about the permutation test that randomizes all $Y$-$Z$-pairings, 
it is well possible that the other randomization-based tests are more reliable for certain situations, at least under the respective restricted null hypotheses and if $Y$ and $Z$ are not stochastically independent.
The performance of all these tests are assessed via an extensive simulation study in Appendix~\ref{sec:simus_add}.

It should also be stressed here that this example of the Pearson correlation coefficient is actually well-known and thoroughly analyzed.
See \cite{diciccio17} for a detailed analysis of the correlation coefficient in combination with random permutation and also coordinate mirrorings.
See also Section~3.8 in \cite{good05} for a brief discussion on a permutation approach to testing for correlation.
However, in Section~\ref{sec:corr}, we will propose a test statistic that differs from DiCiccio and Romano's choice of studentization.
In addition, this example of analyzing correlation coefficients primarily serves for illustrations of the usefulness of the unified approach to randomizing empirical processes;
in essence, it is possible to use basically any reasonable randomization technique as long as integrability conditions are met and the limit distributions are not degenerate.
It is thus appealing that Theorems~\ref{thm:main} and~\ref{thm:delta_meth} below apply simultaneously to all such randomization approaches.

\subsection{Mann-Whitney effect for right-censored paired data}
\label{ex:mw}

\tcb{A classical test for the stochastic superiority of one random variable $Y$ over another, possibly related random variable $Z$ in terms of location parameters is the paired $t$-test: for a sample of i.i.d. pairs, $(Y_i, Z_i), i=1,\dots, n$, with $\mu_Y = E(Y), \mu_Z = E(Z)$, and finite and positive variances $\sigma^2_Y = var(Y), \sigma^2_Z = var(Z)$, the test for $H: \mu_Y \leq \mu_Z = 0$ versus $K: \mu_Y > \mu_Z $ is based on the statistic
$$ t_{n-1} = \frac{\sqrt{n}(\bar Y_n - \bar Z_n)}{\widehat \sigma_{Y-Z,n}} = \frac{\sqrt{n}\bar D_n}{\widehat \sigma_{D,n}}. $$
Here, $D_i = Y_i - Z_i$ are the within-pair differences, $\bar Y_n , \bar Z_n, $ and $\bar D_n$ are the obvious sample means and $\widehat \sigma_{D,n} = \widehat \sigma_{Y-Z,n}$ is the sample standard deviation of the differences.
In general, one cannot assume that $(Y_i, Z_i)$ is bivariate normally distributed.
Thus, the null distribution of the $t$-test is only asymptotically equal to the $t_{n-1}$-distribution, as $n \to \infty$.
Finitely exact testing, however, is still possible if one of the restricted null hypotheses $H_0: (Y, Z) \stackrel d = (Z,Y)$ or $H_0': (Y,Z) \stackrel{d}{=} (-Y, -Z)$ is true; note that each of $H_0$ and $H_0'$ imply that $H$ is true.
Here, $\stackrel{d}{=}$ denotes equality in distribution.
Critical values for such finitely exact tests can be based on randomly interchanging the $Y$-$Z$-labels within each pair (under $H_0$) or by multiplying the differences $D_i$ by random signs $\epsilon_i$, $P(\epsilon_i = 1) = P(\epsilon_i = -1) = \tfrac12$  (under $H_0'$).
In both cases, the resulting \emph{randomized} statistic can be written as 
$$\tilde t_{n-1} = \sqrt{\tfrac{n-1}n} \sum_{i=1}^n (\epsilon_i D_i)/ \Big(\sum_{i=1}^n (\epsilon_i D_i - \tfrac1n\sum_{j=1}^n \epsilon_j D_j)\Big)^{1/2}.$$
Note that, for fixed data points $(Y_i, Z_i)= (y_i, z_i), i=1,\dots,n$, the test statistic $t_{n-1}$ could be considered as one particular realization of $\tilde t_{n-1}$.
One can show that the test based on $t_{n-1}$ (as the test statistic) and on $\tilde t_{n-1}$ (which yields data-dependent critical values) is not only finitely exact under $H_0$ or $H_0'$ but also asymptotically exact under the more general null hypothesis $H$.}

\tcb{In other contexts, if the differences $Y_i-Z_i$ are not meaningful or if the observations are censored, the statistical analysis of the stochastic superiority of $Y$ becomes more involved.
As long as the data have the ordinal type of measurement, the following probability is still a meaningful parameter to indicate superiority, also in the two independent sample case: $p=P(Y_1 > Z_2)$.
It} is commonly estimated by the Mann-Whitney $U$-statistic \citep{mann47} or, equivalently, by the Wilcoxon rank sum statistic, and it is an easily interpretable quantity:
\tcb{for instance, let the variables $Y$ and $Z$ model the times until a bad event (e.g.\ cancer relapse) happens after the test subjects have undergone some Treatments $\mathcal{Y}$ and $\mathcal{Z}$, respectively.
If the probability that the $Y$-outcome is greater than the $Z$-outcome exceeds 50\%, then Treatment $\mathcal{Y}$ seems preferable.}
For obvious reasons, we are going to call $p$ the \emph{Mann-Whitney effect (size)}.

\tcb{Let us briefly revisit the two independent samples case.}
The two sample Wilcoxon test \tcb{is based on an estimator of $p$ and it} is particularly powerful against shift alternatives.
Allowing for ties in the data, \cite{brunner00} combined a tie-adjusted variant of the statistic with a Satterthwaite-Smith-Welch approximation for critical values that are suitable for small sample sizes and also asymptotically exact; note that they called the parameter $p$ the \emph{relative treatment effect}.
\cite{chung16b} permuted the studentized Wilcoxon test, as an example of a permuted $U$-statistic.
\tcb{Extensions} to the survival analytic context, \tcb{i.e.\ to right-censored data, were developed by}
\cite{gilbert62} and \cite{gehan65}.
Inspired by these works, \cite{efron67} extended the Mann-Whitney effect estimator to the same framework by employing Kaplan-Meier estimates instead of empirical distributions.
\cite{dobler18test} conducted a variant of Efron's test as permutation and bootstrap tests.

In the present article, we are going to extend a similar Kaplan-Meier-based test statistic to the \emph{paired two-sample} right-censored case and combine it with a randomization approach.
\tcb{Even though the classical Wilcoxon test is motivated from the two independent samples case, it can be extended for the incorporation of paired data:
the test is based on estimators for the marginal cumulative distribution functions, and those are naturally still estimable if (part of) the data are paired.
The additional information carried by the dependence within a pair might result in a power increase of the test.
In particular, the resulting test will be suitable for matched pairs studies or self-controlled case series; cf.\
\cite{Peterseni4515} for an overview of the latter study design.
Situations of partially paired data could arise as follows:
if among (paired) test subjects a part of them has flawed measurements, e.g.\ due to surgical mistakes, missing values are the consequence.
Suppose that $n_p$ pairs, and $n_1$ and $n_2$ individual subjects, who received the treatments $\mathcal{Y}$ and $\mathcal{Z}$, respectively, are eligible for an inclusion into a study.
Here, $n_1$ and $n_2$ need not be the same; the only requirement is that $n_p+n_1, n_p+n_2\geq 0$.
The statistical combination of this most general case is done in Appendix~\ref{sec:MW_uneq}.
Because the case of no pairs, $n_p=0$, is rather straightforward, we will in Section~\ref{sec:MW2} focus on the challenging novel case $n_1 = n_2 = 0$.}

We use $S_1(t) = P(T_{i1} > t)$ and $S_2(t) = P(T_{i2} > t)$, $t \geq 0$, to denote the marginal survival functions of positive random variables $T_{i1}$ and $T_{i2}$, $i=1,2$,
that are the components of independent and identically distributed (i.i.d.)\ pairs $(T_{11},T_{12})$ and $(T_{21},T_{22})$ of survival times.
We define the \emph{tie-adjusted Mann-Whitney effect} as
\begin{align}
\label{eq:mw}
p = P(T_{11} > T_{22}) + \frac12 P(T_{11} = T_{22}) = - \int S_1^{\pm} \d S_2
\end{align}
and consider the null hypothesis $H: p = \frac12$ of no treatment effect.
Here, $S_1^\pm(t) = .5 ( S_1(t) + S_1(t-)) $ denotes the normalized survival function, that is, the average of $S_1$ and its left-continuous version; cf.~\cite{dobler18test} for a derivation of the integral representation~\eqref{eq:mw}.
In many cases it seems more natural to consider $p$ instead of the within-pair-related probability $\check p = P(T_{11} > T_{12}) + \frac12 P(T_{11} = T_{12}) $ because this parameter $\check p$ might refer to a counterfactual situation in real life:
for example, \tcb{if both members of a pair relate to differently treated body parts within the same test subject, it would not make sense to treat them differently -- except of course for the purpose of the study.
Because, once it is known from the study which treatment prevails, only the superior treatment should be applied to both members of a pair henceforth; cf. Section~\ref{sec:data} about a study in which both eyes of various persons were treated differently.
Thus, the central question in this case would be ``How good are the chances for \emph{each of my eyes} not to get blind if I received Treatment $\mathcal{Y}$ rather than treatment $\mathcal{Z}$?'', which relates to the parameter $p$, and not $\check p$.
Analyses of $\check p$ might be useful if the paired test subjects consist of different individuals that have been matched based on additional characteristics such as sex, age, weight etc.
The parameter $\check p$ will be touched upon again in the discussion in Section~\ref{sec:disc}.}

The analysis of survival analytic parameters such as $p$ is typically complicated due to right-censored event times:
right-censoring renders some event times unobservable.
In such a case, the only available information is that the event of interest has not yet taken place by the time of the censoring; 
see Section~\ref{sec:MW2} for more details about how this can be dealt with statistically -- an asymptotically normal estimator of $p$ based on right-censored paired observations will be developed in that section.
We also refer to Section~11.5 in \cite{good05} for a different, permutation-based approach in the related problem of testing for stochastic ordering, i.e.\ $H': F \geq G$ against $K': F < G$, in the case of censored matched pairs.
\tcb{Note that null hypotheses formulated in terms of $F$ and $G$ are sharper than those based on $p$.}

A hypothesis test based on an estimator $\wh p$ of the Mann-Whitney effect and asymptotically valid  normal quantiles as critical values is improvable
by means of a suitable randomization technique.
In this case, we will consider random interchanges of the components of the pairs; see \cite{KonietschkePauly2012} for an application of this randomization method to the Mann-Whitney effect in the uncensored case.
Under the sub-hypothesis of exchangeability of both survival and of both censoring time components this technique will provide us with finitely exact tests.
In Corollary~\ref{thm:MW} below it will be shown that, in combination with a suitably studentized test statistic, 
this randomization approach yields critical values that converge to standard normal quantiles under $H$ even if the mentioned exchangeability does not hold.
The practical performance of the corresponding hypothesis test will be assessed in Section~\ref{sec:simus}.
\tcb{As mentioned above, a} real two sample data example about the eyes of patients suffering from diabetic retinopathy -- in which one eye of a patient received a treatment, the other not -- will be analyzed in Section~\ref{sec:data}.

\subsection{Competing risks analysis}
\label{sec:cr}

Competing risks models are often used in medical research to model and analyze the impact of various, exclusive types of events.
For example, hospital patients in intensive care units (ICU) could experience the exclusive events \emph{death in ICU} and \emph{alive discharge out of ICU} \citep[p. 1]{beyersmann11}.
Another example concerns leukemia patients
for which there are two possibilities for bone marrow transplantations:
an allogeneic transplant from a donor with a matching stem cell type, or an autologous transplant, i.e.\ the patient is his or her own donor after stem cells have been harvested.
Allogeneic transplants bear the risk of the so-called graft-versus-host disease (GvHD)  \citep{levinsky89}.
As a consequence, there is a fairly high risk that the allogeneically transplanted patient dies due to GvHD instead of a relapse.
There is thus interest to keep the risk of GvHD at bay.
That is why, even if a new treatment of a generic disease is effective in improving the survival chances, 
it should be investigated whether the risk due to a side effect does not outweigh the original disease's effects.
This can be achieved with the help of a competing risks analysis.

Mathematically speaking, such analyses use information on the event time $T$ and the random event indicator $\varepsilon$.
One is then interested in the analysis of the cumulative incidence functions $ F_j(t) = P(T \leq t, \varepsilon=j)$,
i.e.\ the probability that event type $j=1,\dots, k$ has occurred by time $t\geq 0$, where $k\geq 2$ is the total number of exclusive competing risks.
Usually, some of the event times (and then also the event types) are unobservable due to independent right-censoring.
In such cases the Aalen-Johansen estimator \citep{aalen78} can be used to estimate $F_j$.
For simplicity, let us focus on the case of $k=2$ competing risks.
We wish to analyze the relation of the cumulative incidence functions $F_1$ and $F_2$
with the help of a test for the hypotheses $H: F_1(\tau) \geq F_2(\tau)$ versus $K: F_1(\tau) < F_2(\tau)$, where $\tau>0$ is a final evaluation time-point.
For example, in leukemia research one is often interested in the $\tau = 5$ years (relapse-free) survival probability; see e.g.\ \cite{jernberg03}.
We will revisit the competing risks problem in Section~\ref{sec:cr2}.

\section{Empirical processes and a view towards hypothesis testing}
\label{sec:eprt}

From now on, to simplify notation, we will interpret all null and alternative hypotheses as a collection of distributions that satisfy the claimed property under the hypothesis.
We will primarily focus on the following multi-dimensional one-sample setup:
let $\b X_1, \dots, \b X_n$ be i.i.d.\ $d$-dimensional random vectors with distribution $\P$ and $\Pn = \frac1n \sum_{i=1}^n \delta_{\b X_i}$ be the empirical process based on this sample,
where $\delta_{\b x}$ denotes the Dirac probability measure in $\b x \in \R^d$.
Processes are indexed by a family of functions 
$$\mc F \subset \Big\{f: \R^d \rightarrow \R \text{ measurable}: \P f^2 = \int_{\R^d} f^2(\b x) \d \P(\b x) < \infty \Big\}$$
which is assumed to be a $\P$-Donsker class.
Note that one-dimensional marginals of $\Pn$ take the form
$$ \Pn f = \int_{\R^d} f(\b x) \d \Pn (\b x) =  \frac1n \sum_{i=1}^n f(\b X_i), \ f \in \mathcal{F}. $$
The following ideas are in line with the suggestion by \cite{hall91}
that ``care should be taken to ensure that even if the data might be drawn from a population that fails to satisfy $H_0$, resampling is done in a way that reflects $H_0$''.
In our case, $H_0$ will be an appropriate restriction of a general null hypothesis $H$ of interest.
To carry out the resampling through randomization, we use an algebraic group $\mc G$ acting on $\R^d$.
We assume that $\mc G$ is such that uniform sampling from $\mc G$ is possible; we equip $\mc G$ with a suitable $\sigma$-algebra and denote by $Q$ the uniform distribution on $\mc G$.
In this article, $H_0$ is always the null hypothesis of $\mc G$-invariance of the distribution $\P$, i.e.\ $H_0: \P = \tilde \P$, where 
$$\tilde \P(A) = \int_{\R^d} \int_{\mc G} 1\{g(\b x) \in A\} \d Q(g) \d \P(\b x)$$
for Borel sets $A\in \mc B(\R^d)$ characterizes the mixture distribution
and $1\{ \cdot \}$ denotes the indicator function.
We refer to Sections~6.1--6.3 in \cite{good05} for some theoretical results and examples of invariance under groups of transformations.

Examples of algebraic groups $\mc G$ that have finitely many elements are the cylcic group $\mb Z/m \mb Z$ with $m \in \mb N$ elements, the group of all component permutations, and the group that mirrors none, some, or all components of a vector with respect to the coordinate axes.
In the latter two cases, $H_0$ respectively contains all distributions with component-exchangeability and all distributions  which are symmetric with respect to all coordinate axes. 
In Section~\ref{sec:cr2} we will see an example where $\mb Z/2 \mb Z$ is utilized.
As an example group with infinite cardinality, consider the group of all length-conserving rotations in $\R^2$ around the origin:
$$\mc G = \Big\{ A_\theta = \begin{pmatrix}
              \cos \theta & - \sin \theta \\ \sin \theta & \cos \theta
             \end{pmatrix}  : \theta \in [0,2\pi) \Big\},$$
equipped with matrix multiplication as the group operation.
Here, we may draw $\theta$ uniformly from the interval $[0,2\pi)$ to obtain a random element of $\mc G$. 
In this case, $H_0$ corresponds  to all rotation-invariant bivariate distributions.
Generalizations to higher dimensions are obvious.
 
Now, let $G_1, \dots, G_n \sim Q$ be independent random objects with a uniform distribution on $\mc G$.
We define the \emph{randomization empirical process} $\RPn$ as the empirical process of $G_1(\b X_1), \dots, G_n ( \b X_n)$ which are i.i.d.\ with a distribution denoted $\tilde \P$.
That is,
$$\RPn f = \int_{\R^d} f(\b x) \d \RPn (\b x) =  \frac1n \sum_{i=1}^n f(G_i(\b X_i)) \quad \text{\tcb{and}} \quad \tcb{\tilde \P f = \int_{\R^d} \int_{\mathcal{G}} f(g(\b x)) \d Q(g) \d \P(\b x).} $$
Because every application of a randomization test shall use critical values with given fixed values of $\b X_1 = \b x_1, \dots, \b X_n = \b x_n \in \R^d$, 
we wish to analyze the conditional distribution of $(\RPn f)_{f \in \mc F}$ given $\b X_1, \dots, \b X_n$.
Note that its conditional expectation given $\b X_1, \dots, \b X_n$ \tcb{is given by}
$$ \Pn^{\mc G}f = \tcb{E(\RPn f \ \mid\ \b X_1, \dots, \b X_n) =}  \frac1{n} \sum_{i=1}^n \Big( \int_{\mc G} f(g_i(\b X_i)) \d Q(g_i) \Big) 
= \int_{\mc G} \Pn f(g(\cdot)) \d Q(g), \ f \in \mc F. $$

There are several important characteristics of the process $\RPn$ to analyze and remark.
First, a fundamental point to investigate is the asymptotic behaviour of the normalized process $\sqrt{n}(\RPn - \Pn^{\mc G})$, indexed by $\mc F$,
under both, some \tcb{weak} null hypothesis $H \supset H_0$ and the alternative hypothesis, say $K=H^c$, the complement of $H$, as $n \rightarrow \infty$, while $\b X_1, \b X_2, \dots $ are considered as fixed.
To this end, two individual requisites need to be verified: conditional convergence of all finite-dimensional marginal distributions of the normalized randomization empirical process
and its conditional tightness, both given $\b X_1, \b X_2, \dots $ in outer probability.

Second, the randomization empirical process reduces the restricted null hypothesis $H_0$ of $\mc G$-invariance to a simple hypothesis:
if $\b X_1 = \b x_1, \dots, \b X_n = \b x_n$ is a particular dataset and $t_n = T_n(\b x_1, \dots, \b x_n)$ denotes a realization of the test statistic,      
then $\tilde t_n = T_n(G_1 ( \b x_1), \dots, G_n ( \b x_n))$ can be used to exactly assess whether the number $t_n$ is extreme enough to attest  a violation of $H_0$; \tcb{see \cite{hemerik18} for more details}.

Third, a suitable studentization of the test statistic will be necessary for ensuring the asymptotic exactness of a hypothesis test in cases of no $\mc G$-invariance, i.e.\ under $H \setminus H_0$. 
As we will see in the next section, the reason for this is that the randomization procedure in general alters the (asymptotic) distribution of a statistic.

Without the first mentioned property, i.e.\ the asymptotic \tcb{Gaussianity} of the randomization empirical process irrespective of violations of the \tcb{sharp} null hypothesis $H_0$, 
the applicability of the present theory would be far too restrictive:
randomization group invariance rarely holds in real life problems, \tcb{unless the randomization method exactly reflects the physical randomization of the experiment}.
Yet, randomization methods, \tcb{e.g.\ random permutation}, are known to produce very accurate results, often even for small samples.

\section{Main results}
\label{sec:main}

In this section we analyze the asymptotic properties of the randomization empirical process.
To prepare the main statements, we denote convergence in outer probability as $\oPo$ and weak convergence on $\ell^\infty(\mc F)$ as $\oDo$, 
 as the sample size goes to infinity, i.e.\ $n\to\infty$.
We write $BL_1$ for the space of real-valued, bounded Lipschitz-continuous functions with Lipschitz-constant at most 1,
and $E_G$ denotes the conditional expectation given $\b X_1, \dots, \b X_n$ in which only $G_1, \dots, G_n$ are considered random.
$G$ and $\b X$ denote independent copies of $G_1$ and $\b X_1$, respectively. 
We use the notation $\tcb{Q_{\b x} f} = \int_{\mc G} f(g(\b x)) \d Q(g)$ to define a function \tcb{in $\b x \in \R^d$} 
and we \tcb{define} $\| \P\|_{\mc F} = \sup_{f \in \mc F}| \P f |$.

\subsection{Conditional weak convergence of the randomization empirical process}

The following main theorem explains the convergence of  randomization empirical processes as the sample size goes to infinity.
It lays the foundation for all randomization-based hypothesis tests and it gives a confirmative, yet somewhat surprising result.
For the following result it is only required that it is possible to sample uniformly from the algebraic group $\mc G$ and some Donsker properties.
\begin{thm}
\label{thm:main}
  Let $\mc F$ be $\P$- and $\tilde \P$-Donsker and $\tilde{\mc F} = \{\b x \mapsto \tcb{Q_{\b x}f} : f \in \mc F\}$ be $\P$-Donsker with 
  $\| \P\|_{\mc F}, \|\tilde \P\|_{\mc F}, \| \P\|_{\tilde {\mc F}} < \infty$.
  \tcb{Assume that the uniform distribution $Q$ on the algebraic group $\mc G$ exists.}
 Given $\b X_1, \b X_2,  \dots$, we have, as $n \rightarrow \infty$, 
 $$ \tilde{\mb G}_n = \sqrt{n}(\RPn - \Pn^{\mc G}) \oDo \tilde{\mb G} $$
 on $\ell^\infty(\mc F)$ in outer probability where $\tilde{\mb G}$ is a zero-mean Gaussian process with covariance function
 \begin{align*}
  \sigma: (f,h) & \mapsto 
   \int_{\R^d} \Big[\int_{\mc G} f(g(\b x)) h(g(\b x))  \d Q(g)  - \int_{\mc G} f(g(\b x)) \d Q(g) \int_{\mc G} h(g(\b x)) \d Q(g)\Big] \d \P(\b x) \\
  & = \P \Big[ \int_{\mc G} f(g(\cdot)) h(g(\cdot))  \d Q(g) - \Big( \int_{\mc G} f(g(\cdot)) \d Q(g) \Big) \Big( \int_{\mc G} h(g(\cdot)) \d Q(g) \Big)\Big] 
  = \P (Q_\cdot(fh) - Q_\cdot f Q_\cdot h).
 \end{align*}
\end{thm}
To be more precise, the weak convergence in the above theorem is to be understood in the following sense: as $n \to \infty$,
  $\sup_{h \in BL_1}|E_G h (\sqrt{n} (\RPn - \PnG)) - E h (\tilde{\mb G})| \oPo 0;$
see e.g.\ \cite{vaart96:_weak} for this characterization.

Theorem~\ref{thm:main} reveals that the limit process is no Brownian bridge which typically appears in classical empirical process theory.
Instead, we are here dealing with a mixture of $Q$-Brownian bridge processes.
It is interesting to note that the limit process is also fundamentally different from the Gaussian limit process of the permutation empirical process in the two-sample problem which is a Brownian bridge process; cf.\ Section~3.7.1 in \cite{vaart96:_weak}.
In the special case of exchangeability, where the distributions in both samples coincide, the Brownian bridge limit processes of the empirical process and the permutation empirical process coincide as well.
For the randomization empirical process this is in general not even the case under the \tcb{sharp} null hypothesis $H_0 : \P  = \tilde \P$ of group invariance.

Yet, randomization-based hypothesis tests are still exact under $H_0$ for the same reason why permutation tests are exact under exchangeability: the test statistic $T_n = T_n(\b X_1, \dots, \b X_n)$ and its randomization version $\tilde T_n = T_n(G_1(\b X_1), \dots, G_n(\b X_n))$ 
share the same unconditional distribution.
Hence, if the test is right-tailed, 
$$E_{H_0}(1\{ Q^n(\tilde T_n \geq T_n(\b X_1, \dots, \b X_n) \mid \b X_1, \dots, \b X_n) \leq \alpha\} ) \leq \alpha \in [0,1] $$
because $Q^n(\tilde T_n \geq T_n(\b X_1, \dots, \b X_n) \mid \b X_1, \dots, \b X_n)$ is under $H_0$ stochastically greater or equal to a uniformly distributed random variable on $(0,1)$.
\tcb{In this argument, the assumed algebraic group structure plays a prominent role; cf. Section~3 in \cite{hemerik21} for more details.}
Similarly, one can show for a randomized version of the test that the type I error probability under $H_0$ is exactly equal to $\alpha$.
Consequently, even though the Gaussian limit processes differ, 
this does not cause a problem under $H_0$.
However, if only the weak null hypothesis $H \setminus H_0$ is true, a studentization of the statistic $T_n$ is required.

\begin{rem}
 \label{rem:kiefer}
 For a better understanding of the limit Gaussian process $\tilde {\mb G}$ in Theorem~\ref{thm:main}, 
 another approach to construct this process is insightful.
 Denote by $\mb W_{\tilde \P}$ a $\tilde \P$-Brownian motion on $\R^d$, i.e.\ for $\b x,\b x'\in \R^d$, $E(\mb W_{\tilde \P,\b x} \cdot \mb W_{\tilde \P,\b x'})= \tilde \P  \, 1_{(-\bs\infty, \b x]} \cdot 1_{(-\bs\infty, \b x']} = \tilde \P ( (-\bs\infty, \min(\b x, \b x')]) $, where $-\bs\infty = (-\infty, \dots, -\infty) \in \R^d$ and the minimum $\min(\b x, \b x') \in \R^d$ is to be understood coordinate-wise.
 More generally, it can also be considered as a process with indices in $\ell^\infty(\mc F)$,
 i.e.\ a zero-mean process with the covariance function $(f,h) \mapsto E \mb W_{\tilde \P}f \cdot  \mb W_{\tilde \P} h = \tilde \P (fh)$.
 Then the process $ f \mapsto \mb W_{\tilde \P}( f - Q_{\cdot} f )$ with indices in $\ell^\infty(\mc F)$ has the same distribution as $\tilde {\mb G}$;
 see also \cite{Khmaladze17}, Equation~(3), for a similar representation in the context of empirical processes based on bivariate random variables where one of them is considered a covariate and conditioned upon.
\end{rem}

\begin{rem}
 \label{rem:invariance}
 \tcb{As mentioned above, the two limit processes $\mb G$ and $\tilde{\mb G}$ differ in general, even under randomization group invariance, i.e.\ the sharp null $H_0: \P = \tilde \P$.
 In the proof of Theorem~\ref{thm:main}, we can see the reason for this:
 for functions $f \in \mc F$, 
 $$\P [(Q_\cdot f)^2] = E \Big[ \Big( \int_{\mc G} f ( g (\b X)) \d Q(g) \Big)^2 \Big] =  E(f(G_1(\b X)) f(G_2(\b X)))$$ does in general not reduce to $(\P f)^2  \stackrel{H_0}{=} (\tilde \P f)^2 = (\P (Q_\cdot f))^2 $ under the group invariance assumption.
 Indeed, under $H_0$, $E(f(G_1(\b X)) f(G_2(\b X))) \stackrel{H_0}{=} E(f(\b X) f(G_2(\b X))) = \int f(\b x) (Q_{\b x} f) \d P(\b x)$.
 Thus, in a sense, the dependence between $f(G_1(\b X))$ and $ f(G_2(\b X))$ through $\b X$ is too strong.
 Yet, as studentizations are not strictly required under $H_0$  because the tests are conducted as conditional tests, the dependence between the test statistic and the random critical values still ensure the finite exactness of the test.}
 
 \tcb{In cases without group invariance, i.e.\ $H \setminus H_0$,
 due to the aymptotic normality of the test statistics, the variances will play the most crucial role for the asymptotic exactness of a test.
 As we will see in Section~\ref{sec:excontd}, the asymptotic variances of an unstudentized test statistic and its randomized version may or may not coincide, even under $H_0$.}
\end{rem}

\begin{rem}
 \label{rem:no_group}
 \tcb{As a referee pointed out, Theorem~\ref{thm:main} and Theorem~\ref{thm:delta_meth} below remain valid if the set $\mc G$ of operations on $\R^d$ does not have a group structure.
 The finite exactness of hypothesis tests based on the randomization empirical process with more general $\mc G$ might be lost, however, if one cannot translate the operations defined by $\mc G$ into a certain distributional invariance anymore.
 As \cite{hemerik21} argued, randomization tests that are not based on algebraic groups might still control the type I error exactly, but this depends on the experimental design of a study.}
\end{rem}

\subsection{A conditional delta-method and a studentization}

We consider a real-valued population parameter of interest, $\theta = \varphi(\P)$, for some univariate functional $\varphi: \ell^\infty(\mc F) \to \R$;
multivariate extensions  are beyond the scope of this article and will be treated in the near future.
The general two-sided hypotheses take the form $H: \theta = \theta_0$ versus $K: \theta \neq \theta_0$ 
where $\theta_0 \in \R$ is some hypothetical value that can be established through \tcb{data re-}randomization; 
one-sided tests can be obtained analogously.
Examples are offered in Section~\ref{sec:excontd} below.

For real life applications of the asymptotic result of Theorem~\ref{thm:main},
its conclusion still needs to be transferred to the real-valued parameter of interest.
Take $\hat \theta_n = \varphi(\Pn)$ as an estimator of $\theta$.
By the classical functional delta-method, we obtain \tcb{the following} asymptotically linear expansion if $\varphi$ is Hadamard-differentiable \tcb{at $\P$ with Hadamard-derivative $\varphi_\P'$ which is a continuous and linear functional}:
$$ \sqrt{n} (\varphi(\Pn) - \varphi(\P)) = \varphi_\P'(\sqrt{n}(\Pn - \P)) + o_p(1)
= \frac1{\sqrt{n}} \sum_{i=1}^n (IF_{\varphi,\P}(\b X_i) -  \mu_{\varphi, \P}) + o_p(1) \rightsquigarrow \varphi_\P'(\mb G),  $$
where $\mb G$ is a $\P$-Brownian bridge,
$IF_{\varphi,\P} : \R^d \rightarrow \R^p, \b x \mapsto \varphi_\P'(\delta_{\b x})$ is the so-called \emph{influence function}, 
 $ \mu_{\varphi, \P} = E(IF_{\varphi,\P}(\b X))$, and
$o_p(1)$ is a placeholder for a sequence of random variables 
that converge to zero in outer probability.
We denote the asymptotic variance of the random variable in the previous display by $\sigma^2_{\varphi,\P}$.

We develop a functional delta-method for the randomization empirical process
that transfers the asymptotics from the randomization empirical process $\RPn$ to $\tilde \theta_n = \varphi(\RPn)$.
The statement shall be 
$$ \sqrt{n} (\tilde \theta_n - \theta_n^{\mc G}) = \sqrt{n} (\varphi(\RPn) - \varphi(\Pn^{\mc G})) = \varphi_{\tilde \P}'(\sqrt{n}(\RPn - \Pn^{\mc G})) + o_p(1) \rightsquigarrow \varphi_{\tilde \P}'(\tilde{\mb G}) $$
conditionally on the observations in outer probability.
Here, $\theta_n^{\mc G} = \varphi(\PnG)$ is the randomization average, \tcb{i.e.\ conditional on $\b X_1, \b X_2, \dots$}.
The distribution of $ \varphi_{\tilde \P}'(\tilde{\mb G})$ will be normal.
Denote by $Y^*$ the minimal measurable majorant and by $Y_*$ the maximal measurable minorant of a random quantity $Y$.

\begin{thm}
 \label{thm:delta_meth}
  Let $\mc F$ be $\P$- and $\tilde \P$-Donsker and $\tilde{\mc F}$ be $\P$-Donsker with 
  $\| \P\|_{\mc F}, \|\tilde \P\|_{\mc F}, \| \P\|_{\tilde {\mc F}} < \infty$.
  \tcb{Assume that the uniform distribution $Q$ on the algebraic group $\mc G$ exists.}
 Let $\mb B$ be a normed space and $BL_1(\mb B)$ be the space of bounded Lipschitz-continuous functions from $\ell^\infty(\mc F)$ to $\mb B$ with Lipschitz-constant at most~1.
 Let $\varphi: \mb A_\varphi \subset \ell^\infty(\mc F) \rightarrow \mb B$ be Hadamard-differentiable at $\P$ and $\tilde \P$ tangentially to a subspace $\mb A_0 \subset \mb A$.
 Suppose $\Pn$ and $\RPn$ take values in $\mb A_\varphi$.
 Then the functional delta-method applies to the randomization empirical process in outer probability, i.e.
 \begin{align*}
  \sup_{h \in BL_1(\mb B)}|E_G h (\sqrt{n} (\varphi(\RPn) - \varphi(\PnG)) - E h (\varphi'_{\tilde \P}(\tilde{\mb G}))| & \oPo 0,
  \\
  E_G h (\sqrt{n} (\varphi(\RPn) - \varphi(\PnG))^* - E_G h (\sqrt{n} (\varphi(\RPn) - \varphi(\PnG))_* & \oPo 0 \quad \text{for all } h \in BL_1(\mb B)
 \end{align*}
as $n \rightarrow \infty$.
In addition, if $\varphi'_{\tilde \P}$ is defined and continuous on the whole space $\ell^\infty(\mc F)$, we have
 \begin{align*}
  \sup_{h \in BL_1(\mb B)}|E_G h (\sqrt{n} (\varphi(\RPn) - \varphi(\PnG)) - E_G h (\varphi_{\tilde \P}'(\sqrt{n}(\RPn - \Pn^{\mc G})))| \oPo 0.
  \end{align*}
\end{thm}
\tcb{Note, for $\mb B = \R$,} the limiting normal distribution has zero mean and variance
\begin{align}
\label{eq:sigma_tilde}
 \tilde \sigma^2_{\varphi,{\tilde \P}} = \tcb{\P (Q_\cdot IF^2_{\varphi,{\tilde \P}} - (Q_\cdot IF_{\varphi,{\tilde \P}})^2) =}  \int_{\R^d} \Big[\int_{\mc G} IF^2_{\varphi,{\tilde \P}}(g(\b x)) \d Q(g) 
 - \Big[ \int_{\mc G} IF_{\varphi,{\tilde \P}}(g(\b x)) \d Q(g) \Big]^2 \Big] \d \P(\b x) < \infty.
\end{align}
This delta-method allows the removal of the last obstacle for an asymptotically exact test for $H$ versus $K$:
because of the different limit distributions of the empirical process and the randomization empirical process,
the normally distributed random variables
$ \varphi'_{\P}({\mb G}) $ and $\varphi'_{\tilde \P}(\tilde{\mb G})$
generally also have different variances.
Therefore, it is required 
that the weak limits of $ W_n = \sqrt{n} (\varphi(\Pn) - \varphi(\P))$ and $ \tilde W_n = \sqrt{n} (\varphi(\RPn) - \varphi(\Pn^{\mc G}))$
are studentized with the help of appropriate standard deviation estimators based on $\b X_1, \b X_2, \dots$ and $G_1(\b X_1), G_2(\b X_2), \dots$, respectively.
A suitable studentization will ensure the asymptotic pivotality of the limits under the larger null hypothesis $H$ -- both limit distributions are standard normal by Slutzky's lemma -- 
and it still guarantees the finite sample exactness under the restricted null hypothesis $H_0$ of $\mc G$-invariance.

The asymptotically linear representations from the functional delta-methods motivate the following  studentizations for $ W_n$ and $\tilde W_n$, respectively:
$$
\wh \sigma_{\varphi,\P_n}^2 = \frac1n \sum_{i=1}^n \Big\{ IF_{\varphi, \Pn}(\b X_i) - \varphi'_{\Pn}(\Pn)
\Big\}^2
\quad \text{and} \quad
\tilde \sigma_{\varphi,\RPn}^2 = \frac1n \sum_{i=1}^n \Big\{ IF_{\varphi,\RPn}(G_i ( \b X_i)) - \varphi'_{\RPn}(\RPn)
\Big\}^2. $$
The influence function of a complicated functional that is possibly built up of multiple simpler functionals is derivable with the help of a chain rule; see \cite{reid81} for details.
A sufficient condition for the consistency of $\wh \sigma_{\varphi,\P_n}^2 $ and $\tilde \sigma_{\varphi,\RPn}^2 $ is that, for $k=1,2$,
\begin{align}
  \label{eq:conv_IF}
   \Pn (IF_{\varphi,\Pn}^k - IF_{\varphi,\P}^k) = o_p(1) \quad \text{and} \quad \RPn (  IF_{\varphi,\RPn}^k - IF_{\varphi,\tilde \P}^k) = o_p(1).
  \end{align}
The conditions in \eqref{eq:conv_IF} in turn hold, for example,
if the influence function  satisfies a certain Lipschitz condition; see Appendix~\ref{sec:lipschitz} for details.
Alternatively, one could obviously also verify the consistency of $\wh \sigma_{\varphi,\P_n}^2 $ and $\tilde \sigma_{\varphi,\RPn}^2 $ directly;
\tcb{see Section~\ref{sec:MW2} below for an exemplification of such an approach.}
A combination of all ingredients results in the \tcb{following} randomization test which is in the spirit of the permutation two-sample $t$-test as discussed in Lemma~4.1 of \cite{janssen1997}:
\begin{cor}
 \label{cor:stud}
   \tcb{Assume that the conditions of Theorem~\ref{thm:delta_meth} hold, that}
   $\varphi(\Pn^{\mc G}) = \theta_0$ for all $n \in \N$,
   $\sigma^2_{\varphi,{\P}}, \tilde \sigma^2_{\varphi,{\tilde \P}} > 0$, and that \eqref{eq:conv_IF} holds.
Then, as $n \rightarrow \infty$, the following test \tcb{for $H: \theta = \theta_0$ versus $K: \theta_0 \neq \theta$} has asymptotic level $\alpha \in (0,1)$ under $H$ with finite sample exactness under~$H_0$:
$$ 1\{ |T_n| > \tilde c_n \} + \tilde\gamma_n 1\{ |T_n| = \tilde c_n \} $$
where $T_n = \sqrt{n}(\varphi(\Pn) - \theta_0) / \wh \sigma_{\varphi_{\Pn}'}$, 
$\tilde c_n$ is the conditional $(1-\alpha)$-quantile of $|\tilde T_n| = \sqrt{n}|\varphi(\RPn) - \theta_0| / \tilde \sigma_{\varphi_{\RPn}'}$, and $\tilde \gamma_n =(\alpha - P(|\tilde T_n| > \tilde c_n \mid \b X_1, \b X_2, \dots))/ P(|\tilde T_n| = \tilde c_n \mid \b X_1, \b X_2, \dots)$, with $0/0:=0$.
\end{cor}


\subsection{Combination with permutation tests}
\label{sec:perm}

Multiple algebraic groups can obviously be combined to a larger group to obtain another randomization empirical process.
However, such enlargements lead to more restrictive sub-hypotheses $H_0$ for finitely exact inference.
For example, if one would combine the groups of coordinate mirrorings and rotations around the origin,
the finite exactness of hypothesis tests would only hold if $\P$ is symmetric with respect to the coordinate axes and also rotation invariant.
From this point of view, it seems preferable to choose a rather small group that still yields a non-degenerate asymptotic limit distribution of the randomized estimator and, in particular, finite exactness under a rather large \tcb{sharp} null hypothesis $H_0 \subset H$.
\tcb{On the other hand, \cite{hemerik21} argue that rather large randomization groups (or sets) lead to ``higher resolution $p$-values'' and thus to a better power of the randomization test if $\alpha$ is very small.}

Because of the enormous general interest in permutation tests for two independent samples problems,
we shall discuss possibilities for combinations of algebraic group randomization with random sample group permutations.
\tcb{As discussed above,} finitely exact inference is then achievable only for exchangeable samples that share the same group-invariant distribution.
To be precise, let $\b X_1^{(j)}, \dots, \b X_{n_j}^{(j)}$ be i.i.d.\ random vectors from two independent groups $j=1,2$ with distributions $\P^{(j)}$.
Write $\b X_1, \dots, \b X_N$ for the pooled sample, $N=n_1 + n_2$.
Permutation tests are based on random sample group interchanges:
let $\pi = (\pi(1), \dots, \pi(N))$ be a random permutation of $(1,\dots,N)$, 
then many classical permutation tests use the permuted samples 
$\b X_{\pi(1)}, \dots, \b X_{\pi(n_1)}$ and 
$\b X_{\pi(n_1+1)}, \dots, \b X_{\pi(N)}$.
A combination with group randomization can be achieved based on both permuted randomized samples,
$G_{\pi(1)} ( \b X_{\pi(1)}), \dots, G_{\pi(n_1)} ( \b X_{\pi(n_1)})$ and 
$G_{\pi(n_1+1)} ( \b X_{\pi(n_1+1)}), \dots, G_{\pi(N)} ( \b X_{\pi(N)})$.

Assume that $n_1/N \to \lambda \in (0,1)$ and write $\tilde \P^{(1)}, \tilde \P^{(2)}$ for the distributions of $G_1 ( \b X_1)$ and $G_N ( \b X_N)$, respectively.
\tcb{The randomization empirical process of interest in this two-sample context is
$$ \tilde{\mb G}^{(1,\pi)}_{n_1,n_2} = \sqrt{n_1} (\tilde \P^{(1,\pi)}_{n_1} - \bar{\tilde\P}_{n_1, n_2}^{\mc G}) 
= \sqrt{n_1}\Big(\frac1{n_1} \sum_{i=1}^{n_1} \delta_{G_{\pi(i)} ( \b X_{\pi(i)})} - \frac1N \sum_{i=1}^N \delta_{G_i ( \b X_{i})} \Big); $$
equivalently, the process based on the second permuted sample could be considered, $\tilde{\mb G}^{(2,\pi)}_{n_1,n_2}$.
Let us sketch some ideas about the weak convergence of $\tilde{\mb G}^{(1,\pi)}_{n_1,n_2}$;
the actual weak convergence is a conjecture.
A full analysis including a list of all additional requirements is beyond the scope of this article.
If we conditioned on $G_1, \b X_{1}, G_{2}, \b X_{2}, \dots,$ which is denoted by the conditional expectation $E_{\pi}$, well-known results on the permutation empirical process imply that, stated in terms of the bounded Lipschitz metric, 
\begin{align}
\label{eq:conv_perm}
 \sup_{h \in BL_1} |E_{\pi}h(\mb G^{(1,\pi)}_{n_1,n_2}) - E h (\sqrt{1-\lambda} \tilde{\mb G}_{\bar{\tilde\P}}) |
\end{align}
converges to zero in outer probability; cf.~Section~3.7.1 in \cite{vaart96:_weak}.
Here, $\mb G_{\bar{\tilde\P}}$ is a $(\lambda \tilde \P^{(1)} + (1-\lambda) \tilde \P^{(2)})$-Brownian bridge process.}
Interestingly, the random permutation thus corrects the limit distribution such that it coincides with that of the original normalized empirical processes if $\tilde \P_1 = \P_1=\P_2 = \tilde \P_2$, \tcb{despite the $\mc G$-randomization which altered the process structure in the one-sample case in Theorem~\ref{thm:main}}.
\tcb{Apart from verifying asymptotic measurability, it remains to show that \eqref{eq:conv_perm} holds with $E_{\pi}$ replaced by $E_{\pi,G}$, i.e.\ conditional on $\b X_1, \b X_2, \dots$.
The dominated convergence theorem suggests that the upper bound
\begin{align*}
  \sup_{h \in BL_1} |E_{\pi, G}h(\mb G^{(1,\pi)}_{n_1,n_2}) - E h (\sqrt{1-\lambda} \tilde{\mb G}_{\bar{\tilde\P}}) | 
  \leq E_{\pi,G} \sup_{h \in BL_1} |E_{\pi}h(\mb G^{(1,\pi)}_{n_1,n_2}) - E h (\sqrt{1-\lambda} \tilde{\mb G}_{\bar{\tilde\P}}) |
\end{align*}
converges to zero in outer probability.
Due to outer probabilities, however, care must be exercized in the correct application of the dominated convergence theorem; cf.~Problem~1.2.4 in \cite{vaart96:_weak}.}

\tcb{For a permutation-related example, we reconsider} the correlation coefficient from Section~\ref{ex:pearson}. 
We model the sample with the help of independent and identically distributed random vectors $(Y_1, Z_1), \dots, (Y_n, Z_n)$ and we denote the marginal averages by $\bar Y_n = \frac1n \sum_{i=1}^n Y_i$ and $\bar Z_n = \frac1n \sum_{i=1}^n Z_i$.
The classical permutation approach is to randomly permute only the second coordinates, $Z_1,\dots, Z_n$.
Denote the random permutation vector by $(\pi(1), \dots, \pi(n))$.
Provided that integrability conditions hold, the permuted empirical correlation coefficient converges as follows:
$$ \sqrt{n} (\rho^{\pi}_n - 0) = \sqrt{n} \cdot \frac{\sum_{i=1}^n (Y_i - \bar Y_n)(Z_{\pi(i)} - \bar Z_n)}{[\sum_{i=1}^n (Y_i - \bar Y_n)^2 \sum_{i=1}^n (Z_{\pi(i)} - \bar Z_n)^2]^{1/2}} \rightsquigarrow N(0,1) $$
conditionally on $Y_1, Z_1, Y_2, Z_2,\dots$ in probability; cf.\ Theorem~2.1 in \cite{diciccio17}.
\tcb{As the discussion above suggests, a similar convergence should hold} if one combines \tcb{some $\mc G$}-randomization and permutation;
a deduced test for correlation would thus be finitely exact under the restricted null hypothesis of group invariance and independence of $Y_1$ and $Z_1$.
After a suitable studentization, the test would also be asymptotically exact under the general null hypothesis $H: \rho_{Y,Z}=0$; we refer to Theorem~2.2 in  \cite{diciccio17} for this statement and to Section~\ref{sec:corr} below for a different studentization approach.

\subsection{Relationship to Efron's bootstrap}
\label{sec:bs}

We shall see that the classical bootstrap \citep{efron1979bootstrap} is covered by a variant of the above randomization empirical process approach for more general maps $\mc G\subset \{g: \R^{d\cdot n} \to \R^d\}$ that act on the full sample and not just on the individual random vectors.
However, this greater flexibility comes at the cost of a loss of the algebraic group structure and hence no finitely exact hypothesis tests can be established, not even under \tcb{sharp} null hypotheses $H_0 \subset H$.

To describe how Efron's bootstrap can be established this way,
let $\pi_i = (\pi_i(1), \dots, \pi_i(n))$, $i=1,\dots, n$, be independent random permutations of the numbers $1,\dots, n$, and define the random maps $G_i$ via $G_i ( (\b X_1, \dots, \b X_n))= \b X_{\pi_i(1)}$, $i=1,\dots, n$.
In a certain sense, the asymptotic covariance structure given in Theorem~\ref{thm:main} also covers the structure that results from Efron's bootstrap; consider the following finite sample variant of that covariance:
\begin{align*}
 & \P^n \Big[ \int_{\mc G} f(g(\cdot)) h(g(\cdot))  \d Q(g) - \int_{\mc G} f(g(\cdot)) \d Q(g) \int_{\mc G} h(g(\cdot)) \d Q(g)\Big]  \\
  & = \int_{\R^{d\cdot n}} \Big[  \frac1n \sum_{i=1}^n f(\b x_i) h(\b x_i) - \frac1n \sum_{i=1}^n f(\b x_i) \cdot \frac1n \sum_{j=1}^n h(\b x_j) \Big] \d \P^n(\b x_1, \dots, \b x_n) \\
  & = \frac1n \sum_{i=1}^n E(f(\b X_i)h(\b X_i)) - \frac1{n^2} \sum_{i\neq j} E(f(\b X_i)h(\b X_j)) - \frac1{n^2} \sum_{i=1}^n E(f(\b X_i)h(\b X_i)) = \frac{n-1}{n} (\P (fh) -\P f \P h)
\end{align*}
As $n\to\infty$, the $\P$-Brownian bridge structure of the bootstrap empirical process is re-established; see also Theorem~3.6.1 in \cite{vaart96:_weak} for a conditional Donsker theorem for the bootstrap empirical process.
We thus see that the classical bootstrap \tcb{can be retrieved by extending the transformations  $g \in \mc G$ such that they act on the whole sample}.
This is in contrast to the permutation approach of Section~\ref{sec:perm} because the random permutations there cannot be achieved by means of \emph{independent} random transformations.

\section{Three examples continued}
\label{sec:excontd}

\subsection{Test for correlation}
\label{sec:corr}

In this first example \tcb{we will exercize} an application of the randomization empirical process theory.
It should be kept in mind that this example has been similarly worked on by \cite{diciccio17}, but by means of a permutation test and with a different studentization.
Nevertheless, another detailed discussion here will illuminate the use of our Corollary~\ref{cor:stud}.
Let $ (Y_1, Z_1), \dots, (Y_n, Z_n) =: (Y,Z)$ be i.i.d.\ pairs of random variables with joint distribution $\P$, positive and finite marginal variances, and correlation coefficient $\rho_{Y,Z} \in (-1,1)$.
We wish to apply the developed randomization empirical process theory to test the hypotheses  $H :  \rho_{Y,Z} = 0$ against $K : \rho_{Y,Z} \neq 0$. 
A commonly used estimator for $\rho_{Y,Z}$ is the empirical correlation coefficient
 $$\wh \rho_n = \frac{\sum_{i=1}^n (Y_i - \bar Y_n) (Z_i - \bar Z_n)}
  {[ \sum_{i=1}^n (Y_i - \bar Y_n)^2 \sum_{i=1}^n (Z_i - \bar Z_n)^2 ]^{1/2}}.$$

A candidate for a randomization group is
$\mc G^{rot.} = \{ A_\theta \cdot (y,z)^t: \theta \in [0,2\pi)  \}$, the group of rotations around the origin.
It will give rise to finitely exact tests under the \tcb{sharp} null hypothesis
$ H_0^{rot.} : \{ (Y,Z)^t \stackrel{d}{=} A_\theta \cdot (Y,Z)^t \} \subset H $
of rotation invariance of $\P$.
\tcb{The resulting hypothesis test will thus be exact under the sharp null for all spherically symmetric bivariate distributions.}
Another possible choice is $\mc G^{symm.}=\{g : \R^2 \rightarrow \R^2 : g(y,z) \in \{ (y,z), (y,-z), (-y,z), (-y,-z) \}\} $,
the group of mirrorings with respect to the coordinate axes.
Let the random signs $\varepsilon_Y, \varepsilon_Z \stackrel{i.i.d.}{\sim} 2 \cdot Bin(1,.5) -1$ be independent of $Y,Z$.
Based on this group, we will obtain finite exactness under the \tcb{sharp} null hypothesis
$ H_0^{symm.} : \{ (Y,Z) \stackrel{d}{=} (\varepsilon_Y Y,\varepsilon_Z Z) \} \subset H $
of distributions $\P$ that are symmetric with respect to the coordinate axes.
\tcb{Admittedly, these two randomization groups have been chosen for illustrative purposes rather than for their motivation from a physical randomization procedure.
In Section~\ref{sec:MW2} we will encounter a randomization procedure that reflects physical randomization.}

In our further asymptotic analysis of this example, we assume without loss of generality that $E(Y)=E(Z)=0$ and $var(Y) = var(Z) = 1$ because the empirical Pearson correlation coefficient and its randomized counterpart are independent of location and scale parameters.
Next, we note that $\wh \rho_n$ can be expressed as a Hadamard-differentiable functional $\phi$ of the empirical process $\Pn$ of $(Y_i,Z_i), i=1, \dots, n$, indexed by a combination of canonical projections, $\mc F = \{ p_1, p_2, p_1^2, p_2^2, p_1 p_2 \}$.
Thus, slightly abusing the notation, for 
$$\phi: \R^5 \rightarrow \R, \ (y,z,a,b,c) \mapsto (c-yz)/[(a - y^2)(b - z^2)]^{1/2},$$
we have the representation
$\wh \rho_n = \phi ( \Pn ) $
with
$y = \Pn p_1 = \bar Y_n,$ 
$z = \Pn p_2 = \bar Z_n,$
$a = \Pn p_1^2 = \frac1n \sum_{i=1}^n Y_i^2,$  
$b = \Pn p_2^2 = \frac1n \sum_{i=1}^n Z_i^2,$
$c = \Pn ( p_1 p_2) = \frac1n \sum_{i=1}^n Y_i Z_i.$
\tcb{The delta-method yields}
$$\sqrt{n}(\wh \rho_n - \rho_{Y,Z}) =
\phi'_{\rho_{Y,Z}} (\sqrt{n} (\Pn - \P)) =  \frac1{\sqrt{n}} \sum_{i=1}^n ( Y_i Z_i - \frac{\rho_{Y,Z}}2 Y_i^2 - \frac{\rho_{Y,Z}}2 Z_i^2) + o_p(1)$$
where $\phi'_{\rho_{Y,Z}}(\P) = \int (yz - \frac{\rho_{Y,Z}}{2} y^2 - \frac{\rho_{Y,Z}}{2} z^2 ) \d \P(y,z) = 0$ and 
\begin{align}
 \label{eq:if_corr}
 IF_{\phi_\rho,\P}(Y,Z) = \phi'_{\rho_{Y,Z}}(\delta_{(Y,Z)}) = YZ - \frac{\rho_{Y,Z}}{2} (Y^2 + Z^2).
\end{align}
Now, a simple application of the central limit theorem readily yields the following asymptotic behaviour;
we will use the notation
$\sigma^2_V = var(V)$ and $\sigma_{V,W} = cov(V,W)$ for square-integrable real random variables $V$ and $W$.
\begin{lemma}
\label{lem:cor_weak}
 If $E(Y^4 + Z^4) < \infty$, we have
 $
  \sqrt{n} ( \wh \rho_n - \rho_{Y,Z}) \oDo N(0,\sigma_\rho^2) \ \text{as } n \rightarrow \infty
 $ 
 where 
 \begin{align}
 \label{eq:asy_var_pear}
 \sigma_\rho^2 = 
 \sigma^2_{\check Y \check Z} - \rho_{Y,Z} (\sigma_{\check Y \check Z, \check Y^2} + \sigma_{\check Y \check Z, \check Z^2}) + \frac{\rho_{Y,Z}^2}4 (\sigma^2_{\check Y^2} + \sigma^2_{\check Z^2} + 2 \sigma_{\check Y^2, \check Z^2}) 
 \end{align}
 for the standardized random variables
 $\check Y = \sigma_Y^{-1} (Y - E(Y))$ and $\check Z = \sigma_Z^{-1} (Z - E(Z))$.
\end{lemma}
Even though Lemma~\ref{lem:cor_weak} suffices to build a randomization-based hypothesis test for correlation, there is room for improvement.
In general, an application of the Fisher z-transformation seems appealing because it stabilizes the asymptotic variance under normality: 
by the delta-method, we have that
\begin{align*}
  \sqrt{n} ( \tanh^{-1}(\wh \rho_n) - \tanh^{-1}(\rho_{Y,Z})) \oDo N \Big(0,\frac{\sigma_\rho^2}{(1-\rho_{Y,Z}^2)^2} \Big) \quad \text{as } n \rightarrow \infty,
 \end{align*}
 where the asymptotic variance reduces to 1 if the underlying distribution is bivariate normal, irrespective of the actual value of $\rho_{Y,Z}$;
 see \cite{diciccio17} for similar observations and the recommendation to conduct a permutation test for $H: \rho_{Y,Z} = 0$ based on a studentized version of 
$\sqrt{n} \tanh^{-1}(\wh \rho_n)$.
In their Section 2, they proposed to divide this statistic by 
$$ \widehat \tau_n^2 = \frac{\frac1n \sum_{i=1}^n (Y_i - \bar Y_n)^2 (Z_i - \bar Z_n)^2}{[{\frac1n \sum_{j=1}^n (Y_j - \bar Y_n)^2 \frac1n \sum_{k=1}^n (Z_k- \bar Z_n)^2}]} $$
which, under $H: \rho_{Y,Z} = 0$, results in $\sqrt{n} \tanh^{-1}(\wh \rho_n) / \widehat \tau_n^2 \oDo N(0,1)$ as $n\to\infty$.
Under local alternatives and bivariate normality of the data, 
their resulting permutation test has a pivotal limiting power.
For non-normal data, however, the asymptotic variance in~\eqref{eq:asy_var_pear} reveals that
the statistic $\sqrt{n} ( \tanh^{-1}(\wh \rho_n) - \tanh^{-1}(\rho_{Y,Z}))  / \widehat \tau_n^2$ does in general not have a pivotal asymptotic variance under local alternatives.
In order to achieve just this, we propose to choose the following statistic instead:
$ T_n(\rho) = \sqrt{n} {(1-\wh\rho_n^2)}/{\wh\sigma_{\rho,n}}  \cdot ( \tanh^{-1}(\wh \rho_n) - \tanh^{-1}(\rho) ), $
where $\wh\sigma_{\rho,n}^2$ is an estimator of $\sigma_\rho^2$ in \eqref{eq:asy_var_pear} that involves the obvious moment-type estimators.
Note that no additional moment conditions are required for its consistency.
The test statistic for $H: \rho_{Y,Z} = 0 $ versus $K: \rho_{Y,Z} \neq 0$ is given by $T_n(0)$.

In the next subsections we will determine the asymptotic variances of the randomized empirical Pearson correlation coefficients based on both groups $\mc G^{rot.}$ or $\mc G^{symm.}$. 
Write 
$$\tilde T_n  = \sqrt{n} \frac{1-\tilde\rho_{n}^2}{\tilde\sigma_{\rho,n}} ( \tanh^{-1}(\tilde \rho_{n}) - \tanh^{-1}(\rho_{Y,Z}^{\mc G}) )$$ for the randomization version of $T_n(0)$,
where $\tilde\rho_{n}$ and $\tilde\sigma_{\rho,n}^2$ are derived in the same way as $\wh\rho_{n}$ and $\wh\sigma_{\rho,n}^2$, just based on the randomized random vectors. 
Simple calculations show that  $\rho_{Y,Z}^{\mc G} = \phi(\Pn^{\mc G}) = 0$ for both above-metioned examplary choices of $\mc G$.
For simplifying the presentation, we continue to assume that $E(Y)=E(Z)=0$ and $var(Y) = var(Z) = 1$ without loss of generality. 


\subsubsection{Randomization of the Pearson correlation based on vector rotation}

We first reconsider $\mc G^{rot.}$, the group of rotations around the origin
and express the rotations of the vectors $(y,z)$ more conveniently as $r(y,z) \cdot (\cos \theta, \sin \theta)$, with radius $r(y,z) = \sqrt{y^2 + z^2}$ and angle $\theta = \theta(y,z) \in [0,2\pi)$.
We see that the asymptotic variance~\eqref{eq:sigma_tilde} of $\sqrt{n}\tilde \rho_n$ equals
$$ \int_{\R^2} \Big[ \frac1{2 \pi} \int_0^{2\pi} \cos^2(\theta) \sin^2(\theta) \d \theta -  \Big\{ \frac1{2 \pi} \int_0^{2\pi} \cos(\theta) \sin(\theta) \d \theta \Big\}^2 \Big] \cdot r^4(y,z) \d \P(y,z); $$
keep in mind here the particular form \eqref{eq:if_corr} of the influence function and that $\rho_{Y,Z}^{\mc G} = 0$.
The term in curly brackets vanishes because the integrand is an odd function. 
The remaining inner integral simplifies due to the double-angle formulas: $\cos^2(\theta) \sin^2(\theta) = \sin^2(2 \theta) / 4 = \{1 - \cos(4\theta)\} / 8$.
Hence, the integral above reduces to $\tilde\sigma^{2}_{rot.} := E\{(Y^2 + Z^2)^2\}/8 > 0$.
\tcb{Under the null hypothesis of no correlation, according to Lemma~\ref{lem:cor_weak}, the asymptotic variance of  $\sqrt{n}(\wh \rho_n - \rho_{Y,Z})$ is equal to $E(Y^2 Z^2)$.
Neither under the sharp null hypothesis of rotation invariance nor under the independence of $Y$ and $Z$ does this in general coincide with $\tilde\sigma^{2}_{rot.}$.
However, in the very special case of i.i.d. normally distributed $Y, Z$ the asymptotic equality of variances holds.}

\subsubsection{Randomization of the Pearson correlation based on coordinate mirroring}

Similarly, the group $\mc G^{symm.}$ of coordinate mirrorings leads to the following  asymptotic variance~\eqref{eq:sigma_tilde}:
$$ \int_{\R^2} [ \frac14 \{y^2z^2 + (-y)^2z^2 + y^2(-z)^2 + (-y)^2(-z)^2 \} -
\{ \frac14 (yz -yz -yz + yz) \}^2
] \d \P(y,z)  $$
\tcb{which equals $E(Y^2Z^2) > 0$. Under $H: \rho_{Y,Z} = 0$, this coincides with the asymptotic variance of the normalized sample correlation.
Hence, a studentization is not strictly necessary for the test based on the mirroring randomization procedure, in contrast to the rotation- and even the permutation-based tests; cf. \cite{diciccio17}.
Furthermore, straightforward computations reveal that the involved limiting Gaussian processes $\mb G$ and $\tilde{\mb G}$ even share the same distribution under $H_0^{symm.}$.}

\subsubsection{Final remarks}

Denote by $\tilde c_n(\alpha), \alpha\in(0,1)$, the conditional $(1-\alpha)$-quantile of $| \tilde T_n|$ given $(Y_1, Z_1), (Y_2,Z_2), \dots$, the corresponding randomization probability by $\tilde \gamma_n(\alpha)$,
and the standard normal cumulative distribution function by $\Phi$.
If one of the empirical variance estimators $\wh \sigma_{Y,n}^2$ or $\wh \sigma_{Z,n}^2$ is equal to 0,
set the test statistic $T_n(0)$ to 0, and proceed similarly for $\tilde T_n$ if one of the randomized empirical variances is zero.
We arrive at the following corollary for hypothesis tests for correlation which holds, e.g., for $\mc G = \mc G^{rot.}$ or $\mc G^{symm.}$;
denote by $H_0^{\mc G}$ the corresponding \tcb{sharp} null hypothesis of randomization group invariance. We denote by $G\sim Q$ a random group element and by $\| \cdot \|_2$ the Euclidean norm on $\R^2$.
\begin{cor}
 \label{thm:pearson}
 Assume that $\rho_{Y,Z} \neq \pm 1$, $E(Y^4 + Z^4) < \infty$,  and that
$\mc G$ is such that $E\|G(Y,Z)\|_2^2 < \infty$.
 Then, for $n \rightarrow \infty$, we have under $H \cup K$ that
 $\tilde c_n(\alpha)$ converges to $\Phi^{-1}(1-\alpha/2)$ in outer probability. Furthermore, the test
 $$ \Psi_n = 1\{ |T_n(0)| > \tilde c_n(\alpha)\} + \tilde \gamma_n(\alpha) 1\{ |T_n(0)| = \tilde c_n(\alpha)\}  $$
 satisfies $ E(\Psi_n) \to  1_{K} + \alpha 1_{H}$ as $n \to \infty$. 
 Additionally, under $H_0^{\mc G} \subset H : \{\rho_{Y,Z} = 0\}$, the test has level $\alpha$ for finite sample sizes $n \in \N$.
\end{cor}

\noindent\tcb{Even though not all randomization procedures discussed in this section enjoy a motivation from an underlying physical experiment,
the simulation study in the appendix below demonstrates that most resampling-based tests for $H: \{\rho_{Y,Z}=0\}$ still exhibit a reasonably good type I error control when the sharp hypotheses are not true.}

\tcb{As a last special case, we consider bivariate normally distributed $(Y,Z)$.
Because $H: \{\rho_{Y,Z}=0\}$ implies that $(Y,Z)$ is spherically symmetric, it is symmetrically distributed with respect to the coordinate axes, and also that $Y$ and $Z$ are independent, we conclude that all three considered randomization tests (based on mirroring, rotations, and the permutation test) are finitely exact.}

\subsection{The Mann-Whitney effect for right-censored paired data}
\label{sec:MW2}

Let $(T_{11}, T_{12}), \dots, (T_{n1}, T_{n2})$ be i.i.d.\ pairs of positive survival times and $ (C_{11}, C_{12}), \dots, (C_{n1}, C_{n2})$  i.i.d.\ pairs of positive censoring times.
We again denote the survival functions of $T_{ij}$ by $S_j$, $j=1,2$.
Let $\tau > 0$ be the final evaluation time for which we assume 
\begin{align}
 \label{eq:cond_mw}
 P(\min(T_{1j}, C_{1j}) > \tau) >0, \quad j=1,2.
\end{align}
The actually observable data consist of the survival or the censoring times, whatever comes first, i.e.\ $(X_{i1}, X_{i2}) = (\min(T_{i1}, C_{i1}), \min(T_{i2}, C_{i2}))$, 
and the censoring indicators
$(\delta_{i1}, \delta_{i2}) = (1\{T_{i1} \leq  C_{i1} \}, 1\{T_{i2} \leq C_{i2}\})$.
The marginal Kaplan-Meier estimators for the survival functions $S_j$ are given by
$$ \wh S_{j,n}(t) = \prod_{u \leq t}(1 - {\sum_{i: X_{ij} = u} \delta_{ij}}/{\sum_{i} 1\{X_{ij} \geq u\} }), \ j=1,2. $$
The factors in the above product are different from 1 only for a finite number of different values of~$u$.
It is well-known that the Kaplan-Meier estimator is a Hadamard-differentiable functional of the empirical process of the survival times and the censoring indicators;
see \cite{vaart96:_weak}, Example~3.9.31, for the empirical process-based weak convergence result for $\wh S_{j,n}$.

In our case, where two survival functions need simultaneous estimation, we will use the empirical process $\Pn$ of the data $(X_{i1}, X_{i2}, \delta_{i1}, \delta_{i2}),$ 
$i=1,\dots, n, $ indexed by 
$\mc F = \{ 1_{(0, s]}(p_j) \cdot p_{j+2}, \  1_{[s, \infty)} (p_j): j=1,2, \ s \in [0,\tau]  \}$.
The Mann-Whitney effect introduced in Section~\ref{ex:mw} is then estimated with the help of both Kaplan-Meier estimators.
This quantity, restricted to the time interval $[0,\tau]$, i.e.\
$$p = P(\min(T_{11},\tau) > \min(T_{22},\tau)) + \frac12 P(\min(T_{11},\tau) = \min(T_{22},\tau)) = -\int_0^\tau S_1^\pm(u) \d S_2(u)$$
 is estimated based on the truncated data:
overwriting previous notation,
$(X_{i1}, X_{i2}) = (\min(T_{i1}, C_{i1}, \tau), $ 
$\min(T_{i2}, C_{i2},\tau))$, 
and $(\delta_{i1}, \delta_{i2}) = (1\{\min(T_{i1},\tau) \leq  C_{i1} \}, 1\{\min(T_{i2},\tau) \leq C_{i2}\})$,
and the Kaplan-Meier estimations will be based on these truncated data;
see \cite{dobler18test} for more details on the truncation at $\tau$.
Now, the estimated Mann-Whitney effect is 
 $\wh p_n = - \int \wh S_{1,n}^{\pm}\d \wh S_{2,n};$
see, for the independently right-censored, two independent samples case, \cite{dobler18test} and \cite{efron67}, Section~8, for a similar estimator in the case of continuous $S_1$ and $S_2$ and with $\tau=\infty$.
Note that Efron, in order to achieve a ``self-consistency property'' of the estimators $\wh S_{j,n}$, 
set the Kaplan-Meier estimators at their largest event times to zero, irrespective of whether those were an event or a censoring.
This is actually in agreement with our Kaplan-Meier estimators based on the observations truncated at $\tau$ because all such truncated points are marked as ``uncensored''.
Hence, the Kaplan-Meier estimators are forced to take the value 0 at $\tau$ if there is at least one such truncation.

The estimator $\wh p_n$ results from combining the modified Wilcoxon functional $\psi: (f,g) \mapsto \int f^\pm(u) \d g(u)$ with the pair of both Kaplan-Meier estimators. 
%
\tcb{If one assumes that the underlying sample was obtained from a random treatment assignment, $j=1$ or $j=2$, then a treatment group re-randomization of the data seems most sensible as it reflects the physical randomization procedure. Hence, as}
a randomization group to randomize the Mann-Whitney effect estimator, we propose to use $\mc G^{exch.} = \{(x_1,x_2,d_1,d_2) \mapsto (x_1, x_2, d_1, d_2) , (x_1,x_2,d_1,d_2) $ $\mapsto (x_2,x_1,d_2,d_1) \}. $ 
\tcb{This} allows to interchange the sample group correspondence within each observed pair of survival times and also the corresponding censoring indicators.
See \cite{KonietschkePauly2012} for  a similar approach for inference about $p$  in the uncensored paired case.
This choice results in the restricted null hypothesis of sample group exchangeability 
$$H_0^{exch.} : \{ (X_1,X_2,\delta_1,\delta_2) \stackrel{d}{=} (X_2,X_1,\delta_2,\delta_1) \} \ \subset \  H: \{p=.5\}.$$
$H_0^{exch.}$ is true if, for example, the pairs of survival times and also the pairs of censoring times are exchangeable, i.e.\ $(T_1,T_2) \stackrel{d}{=} (T_2,T_1)$ and $(C_1,C_2) \stackrel{d}{=} (C_2,C_1)$.
We would like to stress at this point that we are not making any smoothness or specific dependence assumptions on the survival times.
It will just be required that 
the distribution of $\wh p$ is not degenerate.

Now, because all required Donsker properties on $\mc F$ and $\tilde {\mc F}$ obviously hold and $\phi$ is Hadamard-differentiable, (conditional) central limit theorems immediately apply if the condition in \eqref{eq:cond_mw} is met.
In particular, for independent $G_1, \dots, G_n$ with a uniform distribution on $\mc G^{exch.}$, 
it follows that the randomization empirical process
$\RPn$ based on $G_i(X_{i1},X_{i2},\delta_{i1},\delta_{i2}), i=1,\dots, n$, is asymptotically Gaussian:
as $n \to \infty$ and
conditionally on $X_{i1},X_{i2},\delta_{i1},\delta_{i2}, i=1,\dots, n$, the process
$ \sqrt{n}(\RPn - \Pn^{\mc G}) $ converges weakly in outer probability to a Gaussian process specified in Theorem~\ref{thm:main}.
Consequently, Theorem~\ref{thm:delta_meth} yields for the randomized Mann-Whitney effect
that $\sqrt{n}(\tilde p_n - \frac12) = \sqrt{n} (\phi(\RPn) - \phi(\Pn^{\mc G})) $ is asymptotically normal with some variance
$\tilde \sigma^2 \in (0,\infty)$, trivial cases ($\tilde \sigma^2=0$) excluded.

The influence function corresponding to the Mann-Whitney functional $\phi$ and consistent variance estimators $\wh \sigma^2_{\phi,\Pn}, \tilde \sigma^2_{\phi,\RPn}$ derived from these influence functions can be found in Appendix~\ref{app:MW}.
Finally, a randomization version of $T_n(p) = {\sqrt{n}(\wh p_n - p)}/{\wh \sigma_{\phi, \Pn}}$ is
$\tilde T_n = {\sqrt{n}(\tilde p_n - \frac12)}/{\tilde \sigma_{\phi,\RPn}}$.
Denote the conditional $(1-\alpha)$-quantile of $|\tilde T_n|$ by $\tilde c_n(\alpha)$, $\alpha \in (0,1)$, and the corresponding randomization probability by $\tilde \gamma_n(\alpha)$.
We obtain the following theorem about the resulting randomization hypothesis test:
\begin{cor}
 \label{thm:MW}
 Excluding trivial cases, we have for $n \rightarrow \infty$ under $H \cup K$ that, as $n \rightarrow \infty$,
 $\tilde c_n(\alpha)$ converges to $\Phi^{-1}(1-\alpha/2)$ in outer probability. Furthermore, the test
 $$ \Psi_n = 1\{ |T_n(.5)| > \tilde c_n(\alpha)\} + \tilde \gamma_n(\alpha) 1\{ |T_n(.5)| = \tilde c_n(\alpha)\} $$
 satisfies $ E(\Psi_n) \to 1_{K} + \alpha 1_{H}$ as $n \to \infty$.
\end{cor}
For the test in Corollary~\ref{thm:MW} one still needs to specify the value of the (randomized) test statistic if there was a division by zero due to a very unfavorable censoring pattern.
This could only happen for extremely small sample sizes in combination with particularly strong censoring rates.
It seems most natural to set the (randomized) test statistic to zero in such a case because nothing really can be concluded then.
Still, excluding trivial cases, the test $\Psi_n$ is finitely exact under exchangeability, that is, $H_0^{exch.} \subset H : \{p = \frac12\}$.

\begin{rem}
 \tcb{It is easy to see that even under $H_0^{exch.}$ the limiting Gaussian processes $\mb G$ and $\tilde {\mb G}$ do not share the same distribution.
 Denote by $H_j(t) = P(C_{ij}>t)$ is censoring survival function, $j=1,2$.
 For the indexing functions $f = 1_{(0,s]}(p_1)\cdot p_3$ and $h = 1_{(0,t]}(p_1)\cdot p_3$ it is straightforward to compute the following covariances:
 \begin{align*}
  \P (fh) - \P f \P h = & - \int_0^{\min(s,t)} H_1(u-) \d S_1(u)  \cdot \Big(1 + \int_0^{\max(s,t)} H_1(u-) \d S_1(u)\Big)
  \\
  \P(Q_\cdot (fh) - Q_\cdot f Q_\cdot h) \stackrel{H_0^{exch.}}{=} &
  - \frac12 \int_0^{\min(s,t)} H_1(u-) \d S_1(u)  - \frac12 P(X_1 \leq s, X_2\leq t, \delta_1=1, \delta_2=1).
 \end{align*}
}
\end{rem}

\begin{rem}
\label{rem:MWuneq}
\tcb{Due to the advantageous structure of the Mann-Whitney effect $p$, the test $\Psi_n$ from Corollary~\ref{thm:MW} can be extended to the case that some data are missing completely at random:
 suppose that the available observations correspond to $n_p$ paired individuals, i.e.\  $(X_{i1}, X_{i2}, \delta_{i1}, \delta_{i2}), i=1, \dots, n_p$, and, independent thereof, unpaired and independent individuals with observations, say, $(\tilde X_{ij}, \tilde \delta_{ij}), i=1, \dots, n_j, $ where $j=1,2$ again specifies the treatment group. We assume that $(X_{ij}, \delta_{ij})\stackrel{d}{=}(\tilde X_{ij}, \tilde \delta_{ij})$ for all $i$ and $j$. This can occur in matched pairs studies if for some test subjects no suitable match could be found or if single measurements are flawed.
 In this case, $p$ is still estimable with the same Kaplan-Meier-based approach, where the Kaplan-Meier estimators $\wh S_{j,n}$ are now based on the obvious $n_p + n_j$ data points, $j=1,2$.
 Thus, the presence of unpaired observations results in Kaplan-Meier estimators with an altered dependence structure.
 In Appendix~\ref{sec:MW_uneq} a variance estimator is developed that works for both cases, equal and unequal sample sizes.
 Additionally, Appendix~\ref{sec:MW_uneq_simu} contains additional simulations results based on unequal sample sizes and different marginal distributions: exponential and Gompertz.}
\end{rem}

\subsection{Randomizing Aalen-Johansen estimators for cumulative incidence functions}
\label{sec:cr2}

Usually, in competing risks survival situations, observability of event times and types are hindered due to independent right-censoring.
Hence, the observable data can be modelled as i.i.d.\ random vectors
$(Z_i, \delta_i \varepsilon_i, \delta_i) = (\min(T_i, C_i),$ $1\{T_i \leq C_i\} \varepsilon_i, 1\{T_i \leq C_i\}), \ i=1,\dots, n $, where
$T_1, \dots, T_n$ are the survival times and $\varepsilon_1, \dots, \varepsilon_n$ the event types which are only observable if the corresponding censoring indicators  $\delta_i=1\{T_i \leq C_i\}$ are equal to one.
The censoring times $C_i$ are assumed to be independent of the $(T_i, \varepsilon_i)$.
Estimation of the cumulative incidence functions $F_j(\tau)$ is commonly done with the Aalen-Johansen estimators,
$$ \widehat F_j(\tau) = \int_0^\tau \widehat S(u-) \d \widehat A_j (u) , \ j=1,2, $$
where the integrand is the left-continuous version of the Kaplan-Meier estimator for the overall survival probability 
and $$\widehat A_j(t) = \sum_{i:Z_i \leq t} \delta_i 1\{ \varepsilon_i =j \} / \sum_\ell 1\{ Z_\ell \geq Z_i\}$$ is the so-called Nelson-Aalen estimator for the cumulative hazard function of risk $j$.
We refer to \cite{aalen78} for the Aalen-Johansen estimator in more general multi-state Markov models and to \cite{lin97} for the above-stated form in competing risks situations.

The following idea will be used to find a suitable randomization method for a test for $H: F_1(\tau) \geq F_2(\tau)$ versus $K: F_1(\tau) < F_2(\tau)$ that is finitely exact under the boundary hypothesis $H_0 : F_1(\tau) = F_2(\tau)$:
for each individuum with an observed event type, i.e.\ when no censoring occurred, the event indicator could be randomized because under $H_0$ both risks are equally likely to happen until time~$\tau$.
Such a randomization scheme can be realized as follows.
Denote by $\mb Z / 2\mb Z$ the cyclic group with two elements, 0 and 1, then a randomization group that acts on the data is given by
$\mc G^{bin.} = \{ (p_1, p_3 \cdot (g+1) , p_3) : g \in \mb Z / 2\mb Z \}.$
Here the $p_j$ again denote the canonical projections.
Hence, for $z \in [0,\tau]$, the orbit of an observation of the form $(z,0,0)$ is just $\{(z,0,0)\}$
whereas the orbits of $(z,1,1)$ and $(z,2,1)$ are both equal to $\{(z,1,1),(z,2,1)\}$.
As a test statistic, one could choose for example $\widehat F_1(\tau) - \widehat F_2(\tau) $ or $\widehat F_1(\tau) / \widehat F_2(\tau) $ so that the null hypothesis is rejected for relatively small values of the statistic.
Because the randomization test will be finitely exact under $H_0$ and also finitely keep the level under $H$, there is no need for a studentization of the statistic.

\begin{rem}
The randomization based on $\mc G^{bin.}$ also has a connection to some kind of multiplier resampling scheme:
the randomized data could equally well be described as $G_i ( \b X_i)  = \Big(\begin{smallmatrix} 1 & & \\ & D_i+1 & \\ & & 1 \end{smallmatrix}\Big) \cdot \b X_i$, where $D_1, \dots, D_n$ are i.i.d.\ Bernoulli-distributed with parameter $p=.5$.
This presentation also makes obvious that other hypotheses can be tested similarly:
suppose for instance that a conventional drug results in the fraction $F_1(\tau) / F_2(\tau) = q \in (0,\infty)$.
For a new drug, that does not worsen the overall survival chances, one might additionally want to test whether the type-1-event probability is reduced in relation to the type-2-event probability, i.e.\
$H^{(q)}: F_1(\tau) / F_2(\tau) \geq q$ versus $K^{(q)}: F_1(\tau) / F_2(\tau) < q$.
A test for this could be based on a similar resampling scheme as the one above, 
except that the $D_i$ now should be Bernoulli distributed with parameter $q/(q+1) \stackrel{H_0^{(q)}}{=} F_1(\tau)/(F_1(\tau) + F_2(\tau))$ which reflects the situation under the boundary hypothesis $H_0^{(q)}: F_1(\tau) / F_2(\tau) = q$.
This approach would still give rise to a test for $H^{(q)}$ versus $K^{(q)}$ that is finitely exact under $H_0^{(q)}$ because the independent censoring assumption ensures that the resampled data exactly reflect the situation under $H_0^{(q)}$.
However, it is not possible anymore to express the randomization scheme with the help of an algebraic group-based randomization approach.
\end{rem}

\section{Simulation study: tests about the Mann-Whitney effect}
\label{sec:simus}

Below we describe a simulation study concerning the testing of hypotheses about the Mann-Whitney effect, i.e.\ $H: p=.5$ against $K: p \neq .5$,  under different data dependence structures, marginal distributions, and censoring intensities.
In addition, Appendix~\ref{sec:simus_add} contains a simulation study on the reliability of tests for correlation based on the empirical correlation coefficient.

In particular, we considered the following \emph{copulae}:
  Clayton with parameter -.6, i.e.\ negatively correlated data;
  Gumbel-Hougaard with parameter 5, i.e.\ positively correlated data;
  independence.
 
 \emph{Marginal distributions}: 
  equal exponential distributions with rate 2;
  an exponential distribution with rate 2 and a 50/50 mixture  of exponential distributions with parameters 3 and 1.316; the latter parameter is such that $H$ is approximately true.

 Three \emph{right-censoring intensities}, based on the minimums of $\tau=1$ and uniformly distributed random variables with minimum parameter 0 and maximum parameters
 2.7, i.e.\ about 24.6\%/26.1\% censorings,
 1.6, i.e.\ about 31.9\%/33.1\% censorings,
 and 1.1, i.e.\ about 40.6\%/41.2\% censorings (exponential / mixture survival distribution).
These censoring rates have been found via simulation of 100,000 individuals.

The sample sizes varied from $n=25$ to $150$ with increments of $25$.
We chose the significance levels $\alpha=1\%, 5\%, 10\%$.
We used the following methods to conduct the tests:
 a randomization method based on randomly interchanging the sample group correspondence within each pair, i.e.\ the randomization group $\mc G^{exch.}$; this corresponds to finite exactness of the resulting test under exchangeability, i.e.\ an exchangeable copula and equal marginals which is satisfied in all simulation configurations involving equal marginal exponential survival distributions;
 Efron's bootstrap for survival data \citep{Efro:cens:1981};
  quantiles of the standard normal distribution.
For the latter two methods no finite exactness is achieved in any simulation setting.
\tcb{Each test was simulated 5,000 times and was based on $B=$2,000 randomization/bootstrap iterations.}
\begin{figure}
 \includegraphics[width=0.34\textwidth]{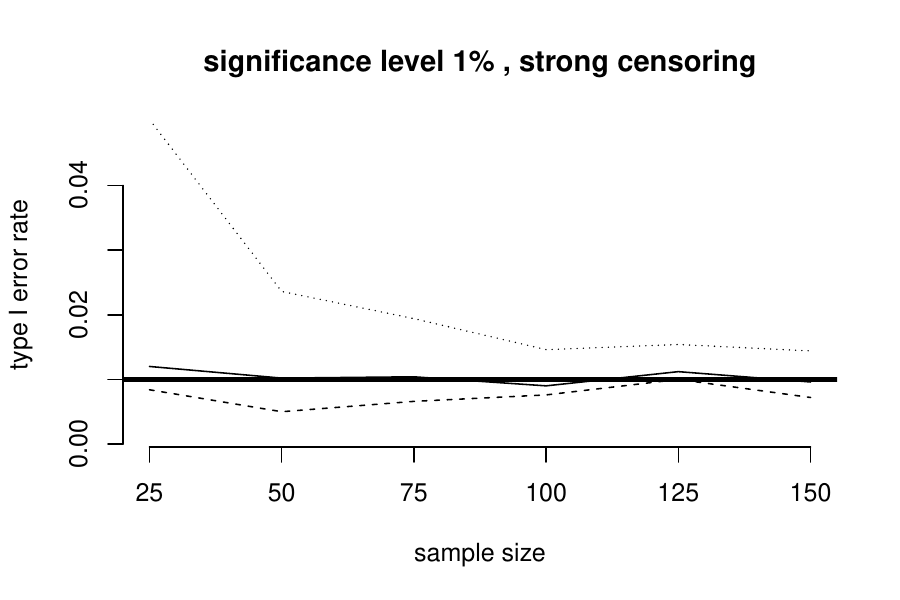}
 \hspace{-0.5cm}
 \includegraphics[width=0.34\textwidth]{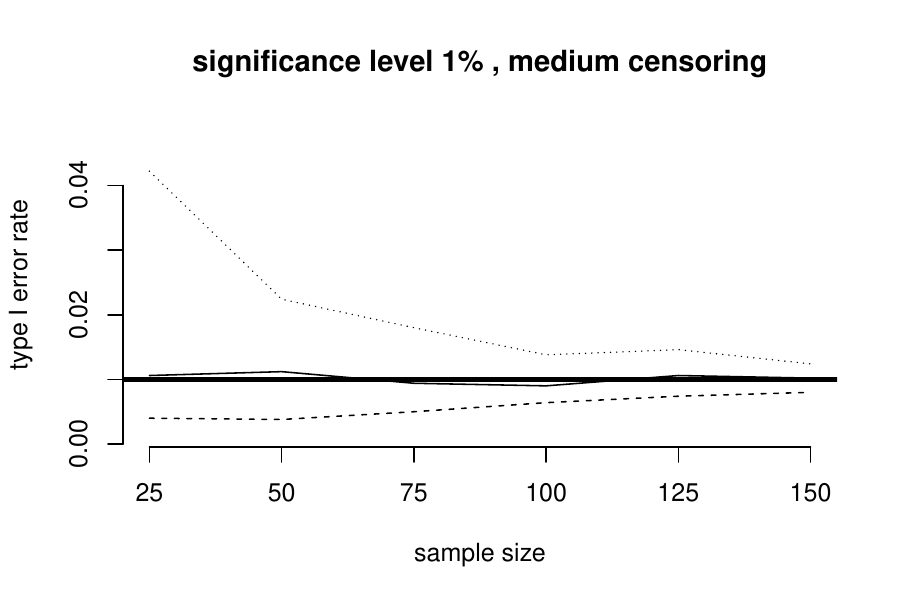}
 \hspace{-0.5cm}
 \includegraphics[width=0.34\textwidth]{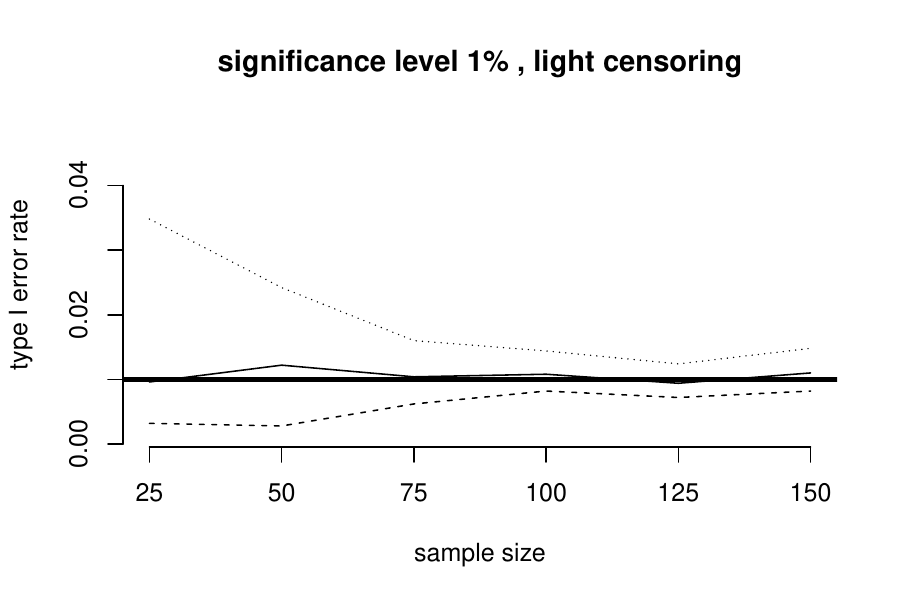} \\[-0.4cm]
 \includegraphics[width=0.34\textwidth]{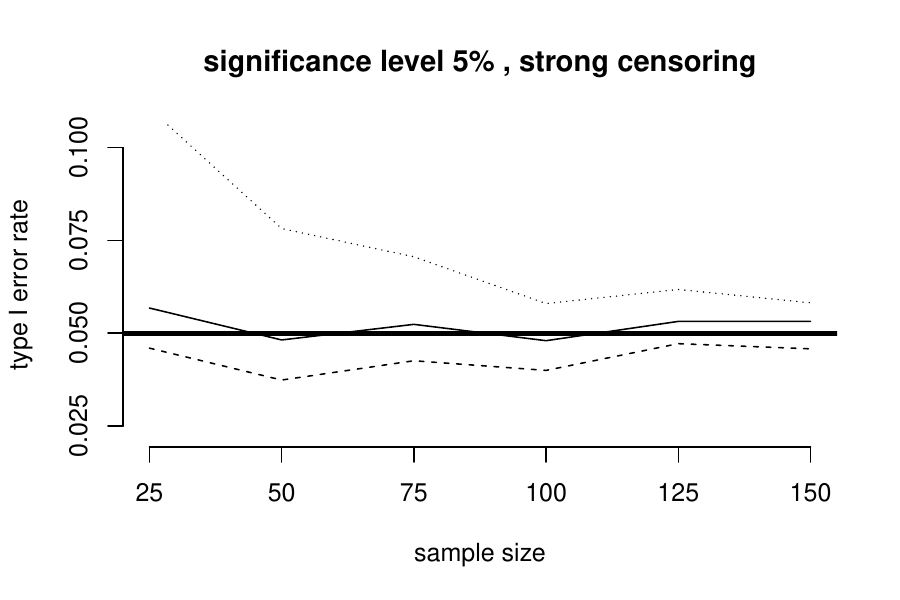}
 \hspace{-0.5cm}
 \includegraphics[width=0.34\textwidth]{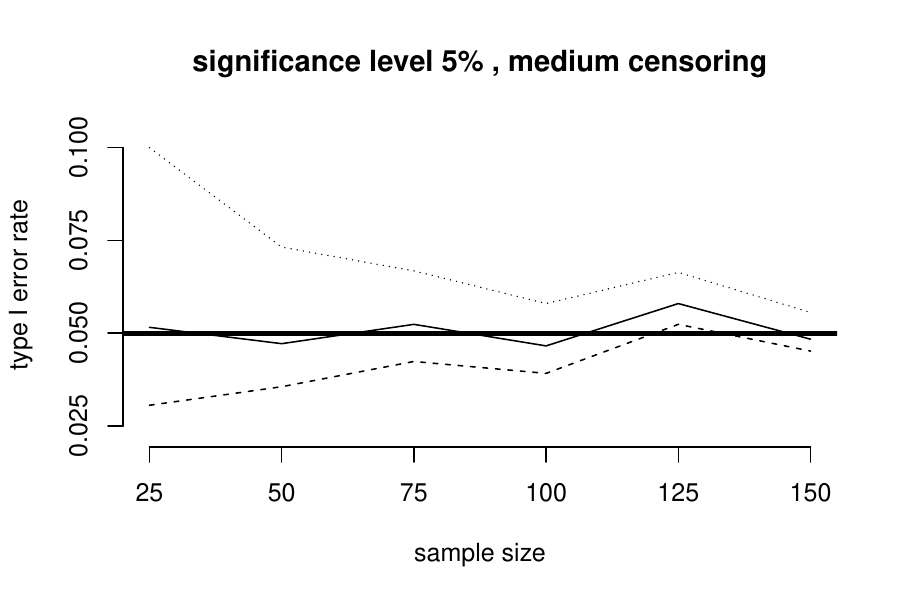}
 \hspace{-0.5cm}
 \includegraphics[width=0.34\textwidth]{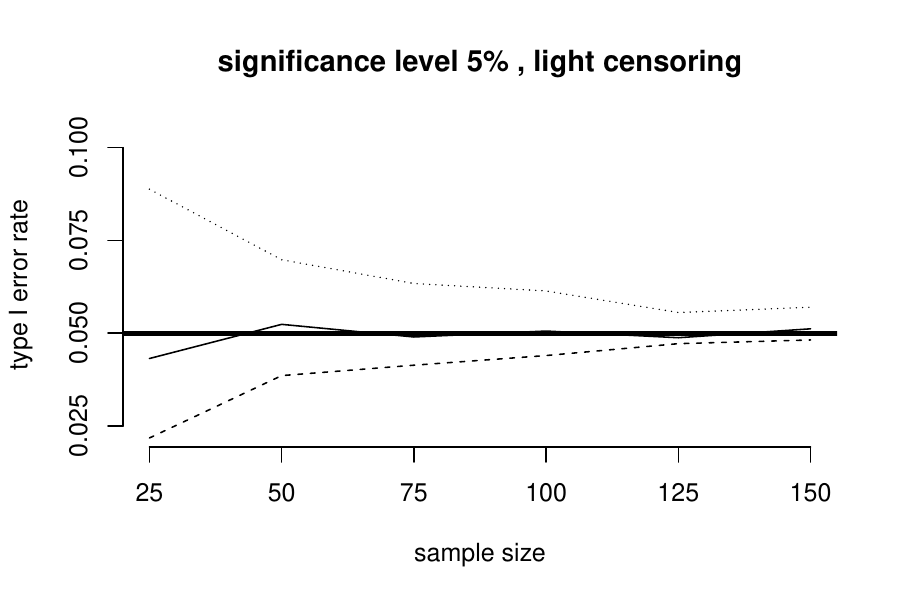} \\[-0.4cm]
 \includegraphics[width=0.34\textwidth]{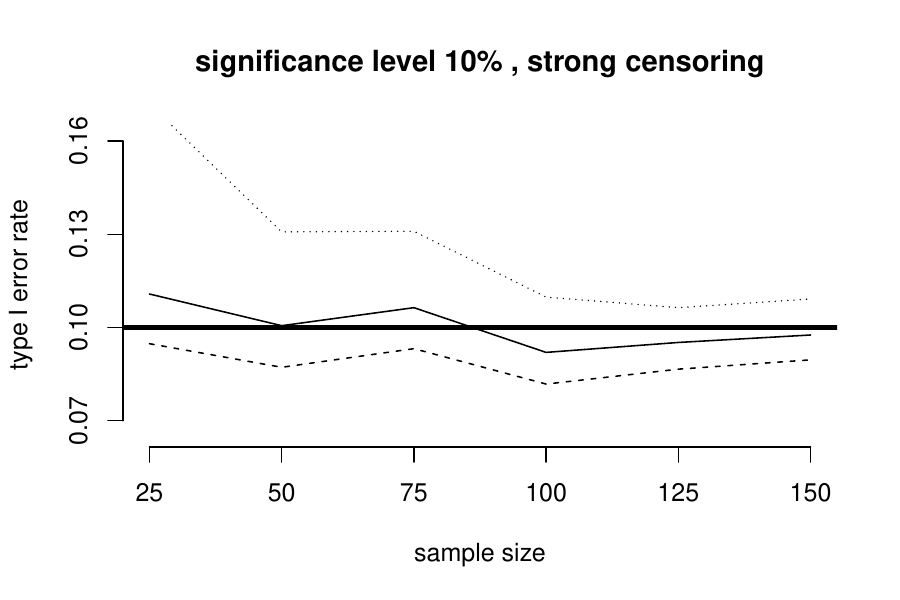}
 \hspace{-0.5cm}
 \includegraphics[width=0.34\textwidth]{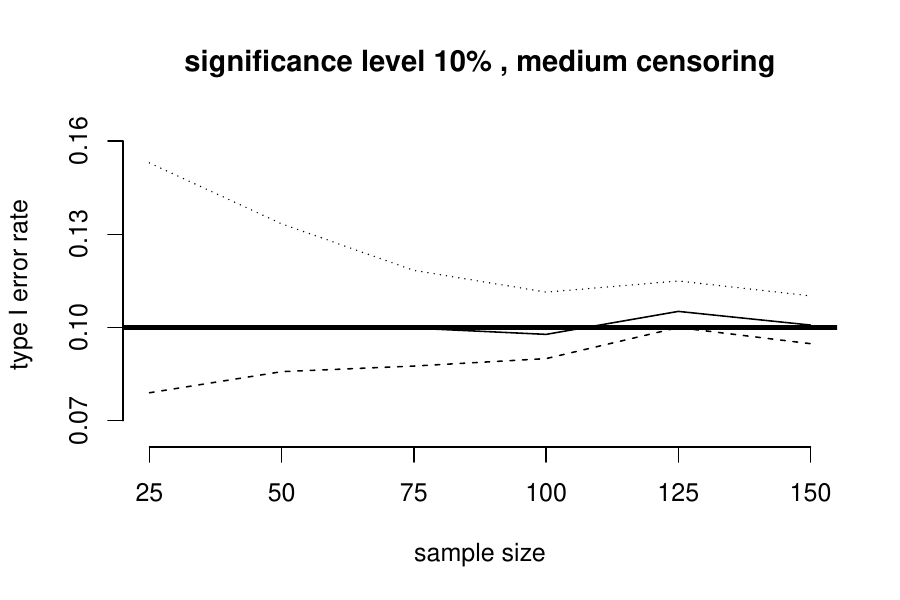}
 \hspace{-0.5cm}
 \includegraphics[width=0.34\textwidth]{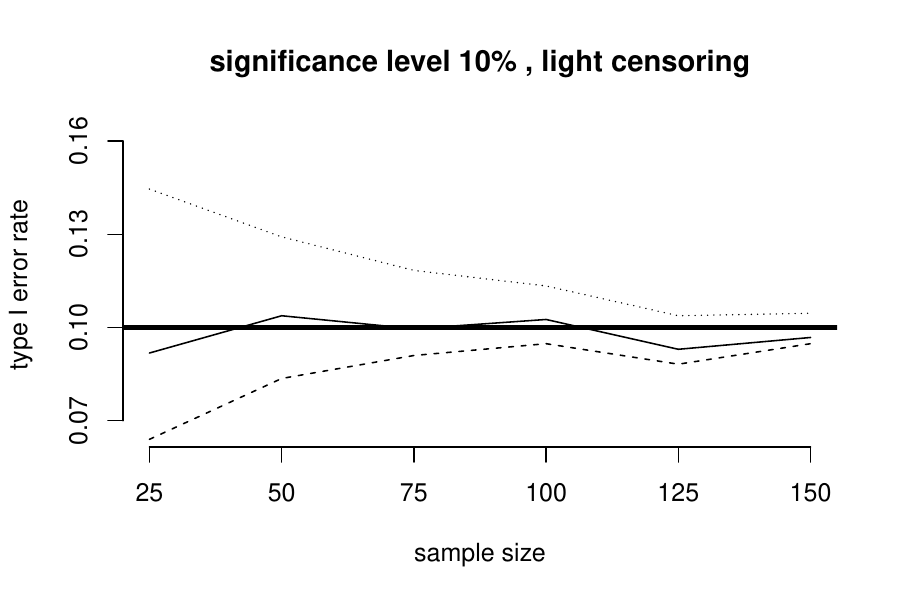} \\[-0.4cm]
 \includegraphics[width=0.34\textwidth]{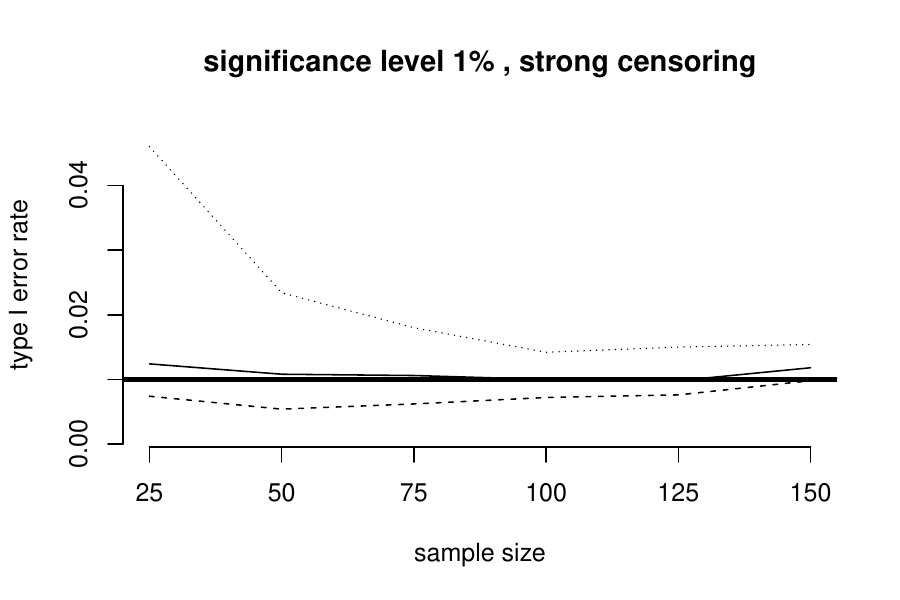}
 \hspace{-0.5cm}
 \includegraphics[width=0.34\textwidth]{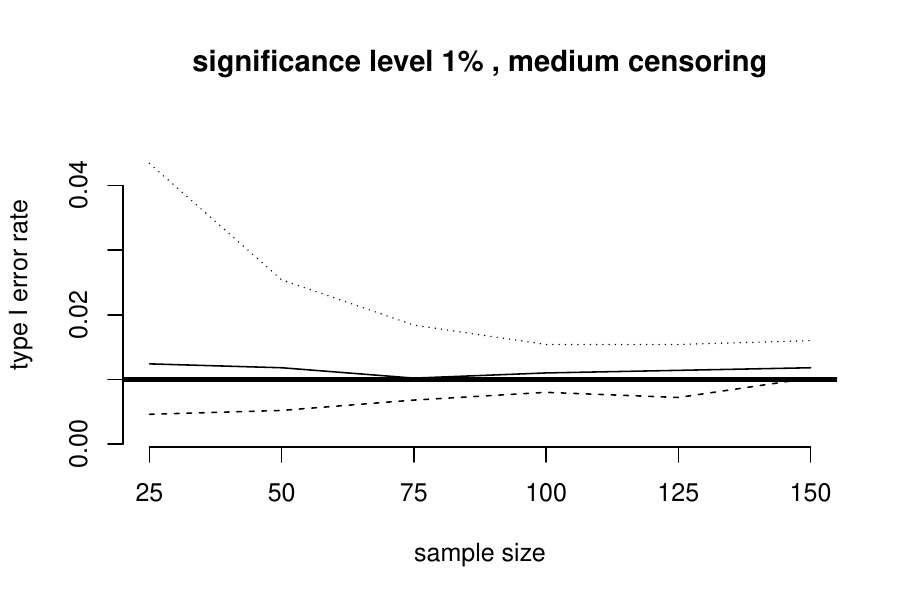}
 \hspace{-0.5cm}
 \includegraphics[width=0.34\textwidth]{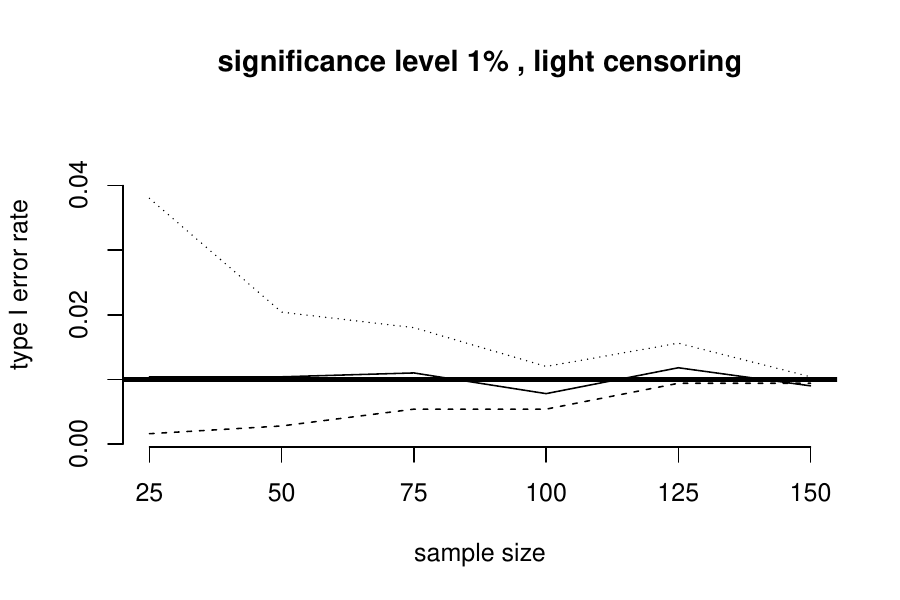} \\[-0.4cm]
 \includegraphics[width=0.34\textwidth]{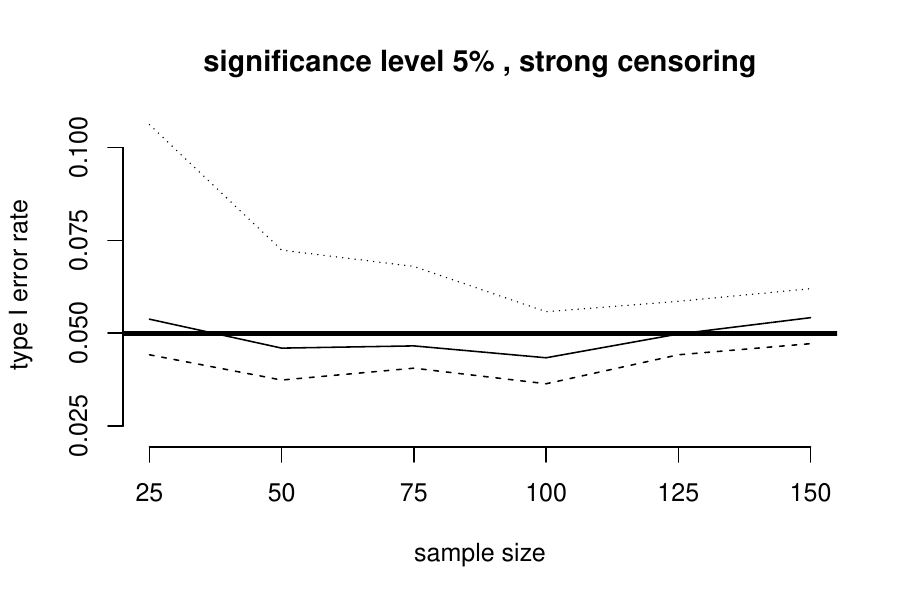}
 \hspace{-0.5cm}
 \includegraphics[width=0.34\textwidth]{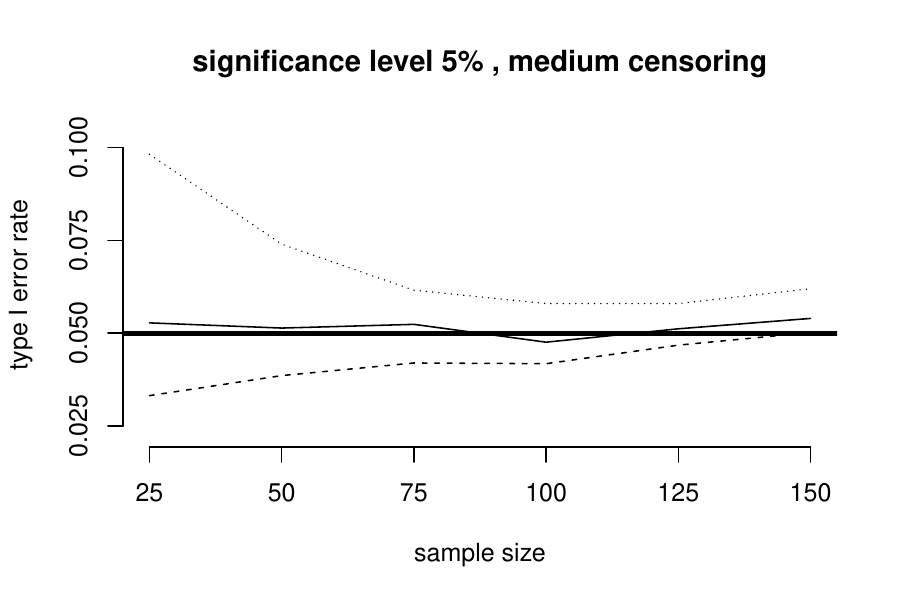}
 \hspace{-0.5cm}
 \includegraphics[width=0.34\textwidth]{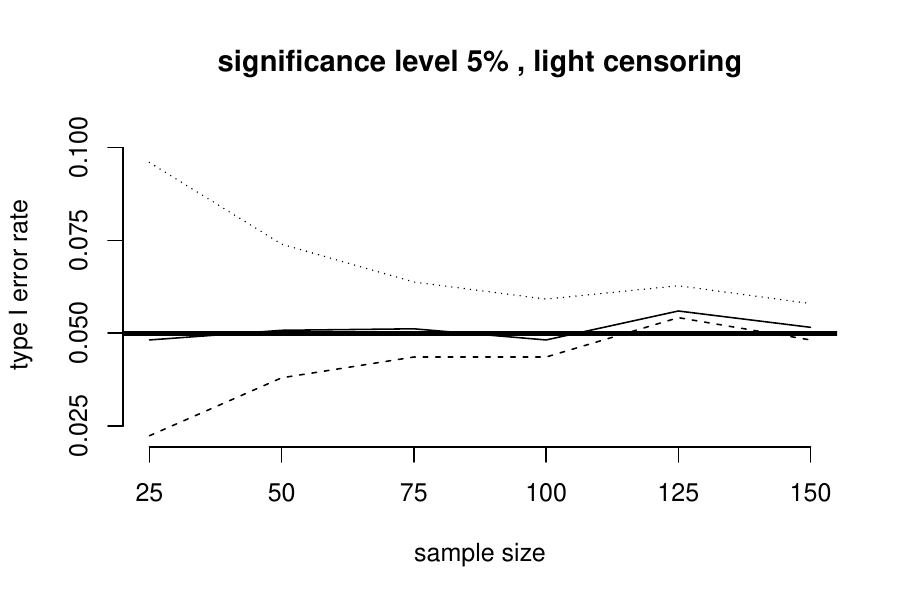} \\[-0.4cm]
 \includegraphics[width=0.34\textwidth]{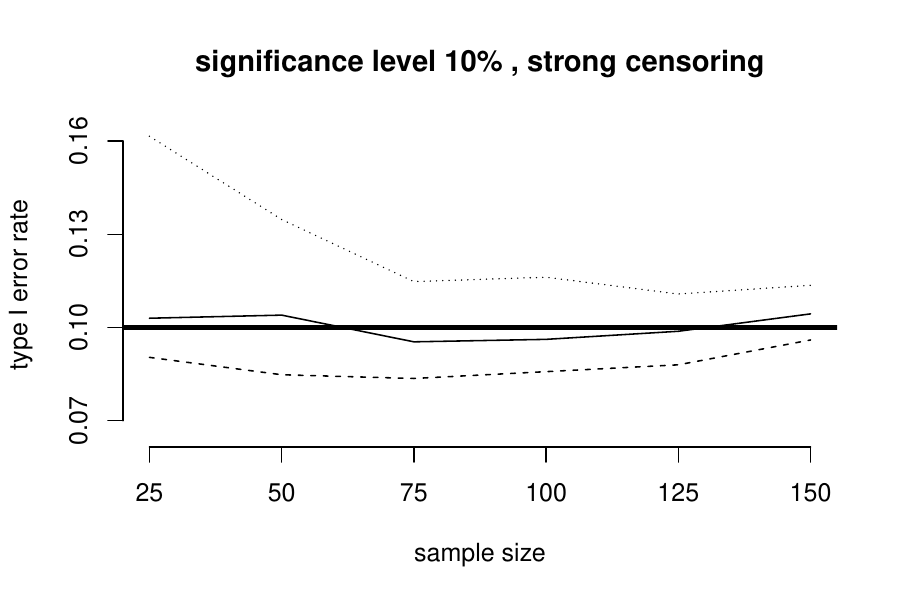}
 \hspace{-0.5cm}
 \includegraphics[width=0.34\textwidth]{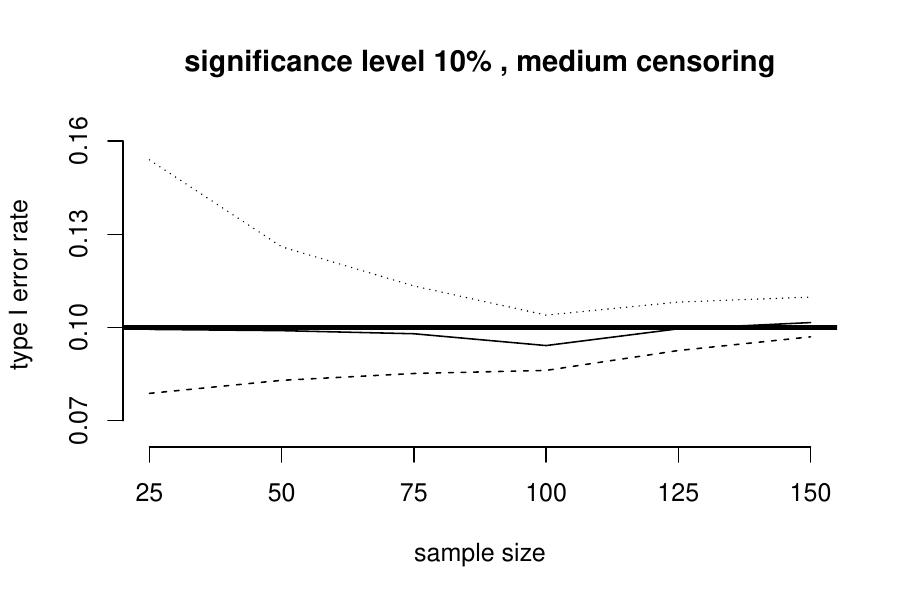}
 \hspace{-0.5cm}
 \includegraphics[width=0.34\textwidth]{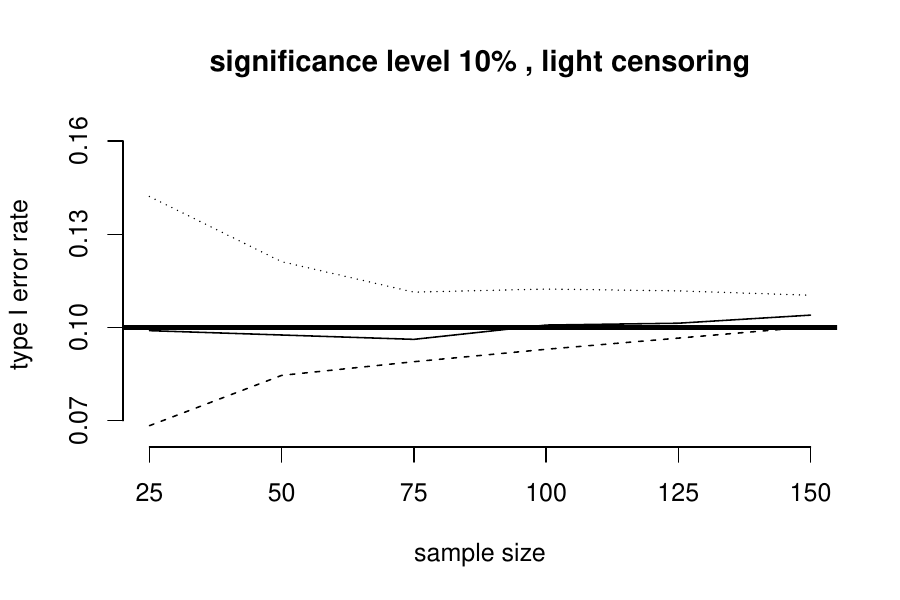}
 \vspace{-0.7cm}
 \caption{Simulated type I error rates of the Mann-Whitney-type tests with the Clayton copula underlying the data under strong (left), medium (middle), and light censoring (right); equal (upper half) and unequal marginal survival distributions (lower half); based on randomization (---), bootstrap (- -), normal quantiles ($\cdots$). 
 The nominal significance level is printed in bold.}
 \label{fig:mw2}
\end{figure}

The results are illustrated in Figure~\ref{fig:mw2} for the Clayton copula and Figures~\ref{fig:mw1} and \ref{fig:mw3} in Appendix~\ref{sec:simus_add} for the Gumbel-Hougaard copula and the independence case. 
Comparing the three considered methods for finding critical values, we find the same overarching picture in all simulation configurations: the tests based on the standard normal quantiles are (very) liberal and the bootstrap-based tests are rather conservative. \tcb{One notable exception is the case of the Gumbel-Hougaard copula with medium or light censoring rates; all tests produce reliable results here; cf.\ Figure~\ref{fig:mw1} in the appendix.}

In contrast \tcb{to the general impression of the asymptotic and bootstrap tests}, the randomization-based tests achieve excellent rejection probabilities, i.e.\ they are very close to the nominal significance level in all considered set-ups.
Comparing the results for the different censoring intensities, we do not see a big difference, except for the Gumbel-Hougaard case in the appendix; 
there, the tests apparently get more reliable with stronger censoring rates.

To \tcb{illustrate} the power of the tests, we generated the survival times in the first group according to
\tcb{$T_{i1} \sim Exp(2)$ and in the second group according to $T_{i2} \sim .5 \cdot Exp(3) + .5 \cdot Exp(\lambda)$.
Here, the rate $\lambda$ in the mixture distribution was chosen such that the Mann-Whitney effect took the values $p=.55, .6, .65,.7,.75,.8$ while the terminal time is $\tau =1 $, i.e.\ $\lambda \approx 1.176, 1.045, .923, .809, .701, .600$. 
We considered the sample sizes $n=50, 100, 150$.
\begin{figure}[ht]
 \includegraphics[width=0.7\textwidth]{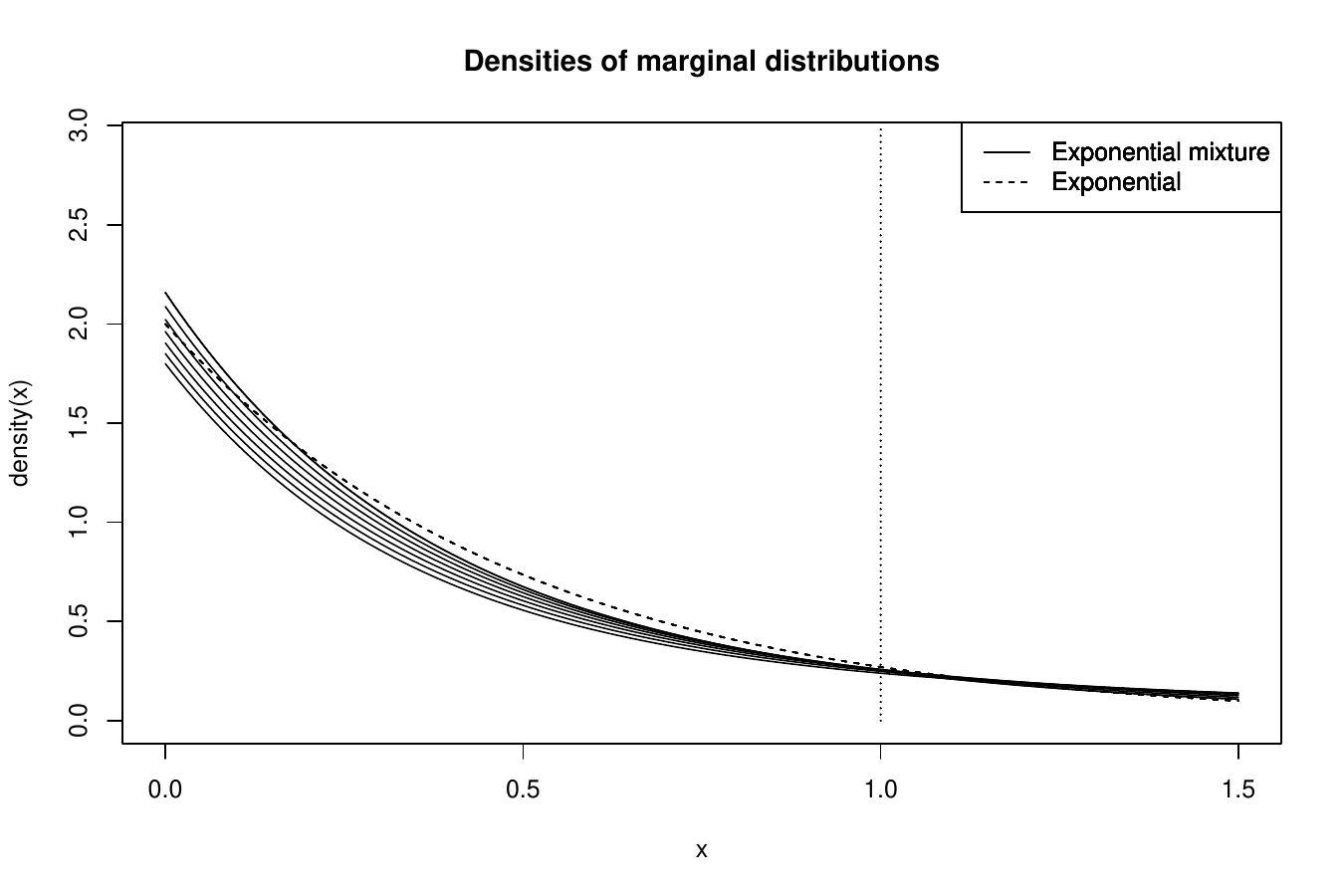}
 \caption{Marginal distributions leading to $p=.5, .55, .6, .65,.7,.75,.8$. The terminal time is $\tau = 1$.}
\end{figure}
\begin{figure}[ht]
 \includegraphics[width=0.34\textwidth]{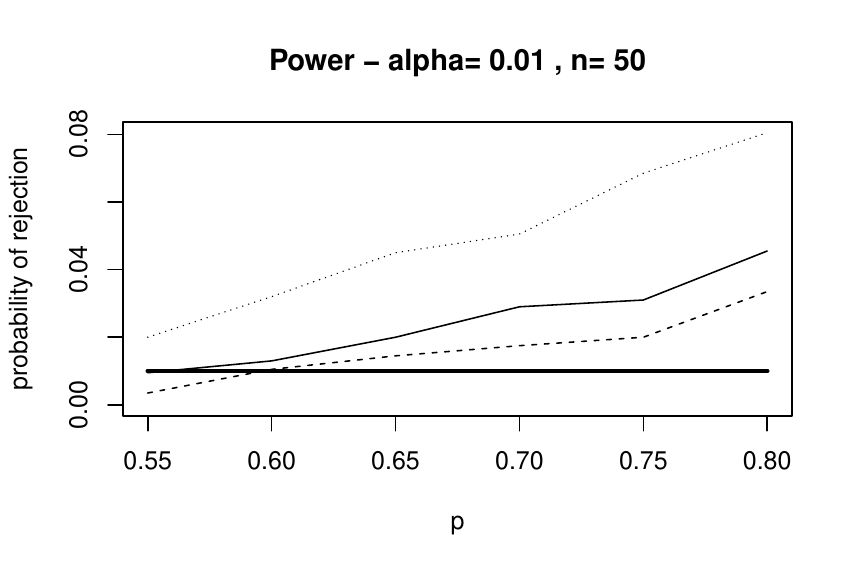}
 \hspace{-0.5cm}
 \includegraphics[width=0.34\textwidth]{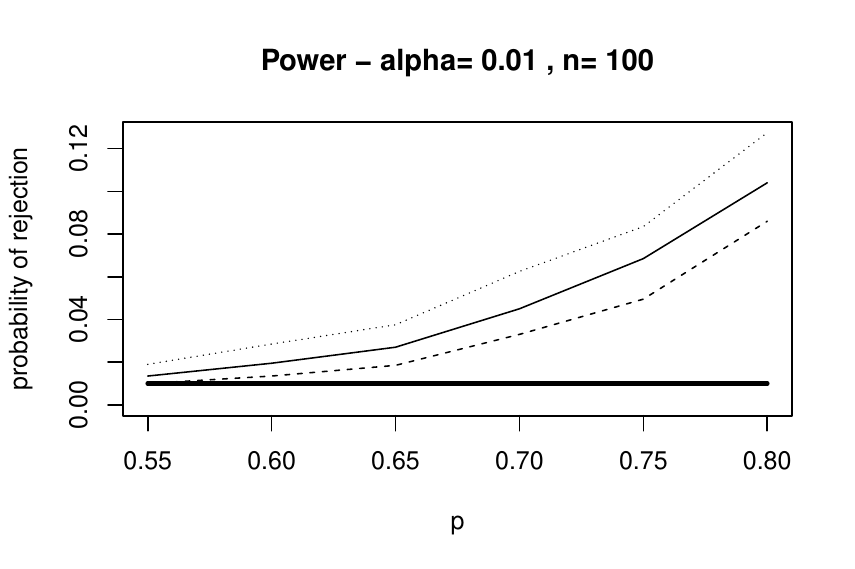}
 \hspace{-0.5cm}
 \includegraphics[width=0.34\textwidth]{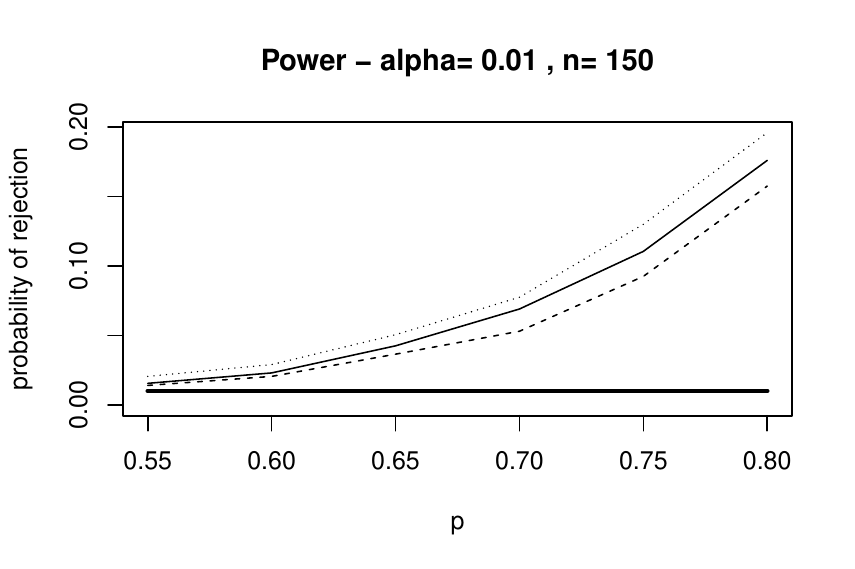}
  \\[-0.4cm]
   \includegraphics[width=0.34\textwidth]{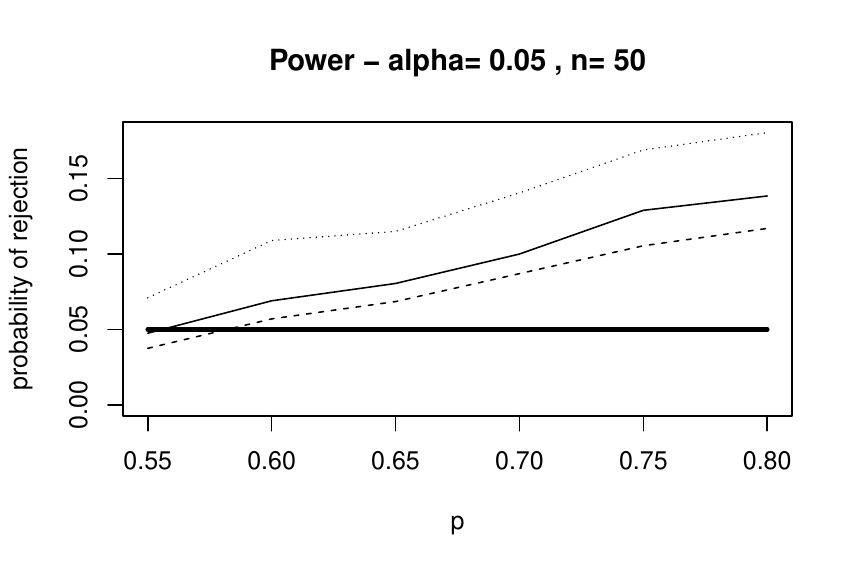}
 \hspace{-0.5cm}
 \includegraphics[width=0.34\textwidth]{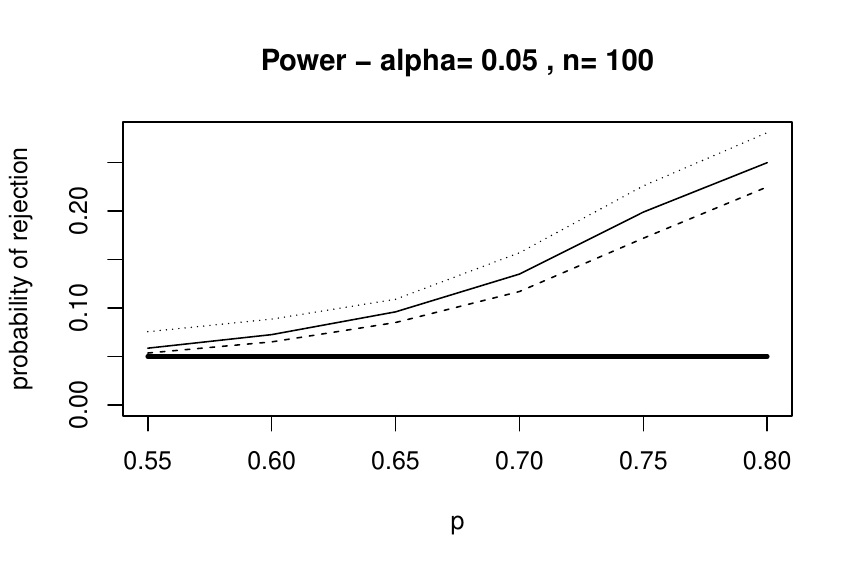}
 \hspace{-0.5cm}
 \includegraphics[width=0.34\textwidth]{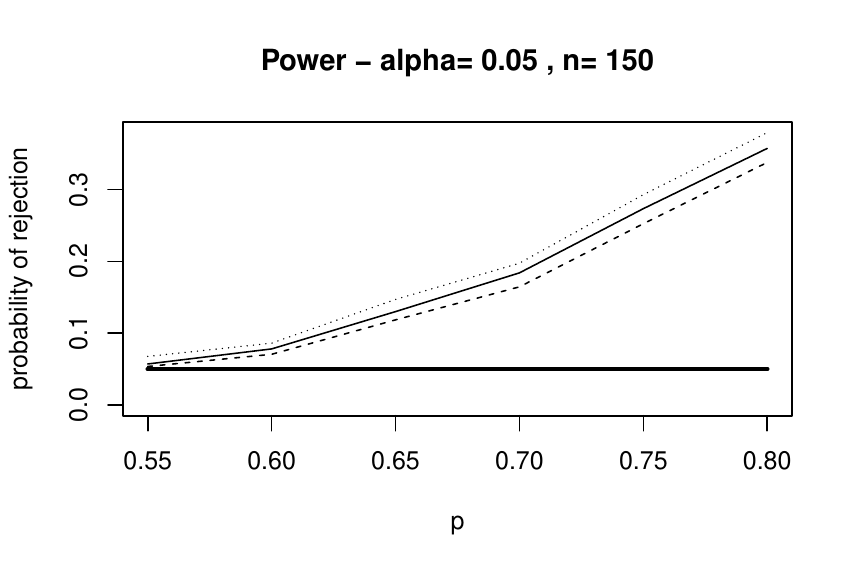}
  \\[-0.4cm] 
  \includegraphics[width=0.34\textwidth]{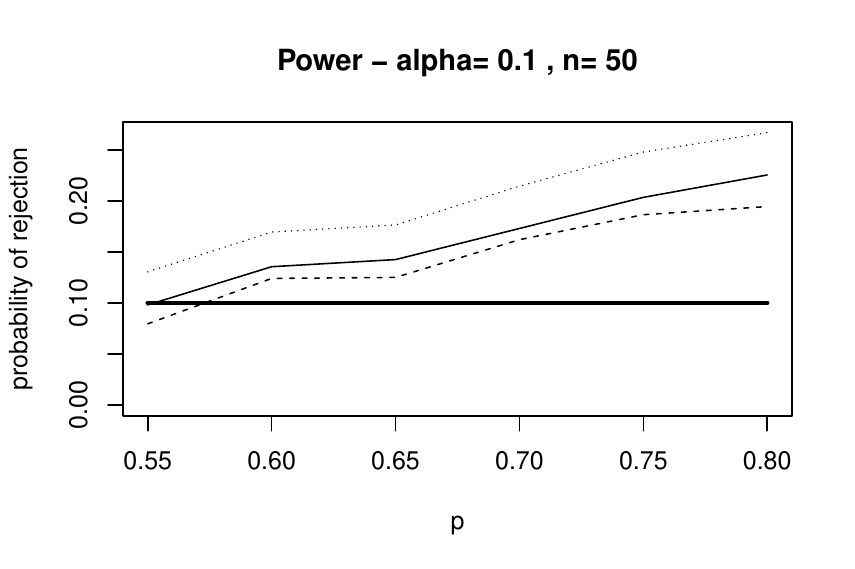}
 \hspace{-0.5cm}
 \includegraphics[width=0.34\textwidth]{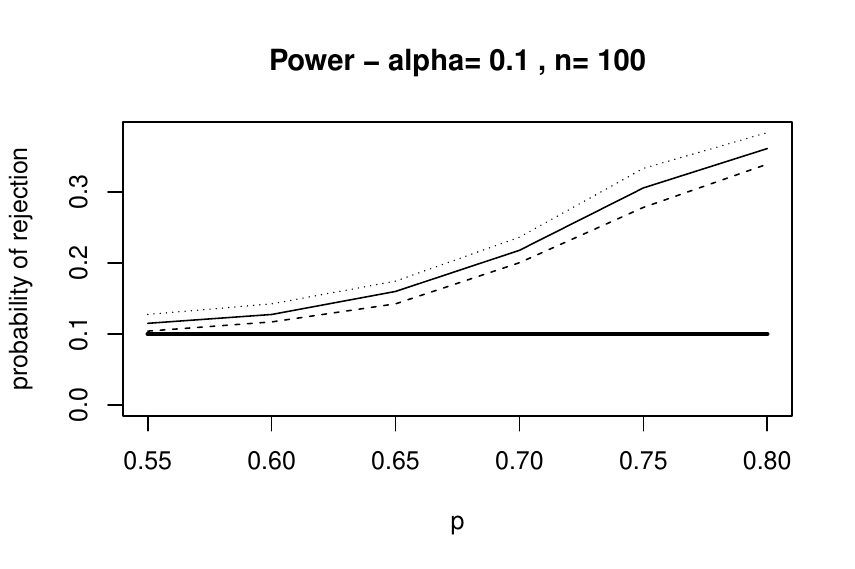}
 \hspace{-0.5cm}
 \includegraphics[width=0.34\textwidth]{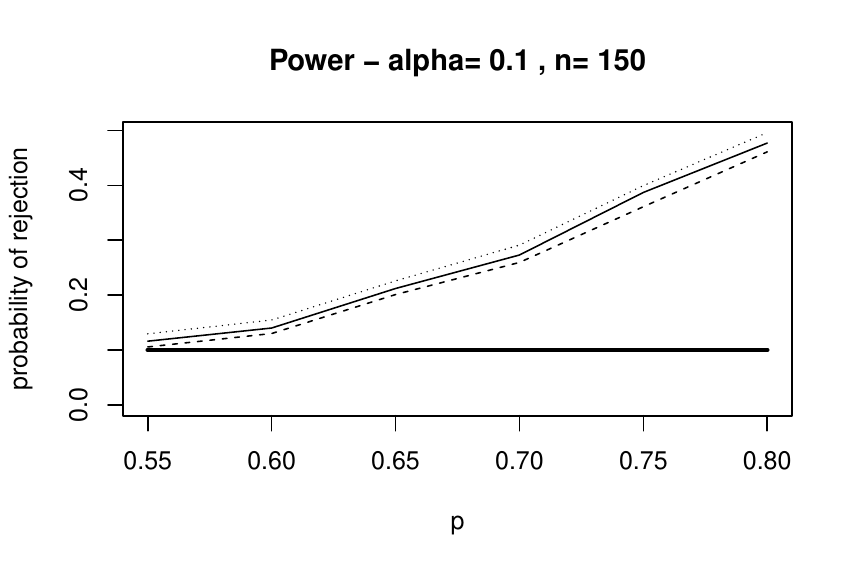}
  \caption{Simulated power of the Mann-Whitney-type tests with the Clayton copula underlying the data under strong censoring  and unequal marginal survival distributions; based on randomization (---), bootstrap (- -), normal quantiles ($\cdots$). 
 The nominal significance level is printed in bold.}
 \label{fig:power}
\end{figure}
The remaining simulation settings were equal to one set-up for simulating the rejection rates under the null: Clayton copula with strong censoring and significance levels $\alpha=.01, .05, .1$.
In contrast to the simulations above, we repeated the tests 2,000 times and chose $B=$1,000 iterations for the randomization and bootstrap-based critical values.}

\tcb{The simulated rejection rates as shown in Figure~\ref{fig:power} are not entirely surprising: 
the conservativeness of the bootstrap test and the anit-conservativeness of the asymptotic test (cf.~Figure~\ref{fig:mw2}) result in the relatively low and high rejection rates, respectively, as compared to those of the randomization-based test.
The power simulation results for the other scenarios are not displayed.}

\section{Data example}
\label{sec:data}

We are going to apply the Mann-Whitney test of Corollary~\ref{thm:MW} to data about patients suffering from diabetic retinopathy. 
The data are available from the \verb=timereg R=-package in the dataset \verb=diabetes=.
It contains information on $N=197$ patients for each of whom a randomly selected eye was treated by means of a laser photocoagulation,
the other eye was observed without a treatment. 
The recorded ``survival times'' are the times to blindness or censoring, whatever came first.
The data can be divided according to the age at onset of diabetes:
for easy reference, we will call these the ``juvenile'' ($n_1=114$) and ``adult'' ($n_2=83$) subgroups;
see \cite{huster89} for a more complete description of the study.

\begin{figure}[ht]
\centering
 \vspace{-0.35cm}
 \includegraphics[width=0.8\textwidth]{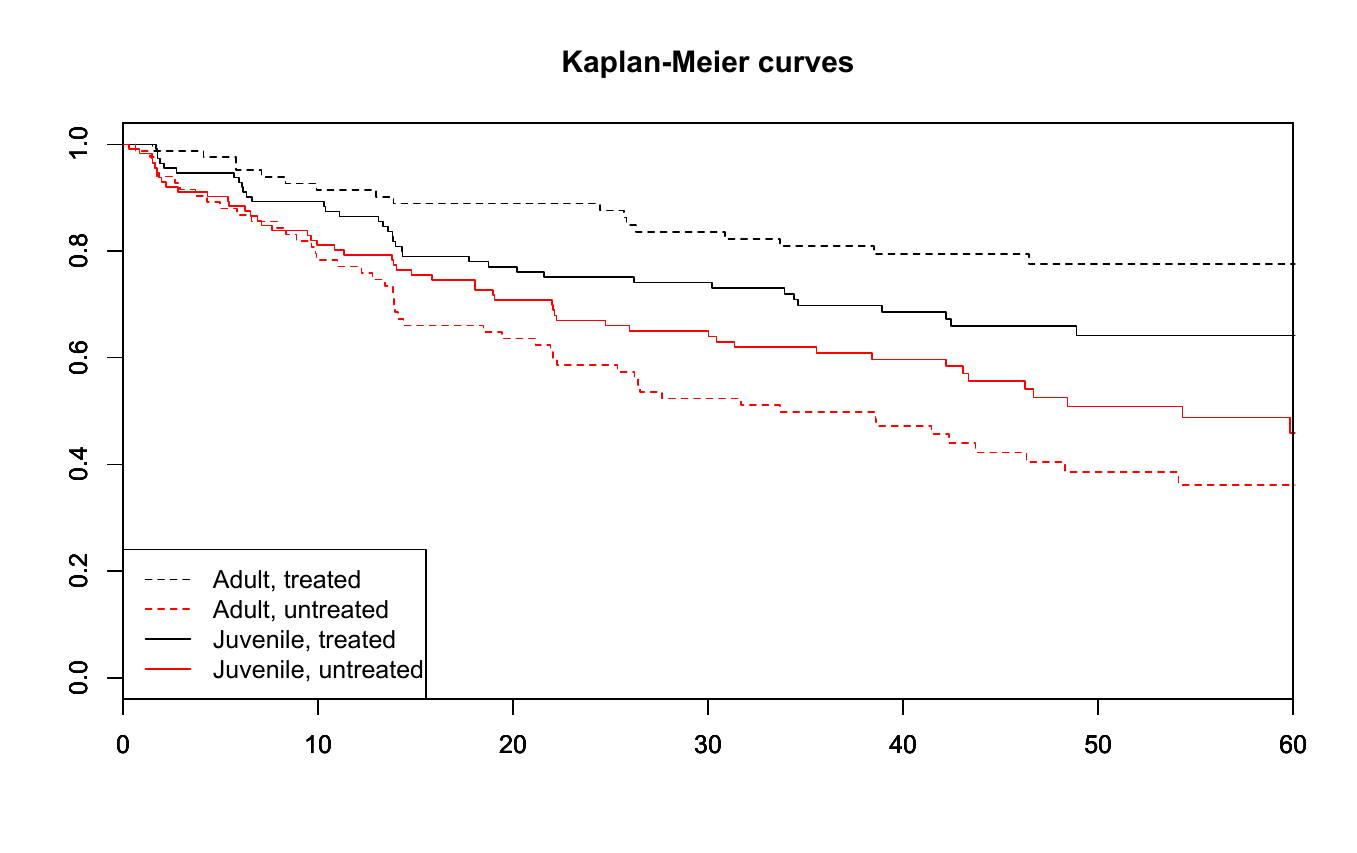}
 \vspace{-1cm}
 \caption{Kaplan-Meier estimates in all considered subgroups in the diabetes dataset.}
 \label{fig:kmes}
\end{figure}

The first research question of interest was whether the treatment was effective in delaying the onset of blindness.
Using parametric models and Wald tests, \cite{huster89} were able to verify this for both subgroups.
They also found a significant interaction effect between treatment and age at onset of diabetes; see Figure~\ref{fig:kmes} for an illustration of this by means of the group-specific nonparametric Kaplan-Meier estimators.
We again refer to their article for more details and additional statistical analyses.

In the following, we will check whether the two-sided randomization-based test developed in Corollary~\ref{thm:MW} arrives at the conclusion that the Mann-Whitney effects are different from .5.
Since the follow-up time of interest was five years, 
we have chosen $\tau = 60$ months as the terminal evaluation time.
This choice leads to the following censoring rates within the subgrous:
\\
 \phantom{X} $\bullet$ juvenile; \quad treated: 68.42\% (78 out of 114);   \quad untreated: 55.26\% (63 out of 114); \\
 \phantom{X} $\bullet$ adult; \quad treated: 79.52\% (66 out of 83); \quad  untreated: 40.96\% (34 out of 83). \\
Hence, we observe quite high right-censoring rates, in particular in the treatment subgroups.
In particular, it seems that the sub-hypothesis $H_0^{exch.}$ of exchangeability between both treatments is not true, neither for the juvenile nor for the adult subgroup.
The Mann-Whitney effect estimates are $.5805$ (juvenile) and $.7074$ (adult).
In words, a treated eye in the adult subgroup has an estimated chance of 70.74\% to evade blindness longer than untreated eyes in the same subgroup.
The effect was weaker in the juvenile group, but still the treatment is favored ($.5805 > .5$).

Applications of the tests based on the Mann-Whitney effect that were  considered in Section~\ref{sec:simus}, i.e.\ the tests based on randomizing the treatment, bootstrapping, and the asymptotic normal distribution,
yielded the following $p$-values for the juvenile group, where $B=$ 2,000 \tcb{randomization/bootstrap} iterations were chosen for the first two tests: 
randomization: \tcb{.0105}; \ bootstrap: \tcb{.011}; asymptotic: \tcb{.0118}.
In the adult group, all $p$-values are less than $.001$.
As a consequence, all tests reject the null hypothesis of no Mann-Whitney effect at the significance level $\alpha=5\%$, even in the juvenile group where the effect was not as large as in the adult group.
The following two-sided confidence intervals were obtained by inverting the two-sided hypothesis tests and provide more information on the effect sizes:

\begin{tabular}{c|ccc|ccc}
 subgroup & \multicolumn{3}{c|}{juvenile} & \multicolumn{3}{c}{adult} \\ 
confidence level & 90\% & 95\% & 99\% & 90\% & 95\% & 99\% \\ \hline
 randomization & $[.528, .633]$ & $[.517, .645]$ & $[.499, .662]$ & $[.650, .765]$ & $[.639, .775]$ & $[.619, .795]$ \\
 bootstrap & $[\tcb{.526}, .634]$ & $[.516, .645]$ & $[\tcb{.497,.664}]$ & $[\tcb{.649, .765}]$ & $[.639, \tcb{.775}]$ & $[\tcb{.616,.799}]$ \\
 asymptotic & $[.528, .633]$ & $[\tcb{.518,.643}]$ & $[.498, .663]$ & $[.652, .763]$ & $[\tcb{.641,.773}]$ & $[.621, .794]$ \\ \hline
\end{tabular}
\\[0.4cm]
All in all, we see that, for each subgroup and nominal confidence level, all three obtained confidence intervals are very similar.
We understand this as an indication that the asymptotic results are taking effect because, as was seen in the simulations of Section~\ref{sec:simus}, the asymptotic tests were quite liberal and the bootstrap tests rather conservative, at least for small sample sizes.
As a consequence, the above confidence intervals and test results seem trustworthy.

\section{Discussion and future research}
\label{sec:disc}

We developed an empirical process theory for randomization-based tests, i.e.\ a conditional weak convergence result and a functional delta-method for the randomization empirical process.
These, in combination with appropriate studentizations,
allowed the construction of asymptotically exact hypothesis tests that are also exact for finite samples under the sub-hypothesis of invariance under the randomization operation.
Future research will focus on the development and application of randomization-based tests in  multivariate testing problems in which the limit distributions of the test statistics might be non-normal.

In the analysis of the dataset about the laser treatment on eyes of diabetic patients we have come to solid conclusions,
without the need to make parametric model assumptions.
It would be interesting to extend the Mann-Whitney effect-based test to a multi-sample test for detecting
an interaction effect between the kind of diabetes and the treatment.
Another future paper will consider statistical inferences on an above-mentioned variant of the Mann-Whitney effect for paired survival data: $\check p = P(T_{11} > T_{12}) + .5 P(T_{11} = T_{12})$ which is a parameter related to a within-pair comparison. 
Even though, as argued in Section~\ref{ex:mw}, a utilization of the Mann-Whitney effect $p$ often seems more natural than the use of $\check p$,
there are situations in which $\check p$ could prove more useful.
For example, an estimate of the possible gain in the expected survival duration, $E(T_{11} - T_{12})$, is probably best accompanied with the related parameter $\check p$.
However, estimation of this parameter requires estimation of (part of) the bivariate survival function of $(T_{11}, T_{12})$. 
It should be noted that such estimation -- based on right-censored paired data -- involves much more complicated functionals than the one involved in the present paper; cf. \cite{gill95}.

One referee suggested to analyze randomization empirical processes that are based on transformations with a more general structure, $G_i(\b X_1, \dots, \b X_n)$ instead of $G_i(\b X_i)$.
However, the different families of limiting Gaussian processes, $\tilde {\mb G}$ (cf.\ Theorem~\ref{thm:main}) which is a $\P$-mixture of $Q$-Brownian motions, and classical Brownian motions resulting from the random permutation approach (see Section~3.7 in \citealt{vaart96:_weak}) already give a taste of the difficulty of handling these within a unified approach.
Nevertheless, it seems interesting to research additional conditions on the algebraic group $\mc G$ that would allow the development of such a theory.
However, this is beyond the scope of the present article.
In this regard, we again wish to point to Sections~\ref{sec:perm} and~\ref{sec:bs} in which combinations of random permutation with randomization and the connection of the bootstrap with the present randomization framework have already been discussed.
These form a first step towards more general randomization procedures in which the randomized observations are not necessarily i.i.d.\ and for which conditional weak convergence theorems might still hold.

\tcb{The present paper focused on the case that studentized test statistics asymptotically have a normal distribution under the null hypothesis.
In future articles the author plans to consider multivariate extensions of the present theory, for instance general randomization methods for tests based on quadratic forms.
If, however, the limit null distribution is not pivotal, i.e.\ if it depends on unknown parameters, constructing an asymptotically exact randomization test will be more cumbersome or even impossible:
in contrast to the bootstrap empirical process, the randomization empirical process has a fundamentally different limit distribution as was seen in Theorem~\ref{thm:main}.
If a test statistic cannot be made pivotal by means of a studentization, the distributional difference between the test statistic and its randomized version will usually persists.
In this context, we again refer to \cite{wu20} who discussed the role of studentization and the applicability of permutation tests based on quadratic forms.
}

\newpage
\appendix
\section*{Appendix}

\section{Further proofs}

\begin{proof}[Proof of Theorem~\ref{thm:main}]
 First, we give a proof of the conditional weak convergence of all finite-dimen\-sio\-nal marginal distributions.
 Let $f_1, \dots, f_m  \in \mc F$  and consider the vector $(\tilde{\mb G}_n f_1, \dots, \tilde{\mb G}_n f_m)$.
 By the Cram\'er-Wold theorem, this vector convergences in distribution to $(\tilde{\mb G} f_1, \dots, \tilde{\mb G} f_m)$ 
 if and only if its canonical scalar product with any vector $\bs \lambda \in \R^m$ converges in distribution to \mbox{$ (\tilde{\mb G} f_1, \dots, \tilde{\mb G} f_m) \cdot \bs \lambda$}.
 Therefore, it is enough to restrict our attention to $\tilde{\mb G}_{n} f  = \tilde{\mb G}_n \sum_{j=1}^m \lambda_j f_j$ for a fixed vector $\bs \lambda = (\lambda_1, \dots, \lambda_m)^t \in \R^m$.
 Exceptional sets do not cause a problem here,
 even though the weak convergence shall be verified for uncountably many vectors $\bs \lambda \in \R^m$.
 This is ensured by the extended Cram\'er-Wold device;
 see Satz~3.19 in \cite{pauly09}.
 The idea of his proof is that the characteristic functions of the above linear combinations are continuous in $\bs \lambda$.
 As a consequence, a verification of the weak convergence for all linear combinations with coefficients in a countable subset, e.g. $\bs \lambda \in \mb Q^m$, suffices.

 We apply Hoeffding's Theorem \citep[Theorem~3.2]{hoeffding52} to verify the desired conditional convergence in distribution.
 To this end, let $G_1, G_1', G_2, G_2', \dots$, be independent random variables with a uniform distribution on $\mc G$
 and define, in addition to $\tilde{\mb G}_{n}f $, 
 a conditionally independent copy thereof, $\tilde{\mb G}'_{n}f = \sqrt{n}(\RPn' - \Pn^{\mc G})f = \frac{1}{\sqrt{n}} \sum_{i=1}^n ( f ( G_i' (\b X_i)) - \int_{\mc G} f(g_i(\b X_i)) \d Q(g_i) )  $.
 We need to analyze the unconditional asymptotic behaviour of the pair
  $(\tilde{\mb G}_{n} f, \tilde{\mb G}'_{n} f)  = \Big(\tilde{\mb G}_n \sum_{j=1}^m \lambda_j f_j, \tilde{\mb G}_n' \sum_{j=1}^m \lambda_j f_j \Big) $.
  Note that this is a sum of i.i.d.\ random variables in mappings of the triples
  $(G_i(\b X_i), G_i'( \b X_i ), \b X_i)$, $i=1, \dots, n$.
  Write $\b X = \b X_1$ and $G=G_1$.
  Thus, as $n \rightarrow \infty$, the classical multivariate central limit theorem yields its convergence in distribution  to a  bivariate normal distribution with expectations 
  \begin{align*}
    E \Big( f ( G (\b X)) - \int_{\mc G} f(g(\b X)) \d Q(g) \Big) 
   & = E \Big( E \Big( f ( G (\b X))  \mid \b X \Big) - \int_{\mc G} f(g(\b X)) \d Q(g) \Big) = 0,
  \end{align*}
  variances
  \begin{align*}
    E \Big[\Big( f ( G (\b X)) & - \int_{\mc G} f(g(\b X)) \d Q(g) \Big)^2 \Big]
       =  E \Big[ E \Big( \Big( f ( G (\b X)) - \int_{\mc G} f(g(\b X)) \d Q(g) \Big)^2  \mid \b X \Big) \Big] \\
      & = E \Big[ \int_{\mc G} f ( g (\b X))^2 \d Q(g) 
      - \Big( \int_{\mc G} f ( g (\b X)) \d Q(g) \Big)^2 \Big] \\
      & = \P \Big[ \int_{\mc G} f(g(\cdot))^2 \d Q(g) 
	- \Big( \int_{\mc G} f(g(\cdot)) \d Q(g) \Big)^2 \Big]
	= \P (Q_{\cdot}f^2 - (Q_{\cdot}f)^2),
  \end{align*}
  and covariances
  \begin{align*}
   & E \Big[\Big(  f ( G (\b X))  - \int_{\mc G} f(g(\b X)) \d Q(g) \Big)
    \Big( f ( G' (\b X))  - \int_{\mc G} f(g(\b X)) \d Q(g) \Big)\Big] \\
     & = E \Big( E \Big( \Big(  f ( G (\b X))  - \int_{\mc G} f(g(\b X)) \d Q(g) \Big)
    \Big( f ( G' (\b X))  - \int_{\mc G} f(g(\b X)) \d Q(g) \Big) \mid \b X\Big)  \Big) \\ 
    & = E \Big( E \Big(  f ( G (\b X))  - \int_{\mc G} f(g(\b X)) \d Q(g) \mid \b X \Big)
    E\Big( f ( G' (\b X))  - \int_{\mc G} f(g(\b X)) \d Q(g) \mid \b X \Big)  \Big)  = 0 ,
  \end{align*}
  where the second equality in the previous display is due to the conditional independence of both components given $\b X$.
 Thus, Hoeffding's Theorem implies that the conditional randomization distribution converges in probability
 to a centered normal distribution with variance $\P (Q_{\cdot}f^2 - (Q_{\cdot}f)^2)$.
 
 Next, we are going to prove the conditional tightness of the randomization empirical process in outer probability by verifying an asymptotic equicontinuity condition.
 Here, conditional tightness means, as defined in Theorem~2.9.6 in \cite{vaart96:_weak}, that
 $$ \sup_{h \in BL_1} | E_G h(\tilde{\mb G}_n) - E h(\tilde{\mb G}) |$$
 goes to zero in outer probability and that the sequence $\tilde{\mb G}_n$ is asymptotically measurable.
 Here, $BL_1$ again denotes the class of all bounded Lipschitz-continuous functions $h: \ell^\infty(\mc F) \rightarrow [0,1]$  with Lipschitz constant at most 1 and $E_G$ means integration with respect to $\d Q(g)$.
 Parts of the proof of $(i) \Rightarrow (ii)$ of the just mentioned Theorem~2.9.6 can be paralleled; for example, $\tilde{\mb G}_n$ converges unconditionally to a tight limit because this holds for both of the normalized processes in the following difference:
 \begin{align*}
  \tilde{\mb G}_n = \sqrt{n} (\RPn - \tilde \P) - \sqrt{n} (\Pn^{\mc G} - \tilde \P).
 \end{align*}
 The first process, $ \sqrt{n} (\RPn - \tilde \P)$, converges weakly because $\mc F$ is assumed to be $\tilde \P$-Donsker,
 and the second, $\sqrt{n} (\Pn^{\mc G} - \tilde \P)$, because $\tilde {\mc F}$ is assumed to be $\P$-Donsker.
 Hence, $\tilde{\mb G}_n$ is asymptotically measurable.
 
 Comparing with the other arguments in the proof of Theorem~2.9.6, 
 it only remains to show that 
 $E^*\|\tilde{\mb G}_n \|_{\mc F_\delta}$ goes to zero as $n \rightarrow \infty$ followed by $\delta \rightarrow 0$.
 Here, $\mc F_\delta = \{f-g : f,g \in \mc F, \rho_\P(f,g) < \delta\}$, where $\rho_\P$ is a suitable seminorm on $\mc F$. 
 The desired convergence holds because, 
 as explained above, $\tilde{\mb G}_n$ converges unconditionally weakly to a tight limit;
 and this is equivalent to the mean version of the asymptotic equicontinuity condition; see Lemma~2.3.11 in \cite{vaart96:_weak}.
\end{proof}

\begin{proof}[Proof of Remark~\ref{rem:kiefer}]
 Obviously, the transformed Brownian motion process has mean zero.
 For $f,h \in \ell^{\infty}(\mc F)$, 
 \begin{align*}
  E \mb W_{\tilde \P}( f - Q_{\cdot} f ) \cdot \mb W_{\tilde \P}( h - Q_{\cdot} h )
  & = \tilde \P(fh) - \tilde \P (f Q_{\cdot} h ) - \tilde \P ((Q_{\cdot} f)  h ) + \tilde \P ((Q_{\cdot} f) (Q_{\cdot} h)) \\
  & = \P (Q_\cdot(fh)-Q_\cdot fQ_\cdot h),
 \end{align*}
 \tcb{since $(Q_{g(\cdot)}f )(Q_{g(\cdot)}h) =(Q_{\cdot}f )(Q_{\cdot}h) $ and 
 $ Q_{\cdot} (f (Q_{\cdot} h) ) =(Q_{\cdot}f )(Q_{\cdot}h) $ due to the group structure of $\mc G$.}
\end{proof}

\begin{proof}[Proof of Theorem~\ref{thm:delta_meth}]
Large parts of the proof of Theorem~3.9.11 in \cite{vaart96:_weak} apply as it does not make use of the particular structure of the bootstrap empirical process considered there, except that the randomization empirical process is centered at $\PnG$ and not at $\Pn$.
As in the proof the just mentioned theorem, we assume without loss of generality that $\varphi'_{\tilde \P}$ is defined on the whole space $\ell^\infty(\mc F)$.

As explained in the proof of Theorem~\ref{thm:main},
the weak convergences of both processes $\sqrt{n}(\RPn - \tilde \P)$ and $\sqrt{n}(\PnG - \tilde \P)$ hold unconditionally.
As a consequence, both sequences $\sqrt{n}(\varphi(\RPn) - \varphi(\tilde \P) )$ and $\sqrt{n}(\varphi(\PnG) - \varphi(\tilde \P) ) $ also converge unconditionally because the classical functional delta-method applies:
\begin{align*}
 \sqrt{n}(\varphi(\RPn) - \varphi(\tilde \P) ) = \varphi'_{\tilde \P}(\sqrt n (\RPn - \tilde \P)) + o^*_p(1) & \quad \text{and}
 \\
 \sqrt{n}(\varphi(\PnG) - \varphi(\tilde \P) ) = \varphi'_{\tilde \P}(\sqrt n (\PnG - \tilde \P)) + o^*_p(1) &
 \end{align*}
 unconditionally as $n \rightarrow \infty$.
The rest of the proof again continues along the lines of Theorem~3.9.11 in \cite{vaart96:_weak}:
a subtraction of both equations in the previous display gives that $\sqrt{n}(\varphi(\RPn) - \varphi(\PnG) ) - \varphi'_{\tilde \P}(\sqrt n (\RPn - \PnG))$ converges unconditionally to zero in outer probability from which the desired result follows. 

The asymptotic variance is the limit of the following conditional variances, 
\begin{align*}
& var(\sqrt{n} \cdot \varphi_{\tilde \P}'(\RPn - \Pn^{\mc G}) \mid \b X_1, \b X_2, \dots )  \\
& = \frac1n \sum_{i=1}^n var( \varphi_{\tilde \P}'( \delta_{G_i ( \b X_i}) -  \int_{\mc G} \delta_{ g  ( \b X_i) } \d Q(g) ) \mid \b X_1, \b X_2, \dots ) \\
 & = 
 \frac1n \sum_{i=1}^n \Big\{ \int_{\mc G}  IF^2_{\varphi,\tilde \P}(g(\b X_i)) \d Q(g) 
 - \Big[ \int_{\mc G} IF_{\varphi,\tilde \P}(g(\b X_i)) \d Q(g) \Big]^2 \Big\}.
\end{align*}
By the strong law of large numbers, this converges to $$\int_{\R^d} \Big\{\int_{\mc G} IF^2_{\varphi,\tilde \P}(g(\b x)) \d Q(g) 
 - \Big[ \int_{\mc G} IF_{\varphi,\tilde \P}(g(\b x)) \d Q(g) \Big]^2 \Big\} \d \P(\b x)$$ almost surely as $n \rightarrow \infty$.
\end{proof}

\begin{proof}[Proof of Corollary~\ref{cor:stud}]
A combination of the assumed convergences in~\eqref{eq:conv_IF} with 
$$ \Pn (IF_{\varphi, \P}^2) = \Pn (\varphi'_{\P})^2 \oPo \P (\varphi'_{\P})^2 \quad \text{and} \quad  \RPn (IF_{\varphi, \tilde \P}^2) = \RPn (\varphi'_{\tilde \P})^2 \oPo \tilde \P (\varphi'_{\tilde \P})^2 $$
by the law of large numbers implies that  
$\wh \sigma_{\varphi,\P_n}^2 \oPo \sigma^2_{\varphi'_\P} \in (0, \infty)$  
and $\tilde \sigma_{\varphi,\RPn}^2 \oPo  \tilde \sigma^2_{\varphi'_{\tilde \P}} \in (0, \infty) $.
 The asymptotic exactness of the proposed test now follows by combining the consistency of the variance estimators with Theorem~\ref{thm:main} through Slutzky's lemma.
 Note here that convergence in probability of $\tilde \sigma_{\varphi,\RPn}^2$ is equivalent to conditional convergence in probability given $\b X_1, \b X_2, \dots.$
 
 The finite sample exactness of such randomization tests under restricted null hypotheses $H_0$ of $\mc G$-invariance is well-known and not further discussed here; see e.g.\ Theorem~1 in \cite{hemerik18}.
 We just remark that $H_0$ implies that the randomized studentized test statistic has the same unconditional distribution as the studentized test statistic.
 The condition $\varphi(\Pn^{\mc G})=\theta_0$ is necessary to ensure this.
\end{proof}

\section{Lipschitz condition for the consistency of variance estimators}
\label{sec:lipschitz}

The conditions in \eqref{eq:conv_IF} in the main manuscript \cite{dobler20} hold, for example,
if the influence function  satisfies a pointwise Lipschitz condition with square-integrable Lipschitz constants $L(\b X)$ and $L(G(\b X))$:
$$ 
\Pn| IF_{\varphi,\Pn} - IF_{\varphi,\P} | 
= \frac1n\sum_{i=1}^n | IF_{\varphi,\Pn}(\b X_i) - IF_{\varphi,\P}(\b X_i)|
\leq \frac1n\sum_{i=1}^n L(\b X_i) \cdot d(\Pn, \P)
= (\Pn L) \cdot d(\Pn,\P),$$
 where the metric $d(\cdot, \cdot)$ metrizes weak convergence. 
 Indeed, by Jensen's and the Cauchy-Schwarz inequality,
\begin{align*}
  & | \Pn IF^2_{\varphi,\Pn} - \Pn IF^2_{\varphi,\P} | 
\leq \Pn| IF^2_{\varphi,\Pn} - IF^2_{\varphi,\P} | 
= \Pn| IF_{\varphi,\Pn} - IF_{\varphi,\P} | \cdot | IF_{\varphi,\Pn} + IF_{\varphi,\P}| \\
& \leq \Big[ \Pn | IF_{\varphi,\Pn} - IF_{\varphi,\P} |^2 \Pn | IF_{\varphi,\Pn} - IF_{\varphi,\P} + 2  IF_{\varphi,\P} |^2 \Big]^{1/2} \\
& \leq \Big[ \Pn L^2 \cdot  d^2(\Pn,\P) \cdot 2 \Big( \Pn L^2 \cdot d^2(\Pn,\P) + 2 \Pn IF^2_{\varphi,\P} \Big) \Big]^{1/2};
\end{align*}
 analogous inequalities hold for $\Pn$ and $\P$ replaced by $\RPn$ and $\tilde \P$, respectively.

\section{Influence function for the Mann-Whitney effect estimator and consistent variance estimates}
\label{app:MW}

\tcb{We will first discuss the case in which all individuals were matched and thus $n$ pairs were formed. The case of unequal sample sizes will be discussed in the next subsection.}

\subsection{The completely paired data case}

For estimating the asymptotic variances, it remains to derive and estimate the influence function of the Mann-Whitney effect estimator.
Denote by $y_j(t) = P(X_{ij} > t)= S_j(t) H_j(t)$, where $H_j(t) = P(C_{ij}>t)$ is the censoring survival function, $j=1,2$.
Furthermore, we are going to use the cumulative hazard functions $\Lambda_j(t) = - \int_0^t \frac{\d S_j(u)}{S_j(u-)}$.
As references for the following Hadamard-derivatives, see Example~3.9.19 and Lemma~3.9.30 in \cite{vaart96:_weak}.
The influence function of the $j$th Nelson-Aalen estimator evaluated at $t \in [0,\tau]$ is given by
\begin{align*}
 & \int_0^t \frac{\d 1\{X_{ij} \leq u , \delta_{ij} = 1\}}{y_j(u)} - \int_0^t \frac{1\{ X_{ij} \geq u\}}{y_j^2(u)} \d P(X_{ij} \leq u , \delta_{ij} = 1) 
 & = \frac{\delta_{ij} 1\{ X_{ij} \leq t \}}{y_j(X_{ij})} - \tilde\sigma_j^2(t \wedge X_{ij}).
\end{align*}
Here we used $P(X_{ij} \leq u , \delta_{ij} = 1) = - \int_0^u H_j(v-) \d S_j(v) $ and the abbreviation $\tilde\sigma_j^2(t) = - \int_0^t \frac{\d S_j(u)}{H_j(u-) S^2_j(u-)} $.
Defining $\frac 0 0 = 0$, 
the influence function of the Kaplan-Meier estimator is given by
$$  S_j(t) \int_0^t \frac{\d \Big( \frac{\delta_{ij} 1\{ X_{ij} \leq u \}}{y_j(X_{ij})} - \tilde\sigma_j^2(u \wedge X_{ij}) \Big) }{1-\Delta \Lambda_j(u)}
=  S_j(t) \Big[ \frac{\delta_{ij} 1\{ X_{ij} \leq t \}}{y_j(X_{ij})(1-\Delta \Lambda_j(X_{ij}))} - \int_0^{t \wedge X_{ij}} \frac{\d \Lambda_j(u)}{H_j(u-) S_j(u)  }  \Big];$$
see \cite{reid81} for a similar representation of the influence function.
Here we used the notation $\Delta f(t) := f(t) - f(t-)$ to denote the jump size of a right-continuous function $f$ at $t$.
For future use, we abbreviate the integral on the right-hand side in the previous display by $\sigma_j^2(t \wedge X_{ij})$.

The final map to obtain the Mann-Whitney effect is the modified Wilcoxon functional.
Its Hadamard-derivative as derived 
in the supplementary material to \cite{dobler18test} will be used in the following form: the derivative at $(S_1, S_2)$ is given by
\begin{align}
\label{eq:had_wilc}
 (h_1, h_2) \longmapsto \frac12 \Big[ - \int_{[0,\tau)} h_1 \d S_2 + \int_{[0,\tau)} S_2 \d h_1 + \int_{[0,\tau)} h_2 \d S_1 - \int_{[0,\tau)} S_1 \d h_2 \Big],
\end{align}
where the integrals with respect to $h_j$ are defined via integration by parts if $h_j$ has unbounded variation.
Now, from the preparations above it follows that the influence function of the functional $\phi$, which maps the empirical process to the Mann-Whitney effect estimate, is $IF_{\phi, \P}(X_{i1}, X_{i2}, \delta_{i1}, \delta_{i2})$
\begin{align*}
 =  \frac12 \sum_{j=1}^2 (-1)^{j} \Big\{ & \int_{[0,\tau)} S_j(t) \Big[ \frac{\delta_{ij} 1\{ X_{ij} \leq t \}}{y_j(X_{ij})(1-\Delta \Lambda_j(X_{ij}))} - \int_0^{t \wedge X_{ij}} \frac{\d \Lambda_j(u)}{H_j(u-) S_j(u)  }  \Big] \d S_{3-j}(t) \\
& - \int_{[0,\tau)} S_{3-j}(t) \d \Big[ \frac{\delta_{ij} S_j(t) 1\{ X_{ij} \leq t \}}{y_j(X_{ij})(1-\Delta \Lambda_j(X_{ij}))} - S_j(t) \int_0^{t \wedge X_{ij}} \frac{\d \Lambda_j(u)}{H_j(u-) S_j(u)  }  \Big] \Big\}.
\end{align*}
Writing $k = k(j) = 3-j$, the first integral from $0$ to $\tau$ can be simplified to 
\begin{align*}
 & \frac{\delta_{ij} 1\{X_{ij} < \tau\}}{H_j(X_{ij}-) S_j(X_{ij})} 
 \Big[
 \int_{(X_{ij},\tau)} S_j(t)   \d S_{k}(t)
 +  S_j(X_{ij}) S_k(X_{ij})
 - S_j(X_{ij}) S_k(X_{ij}-)
 \Big] \\
 & - \int_{[0,X_{ij})} S_j(t) \sigma^2_j(t) \d S_{k}(t) 
 - \sigma^2_j(X_{ij}) \int_{[X_{ij},\tau)} S_j(t) \d S_{k}(t).
\end{align*}
Likewise, the second integral  from $0$ to $\tau$ equals
\begin{align*}
   & \frac{\delta_{ij} 1\{X_{ij} < \tau\}}{H_j(X_{ij}-) S_j(X_{ij})} \Big[ \int_{(X_{ij},\tau)} S_{k}(t) \d S_j(t) + S_j(X_{ij}) S_k(X_{ij}) \Big] \\
   & - \int_{[0,X_{ij})} S_k(t) \sigma_j^2(t) \d S_j(t)
   - \sigma_j^2(X_{ij}) \int_{[X_{ij},\tau)} S_k(t) \d S_j(t) 
   - \int_{[0,X_{ij}]} 1\{t < \tau \} \frac{S_k(t) \d \Lambda_j(t)}{H_j(t-)}
\end{align*}
Brought together, the influence function $IF_{\phi, \P}$ simplifies to 
\begin{align}
\begin{split}
\label{eq:IFMW}
  IF^{(1)}_{\phi, \P}(X_{i1}, \delta_{i1}) & - IF^{(2)}_{\phi, \P}(X_{i2}, \delta_{i2}) \\
  = \frac12 \sum_{j=1}^2 (-1)^{j} \Big\{ 
  & \frac{\delta_{ij}  1\{X_{ij} < \tau\}}{H_j(X_{ij}-) S_j(X_{ij})} 
 \int_{(X_{ij},\tau)}  [ S_j(t)   \d S_{k}(t) - S_k(t)   \d S_{j}(t) ]
 - \delta_{ij} 1\{X_{ij} < \tau\} \frac{S_k(X_{ij}-)}{H_j(X_{ij}-)} \\
 & - \int_{[0,X_{ij})} \sigma_j^2(t) [ S_j(t)   \d S_{k}(t) - S_k(t)   \d S_{j}(t) ] \\
 & - \sigma_j^2(X_{ij}) \int_{[X_{ij}, \tau)}  [ S_j(t)   \d S_{k}(t) - S_k(t)   \d S_{j}(t) ]
 +  \int_{[0,X_{ij}]} 1\{t < \tau \} \frac{S_k(t) \d \Lambda_j(t)}{H_j(t-)} \Big\}.
\end{split}
\end{align}

For simplifying readability we now omit the notion of $n$ in subscripts.
The variance of the above influence function can be estimated by replacing $S_j, H_j$, and $\sigma_j^2(t)$ with 
$\wh S_j, \wh H_j$, and 
$$\wh \sigma_j^2(t) =n \int_0^t \frac{\wh \Lambda_j(\d u)}{ Y_j(u) (1-\Delta\wh \Lambda_j(u))},$$
respectively, 
where $\wh H_j$ are the Kaplan-Meier estimators for the censoring survival functions $H_j$,
$Y_j(t) = \sum_{i=1}^{n} 1\{X_{ij} \geq t\}$ are the number at risk processes and $\wh \Lambda_j(t) = - \int_0^t \frac{\d \wh S_j(u)}{\wh S_j(u-)}$ are the Nelson-Aalen estimators, $j=1,2$.
Note that the estimators $\wh S_1, \wh S_2 $, that appear in the denominators, have to be slightly adjusted to prevent division by $\wh S_j(X_{ij}) =0$ if the latest observed time in group $j \in \{ 1,2\}$ is uncensored and less than $\tau$.
However, due to the condition in \eqref{eq:cond_mw}, the probability that such a case occurs rapidly decreases with increasing sample size.

Due to the uniform consistency of the involved estimators, the asymptotic boundedness of $S_j$ and $H_j$ away from 0, and the continuity of the above influence function as a functional in $\wh S_j, \wh H_j$, and $\wh \sigma_j^2(t), j=1,2, $
we have established the consistency of the following  estimator for the variance of $\sqrt{n} (\wh p-p)$:
$$\wh \sigma^2_{\phi,\Pn} = \frac1n \sum_{i=1}^n \Big\{ IF_{\phi, \Pn}(X_{i1}, X_{i2}, \delta_{i1}, \delta_{i2}) - \frac1n \sum_{j=1}^n IF_{\phi, \Pn}(X_{j1}, X_{j2}, \delta_{j1}, \delta_{j2}) \Big\}^2. $$

Now, for the randomized Mann-Whitney effect estimator, we similarly receive the following consistent estimator for the asymptotic variance of $\sqrt{n}(\tilde p - \frac12)$:
$$\tilde \sigma_{\phi,\RPn}^2 = \frac1n \sum_{i=1}^n \Big\{ IF_{\phi, \RPn}(G_i  (X_{i1}, X_{i2}, \delta_{i1}, \delta_{i2})) - \frac1n \sum_{j=1}^n IF_{\phi, \RPn}(G_j  (X_{j1}, X_{j2}, \delta_{j1}, \delta_{j2})) \Big\}^2. $$

\subsection{The case of unequal sample sizes}
\label{sec:MW_uneq}
\tcb{If some individuals in the data set could not matched to anyone else, we obtain a situation as described in Remark~\ref{rem:MWuneq} in the main manuscript, with $n_p$ paired observations and $n_j$ single group $j$ observation observations; $j=1,2$.
The available data may then be represented as
$$ (X_{i1}, X_{i2}, \delta_{i1}, \delta_{i2}), \quad i=1, \dots, n_p + n_1 + n_2 =:n ,  $$
where $X_{i1} = \delta_{i1} \equiv 0$ if $i > n_p + n_1$ and $X_{i2} = \delta_{i2} \equiv 0$ if $n_p < i \leq n_p+n_1$; if the Kaplan-Meier estimator incorporates a data point $(0,0)$, the estimate is not changed at all..
A comparison with \eqref{eq:IFMW}  reveals that the asymptotic linearization of $\wh p_n$ can be written with the help of the three case-specific empirical processes,
$\P^{(p)}_{n_p} = \frac1{n_p} \sum_{i=1}^{n_p} \delta_{(X_{i1}, X_{i2}, \delta_{i1}, \delta_{i2})}, 
\P^{(1)}_{n_1} = \frac1{n_1} \sum_{i=1}^{n_1} \delta_{(X_{i1}, 0, \delta_{i1}, 0)},
\P^{(2)}_{n_2} = \frac1{n_2} \sum_{i=1}^{n_2} \delta_{(0,X_{i2}, 0, \delta_{i2})}$.
%
%
Thus, because both Kaplan-Meier estimators are in general based on different sample sizes $n_p + n_1$ and $n_p + n_2$, the following asymptotic linearization holds:
\begin{align}
\begin{split}
 \label{eq:MWuneq2}
 \sqrt{n}(\wh p_n - p)  & = \sqrt{n} (\psi(\wh S_{1,n}, \wh S_{2,n}) - \psi(S_{1}, S_{2})) \\ 
 & = n^{1/2} \Big( \frac{1}{n_p + n_1}(n_p \P_{n_p}^{(p)} + n_1  \P_{n_1}^{(1)} ) IF^{(1)}_{\phi, \P} 
 - \frac{1}{n_p + n_2}(n_p \P_{n_p}^{(p)} + n_2  \P_{n_2}^{(2)} ) IF^{(2)}_{\phi, \P}   \\
 & \quad - E\Big( \frac{1}{n_p + n_1}(n_p \P_{n_p}^{(p)} + n_1  \P_{n_1}^{(1)} ) IF^{(1)}_{\phi, \P} 
 - \frac{1}{n_p + n_2}(n_p \P_{n_p}^{(p)} + n_2  \P_{n_2}^{(2)} ) IF^{(2)}_{\phi, \P} \Big) \Big).
 \end{split}
\end{align}
%
%
%
%
%
Here we have slightly abused the notation by letting $IF^{(1)}_{\phi, \P}$ and $IF^{(2)}_{\phi, \P}$ be functions of four arguments while, in fact, they are constant with respect to the arguments $(2,4)$ and $(1,3)$, respectively.
%
%
%
%
%
%
%
%
Its variance is equal to
\begin{align*}
& n \Big( \frac{n_p}{(n_p + n_1)^2} var(\P_{1}^{(p)} IF^{(1)}_{\phi, \P})
- 2 \frac{n_p}{(n_p + n_1)(n_p + n_2)} cov (\P_{1}^{(p)} IF^{(1)}_{\phi, \P}, \P_{1}^{(p)} IF^{(2)}_{\phi, \P})
+ \frac{n_p}{(n_p + n_2)^2} var(\P_{1}^{(p)} IF^{(2)}_{\phi, \P}) \\
& + \frac{n_1}{(n_p + n_1)^2} var(\P_{1}^{(1)} IF^{(1)}_{\phi, \P})
+ \frac{n_2}{(n_p + n_2)^2} var(\P_{1}^{(2)} IF^{(2)}_{\phi, \P}) \\
& = n \cdot \Big( n_p var( \P_{1}^{(p)} \tilde{IF}_{\phi, \P} ) + n_1 var( \P_{1}^{(1)} \tilde{IF}_{\phi, \P} ) + n_2 var( \P_{1}^{(2)} \tilde{IF}_{\phi, \P} ) \Big),
\Big)
\end{align*}
where $\tilde{IF}_{\phi, \P} = IF^{(1)}_{\phi, \P} / (n_p + n_1) - IF^{(2)}_{\phi, \P}/ (n_p + n_2)$.
Due to the structure of the Hadamard-derivative of the Wilcoxon functional, cf.\ \eqref{eq:had_wilc}, it is clear that each $E(IF^{(j)}_{\phi, \P} (X_{i1}, X_{i2}, \delta_{i1}, \delta_{i2}) ) = 0$, $j=1,2$.
Consequently, we propose to use the pooled variance estimator,
$ \widehat \sigma^2_{\phi, \Pn} =  n^2 \cdot \Pn ( \tilde{IF}_{\phi, \Pn} - \Pn \tilde{IF}_{\phi, \Pn} )^2 $, which also has a very simple structure.
Similarly, the randomization version of this estimator, $\tilde \sigma^2_{\phi, \RPn}$, is obtained by replacing $\Pn$ with the randomization empirical process $\RPn$.
Sections~\ref{sec:MW_2} and \ref{sec:MW_uneq_simu} contain additional Monte-Carlo simulation results that are based on this variance estimator, in addition to those in the main manuscript.}

\section{Additional simulation results}
\label{sec:simus_add}

\subsection{Pearson correlation coefficient}



For simulating the type I error rates, i.e.\ under the null hypothesis $H: \rho_{Y,Z}=0$, we generated data according to the following distributions:
\begin{itemize}
 \item[I.] bivariate standard normal distribution, i.e.\ symmetry with respect to the axes and independence of both cooordinates;
 \item[II.] bivariate $t_5$-distribution, i.e.\ symmetry, no independence, and heavier tails than under the normal distribution;
 \item[III.] independent $\chi^2_5$ distributions, i.e.\ right-skewed marginal distributions;
 \item[IV.] a mixture distribution of the above bivariate $t_5$- (50\%) and $\chi^2_5$-distributions (50\%), i.e.\ no independence and right-skew marginal distributions.
\end{itemize}

We considered the sample sizes $n = 10, 15, \dots, 100$ and the significance levels $\alpha = 1\%, 5\%, 10\%$.
The following methods for finding critical values and randomization probabilities have been used:
\begin{itemize}
 \item a randomization method based on randomly mirroring the data with respect to the coordinate axes, i.e.\ the randomization group $\mc G^{symm.}$; this corresponds to finite exactness of the resulting test under symmetric distributions such as in the cases I and II;
 \item a randomization method based on randomly rotating the data with respect to the origin, i.e.\ the randomization group $\mc G^{rot.}$; this corresponds to finite exactness of the resulting test under rotation-invariant distributions such as in the cases I and II;
 \item Efron's classical bootstrap method \citep{efron1979bootstrap}, i.e.\ independently drawing the data with replacement; no finite sample exactness under any of the considered distributions;
 \item Random permutation of the first entries of all pairs; this corresponds to finite exactness of the resulting test if both components in the data pairs are independent such as in the cases I and~III;
 \item quantiles of the standard normal distribution; no finite sample exactness under any of the considered distributions.
\end{itemize}
Even though a similar permutation test has been used by \cite{diciccio17}, the test statistic used here is different from theirs, as explained in Section~\ref{sec:corr} in the main manuscript \cite{dobler20}, because we are using a different studentization.

We omitted the presentation of the simulated type I error probabilities based on the untransformed empirical correlation coefficient 
because 
the Fisher z-transformation improved the performance in most of the cases.
Furthermore, we do not display the results for the asymptotic tests which are based on normal quantiles; in non-normal scenarios their behaviour was far too liberal.

\begin{figure}
 \includegraphics[width=0.45\textwidth]{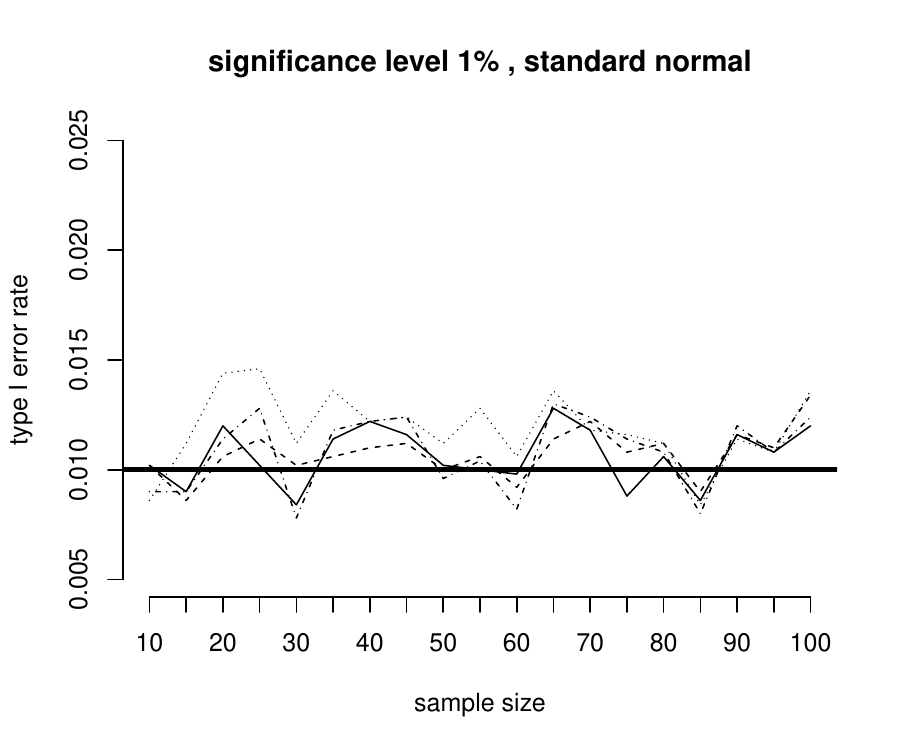}
 \includegraphics[width=0.45\textwidth]{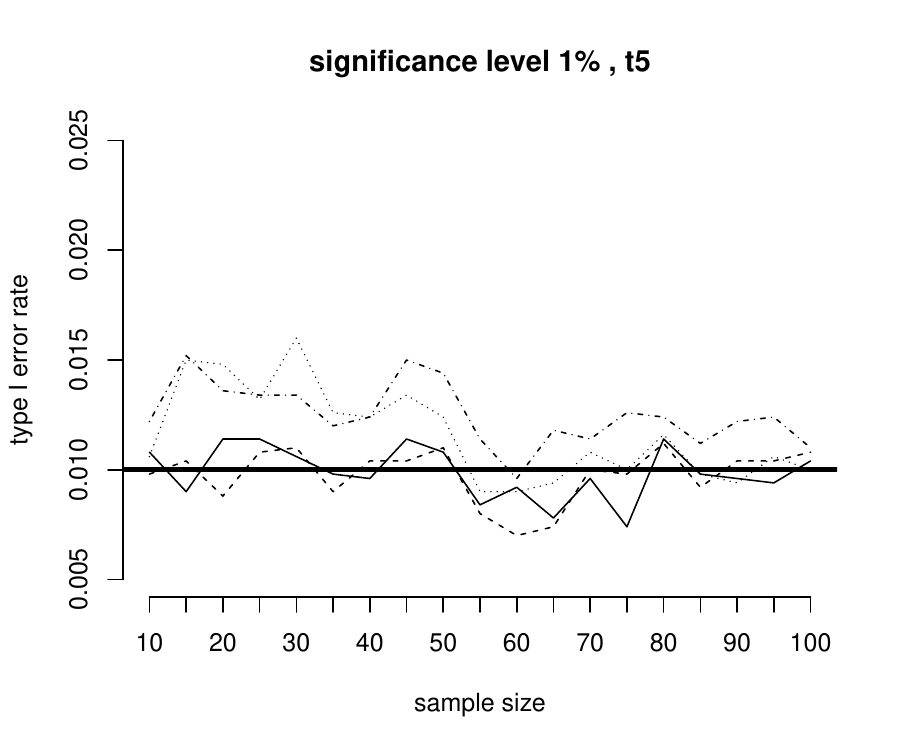}
 \\
 \includegraphics[width=0.45\textwidth]{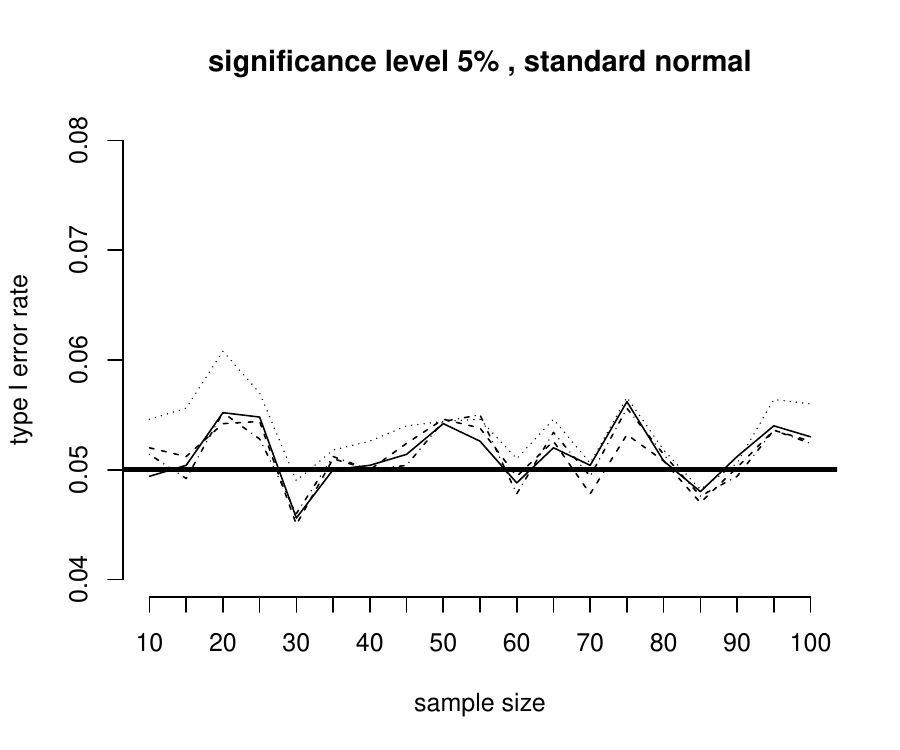}
 \includegraphics[width=0.45\textwidth]{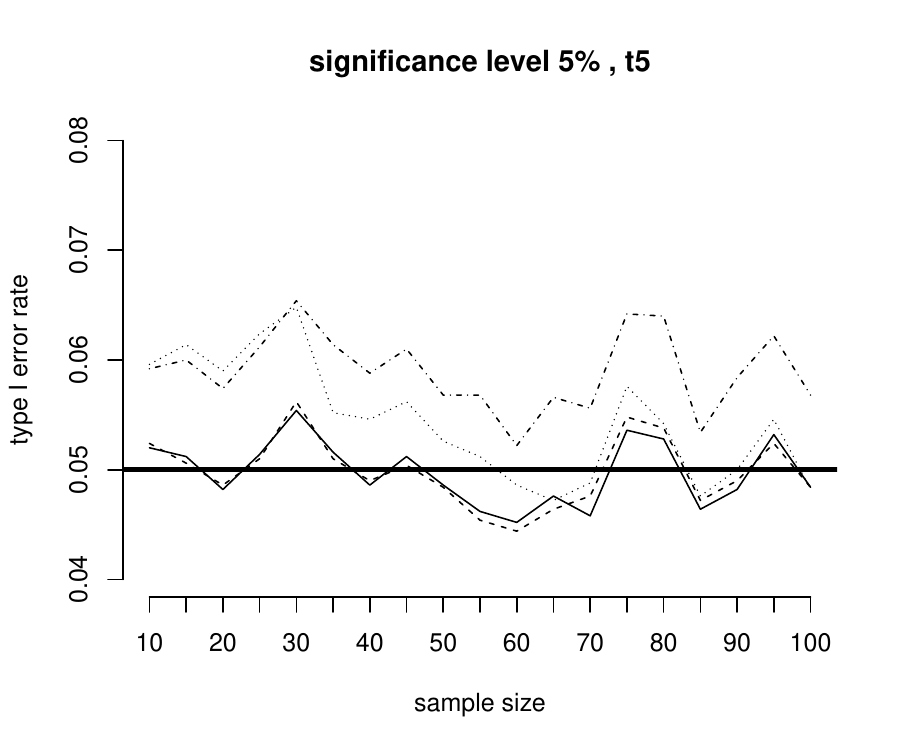}
 \\
 \includegraphics[width=0.45\textwidth]{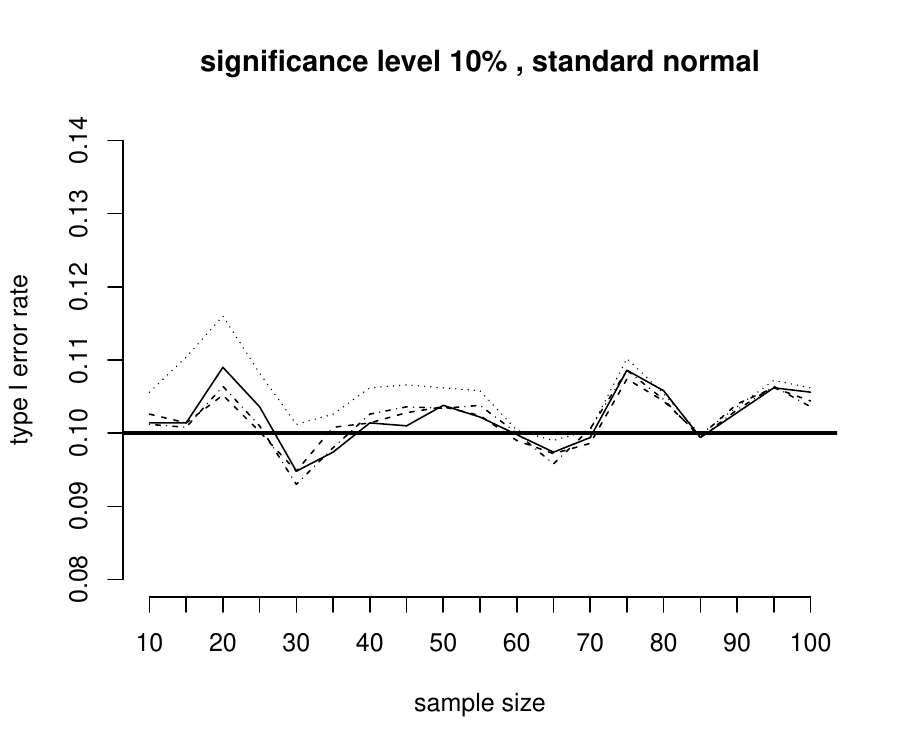}
 \includegraphics[width=0.45\textwidth]{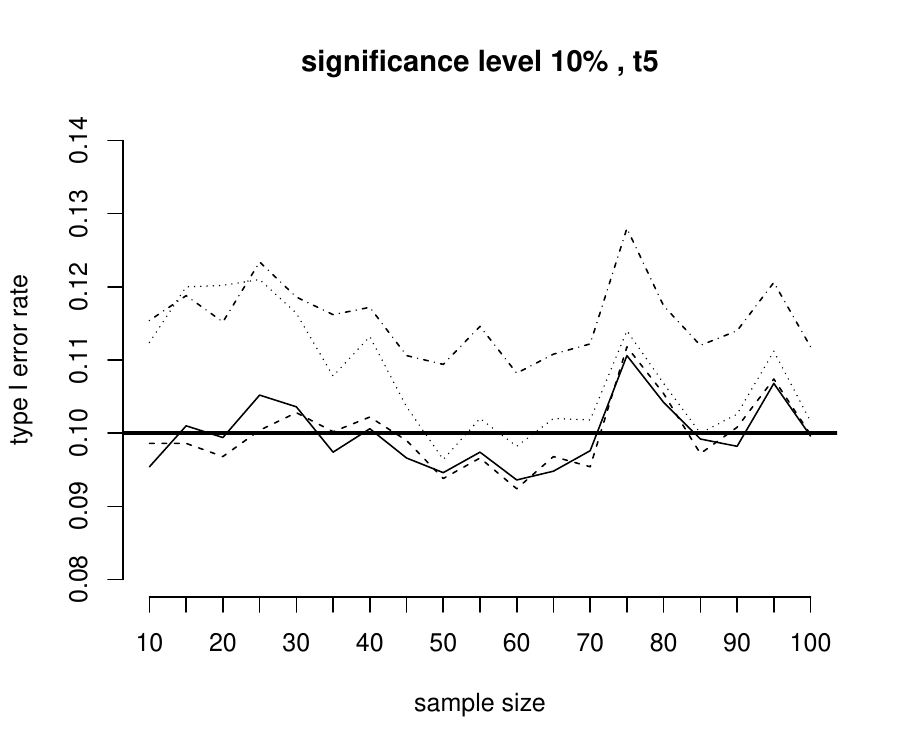}
 \caption{Simulated type I error rates of the correlation test with bivariate normally (left) and bivariate $t_5$ distributed data (right); based on mirroring (---), rotating (--~--), bootstrap ($\cdots$), permutation ($\cdot$ -- $\cdot$). 
 The horizonal line is the nominal significance level.}
 \label{fig:pearson1}
\end{figure}

\begin{figure}
 \includegraphics[width=0.45\textwidth]{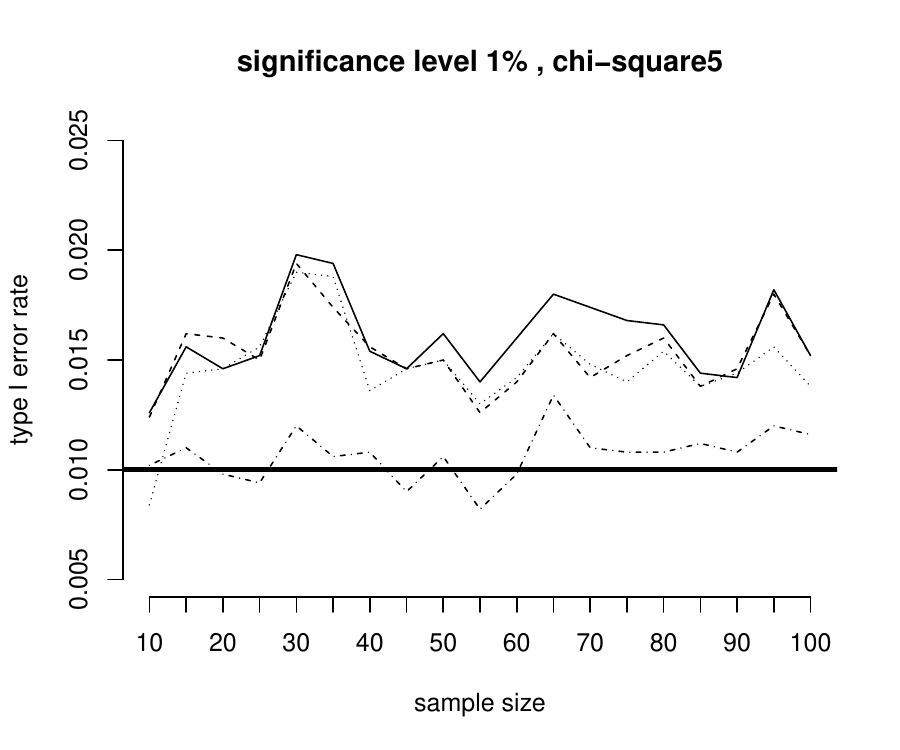}
 \includegraphics[width=0.45\textwidth]{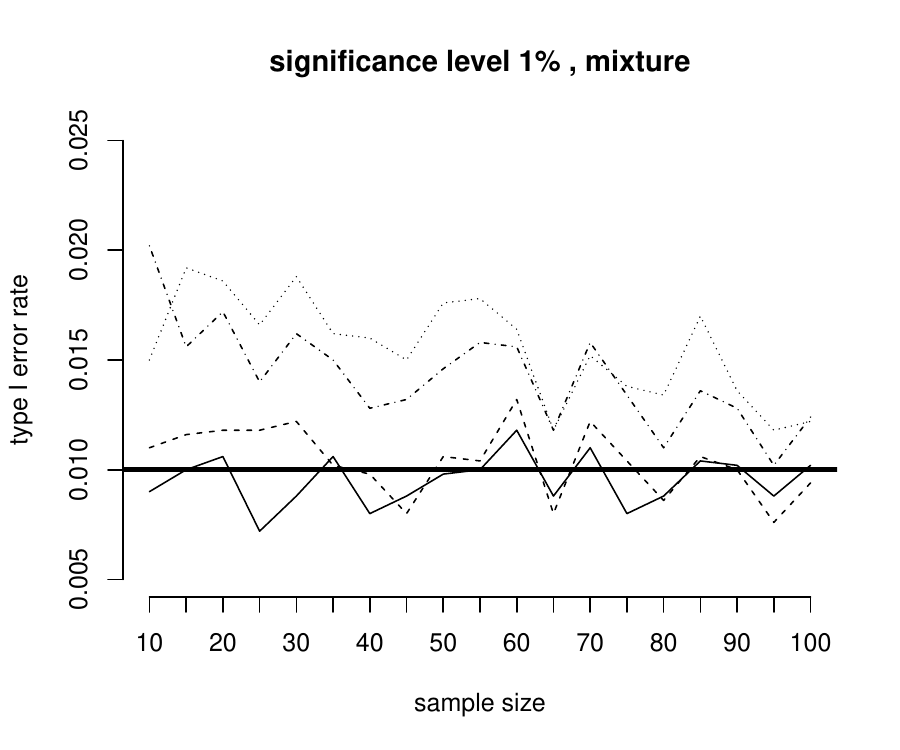}
 \\
 \includegraphics[width=0.45\textwidth]{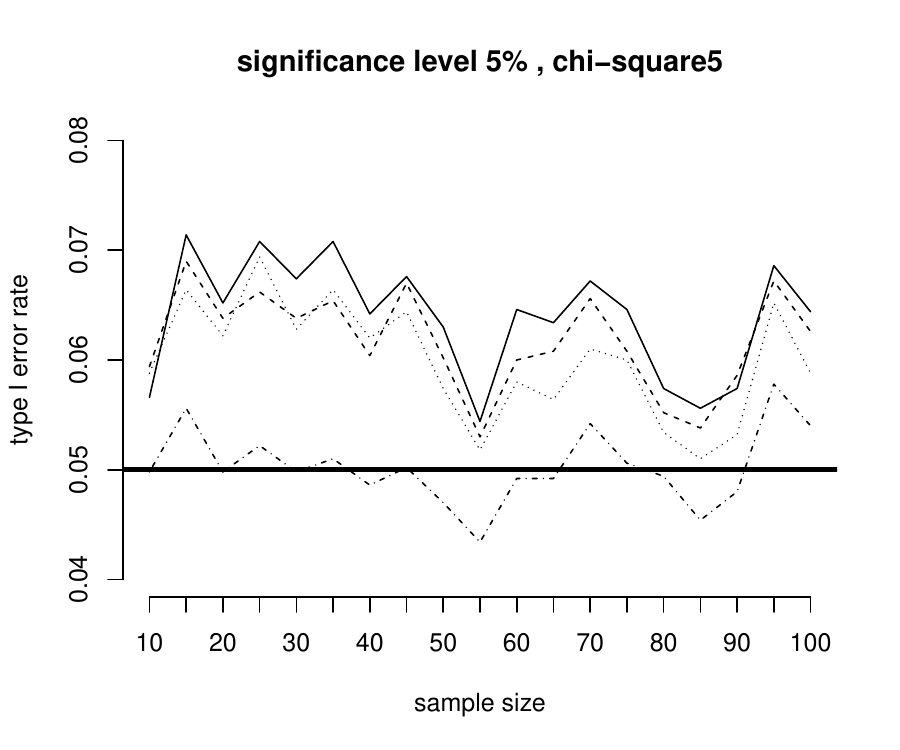}
 \includegraphics[width=0.45\textwidth]{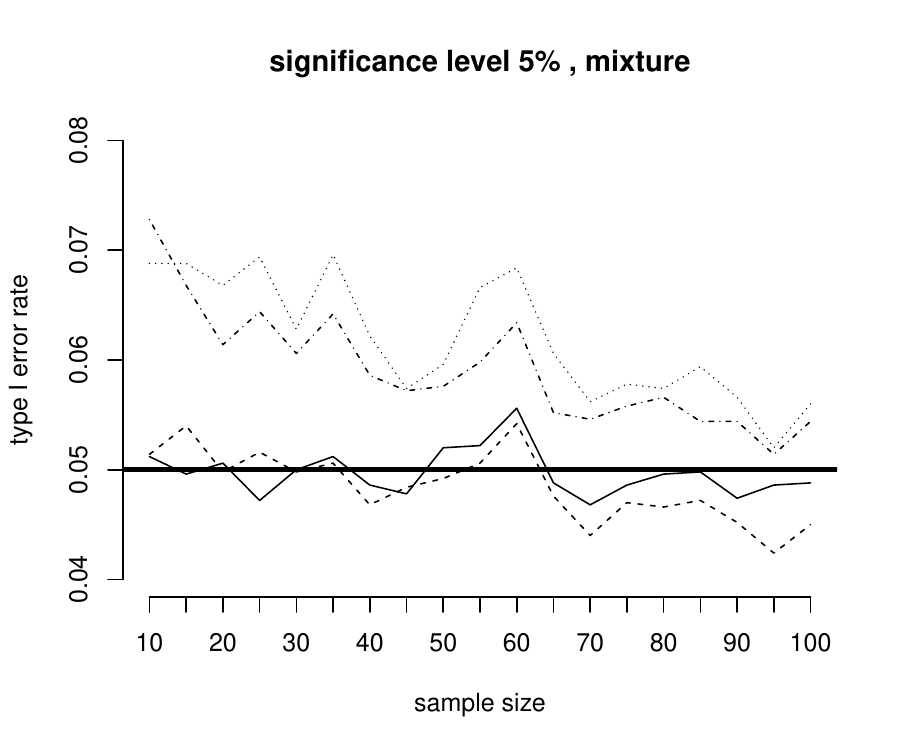}
 \\
 \includegraphics[width=0.45\textwidth]{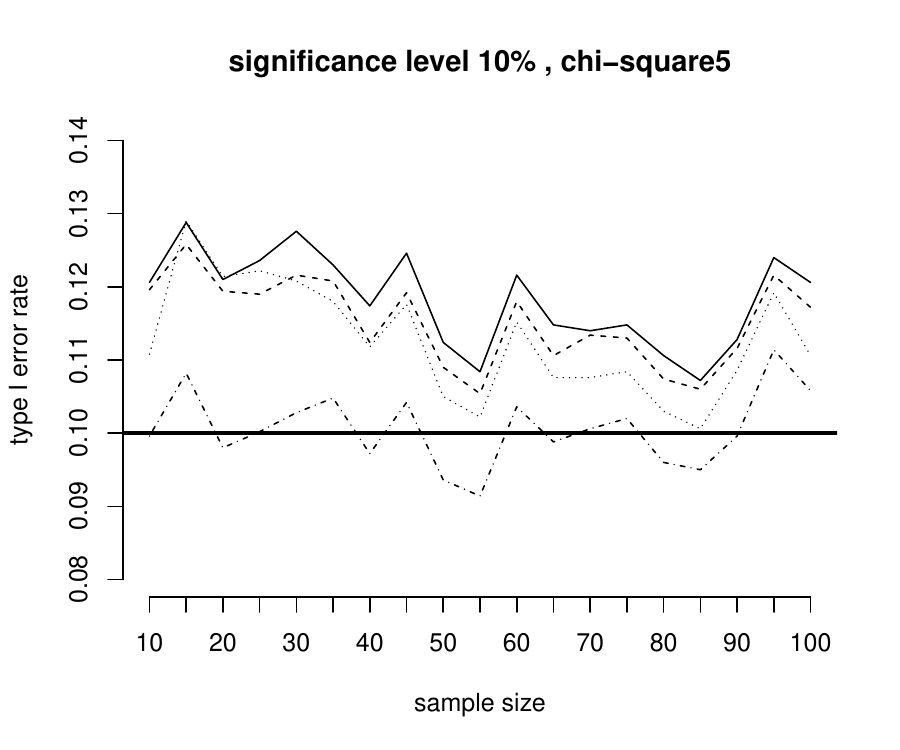}
 \includegraphics[width=0.45\textwidth]{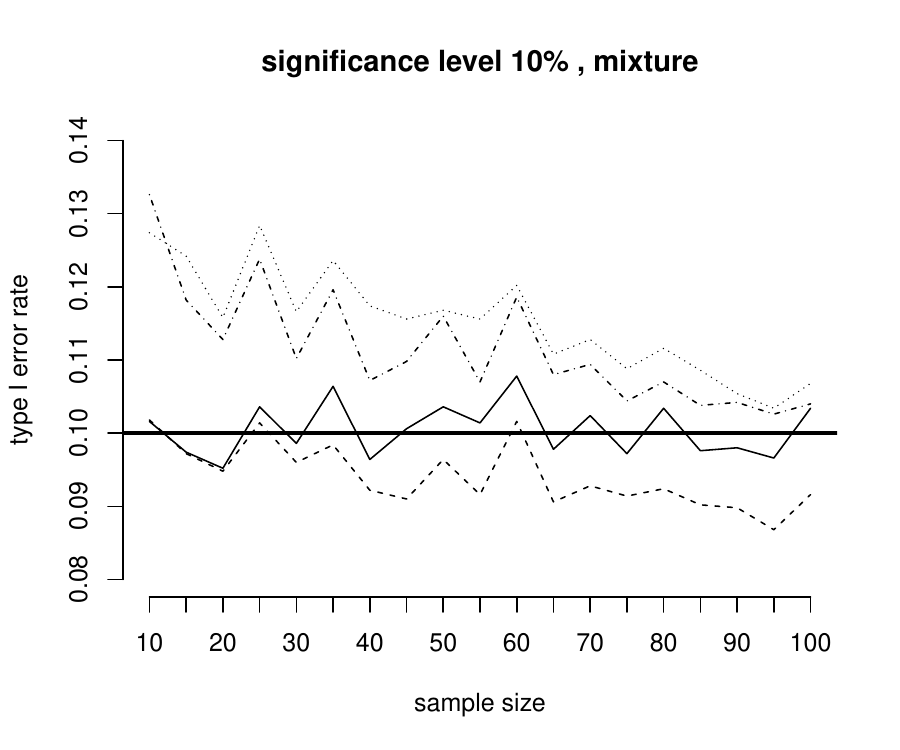}
 \caption{Simulated type I error rates of the correlation test with $\chi^2_5$ distributed data (left) and data from a bivariate $t_5$-$\chi^2_5$ mixture distribution (right); based on mirroring (---), rotating (--~--), bootstrap ($\cdots$), permutation ($\cdot$ -- $\cdot$). 
 The horizonal line is the nominal significance level.}
 \label{fig:pearson2}
\end{figure}

The plots in Figures~\ref{fig:pearson1} and~\ref{fig:pearson2} illustrate, first of all, the finite exactness of the randomization and permutation tests in the respective cases. 
Overall, both randomization-type tests show a similar accuracy.
It is more interesting to compare the randomization and permutation tests' performance with that of the bootstrap test when they are not finitely exact.
In the set-ups I, II, and IV the bootstrap test is more liberal than both randomization tests.
For $\chi_5^2$-distributed data (set-up III), the bootstrap tests behave similarly to the randomization tests, if not less liberal.
For $t_5$-distributed data (set-up II), the permutation tests are even more liberal than the bootstrap tests and for the mixture distribution (set-up IV) they are only slightly less liberal.
In the perhaps most interesting case of the mixture distribution, where none of the tests is finitely exact, it is seen that the mirroring-based tests are most accurate, the rotation-based tests are accurate as well but somewhat conservative, and the permutation tests are too liberal.

\begin{figure}[H]
\includegraphics[width=0.45\textwidth]{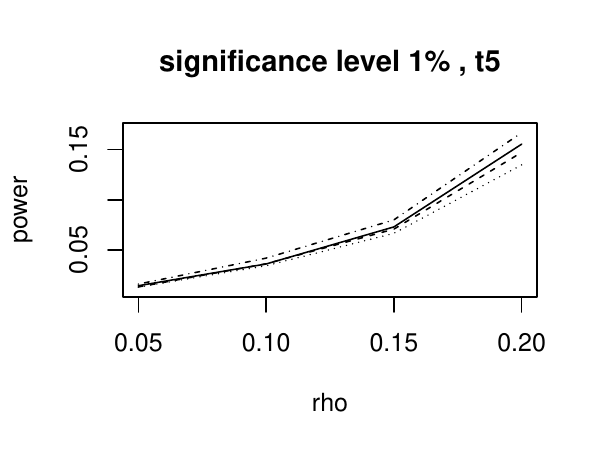}
 \includegraphics[width=0.45\textwidth]{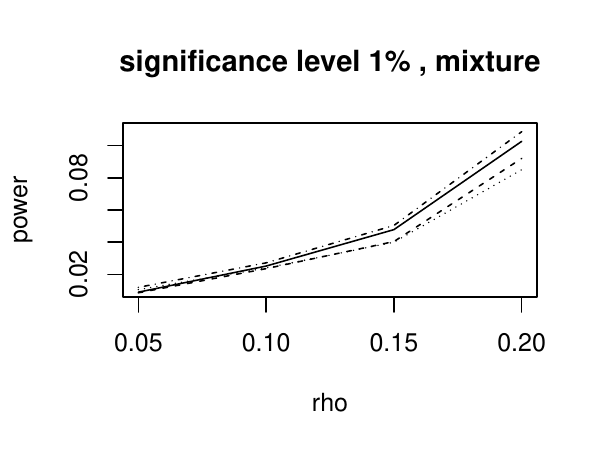}
 \\[-0.2cm]
 \includegraphics[width=0.45\textwidth]{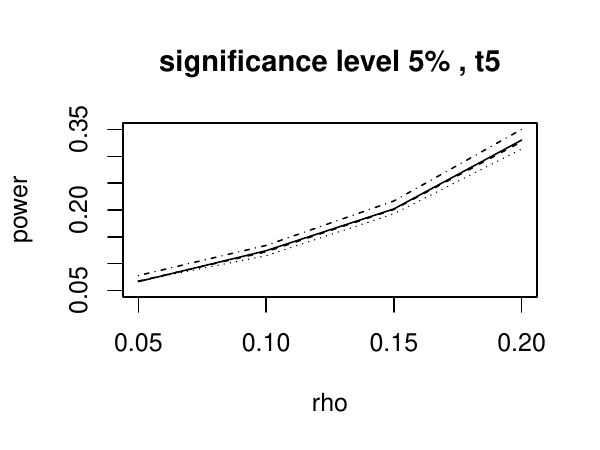}
 \includegraphics[width=0.45\textwidth]{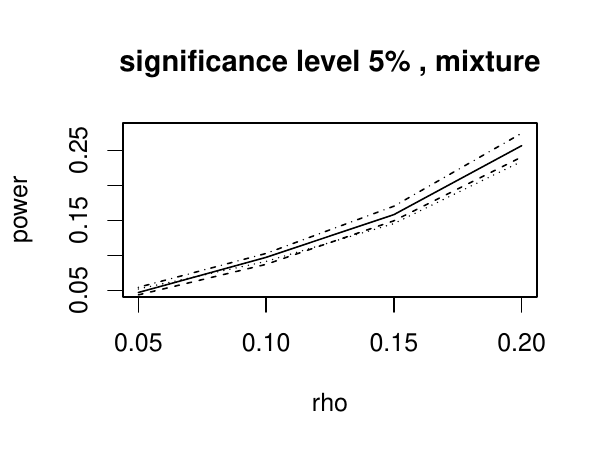}
 \\[-0.2cm]
 \includegraphics[width=0.45\textwidth]{Pearson_Power_t5_alpha5_trans.pdf}
 \includegraphics[width=0.45\textwidth]{Pearson_Power_mix_alpha5_trans.pdf}
 \caption{Simulated power of the correlation test with $t_5$- (left),  $\chi^2_5$-distributed data (middle) and data from a bivariate $t_5$-$\chi^2_5$ mixture distribution (right); based on mirroring (---), rotating (--~--), bootstrap ($\cdots$), permutation ($\cdot$ -- $\cdot$). 
 The sample size is $n=100$.}
 \label{fig:pearsonP1}
\end{figure}

Next, we simulated the power of the correlation tests with true correlation values $\rho= .05, .1, .15,$ $.2$ and sample size $n=100$.
In the multivariate normal case, all tests showed a very similar performance which is why we do not display these results.
In the case of a multivariate $t_5$-distribution (Figure~\ref{fig:pearsonP1}, plots on the left) the permutation tests had the greatest power which is certainly due to its liberality under the null hypothesis.
The case of a multivariate $\chi^2_5$-distribution was realized by generating $X+Z$ and $Y+Z$ where $X,Y,Z$ are independent $\Gamma$-distributed random variables with scale parameters 2 and shape parameters $2.5 \cdot (1 -\rho)$ for $X$ and $Y$ and $2.5 \rho$ for $Z$.
Because most tests were more or less liberal under the null hypothesis under the $\chi^2_5$-setting for $n=100$, we increased the sample size to $n=200$.
There, however, hardly any difference is seen between the power of the tests (results not shown).
In the mixture case (Figure~\ref{fig:pearsonP1}, plots on the right) the permutation test was slightly more powerful than the others.
But again, the permutation test was quite liberal under the null hypothesis.

We conclude that, if the tests are reliable under the null hypothesis, the differences in power are not very large.
Thus, we cannot give a clear recommendation for the choice of test based on the power study.
However, the random permutation approach of course often enjoys the experimental justification that it reflects the situation that is created by random assignments, that is, independence.

\subsection{Mann-Whitney effect test: simulation results for underlying Gumbel-Hougaard and independence copulae}
\label{sec:MW_2}
\begin{figure}[H]
 \includegraphics[width=0.33\textwidth]{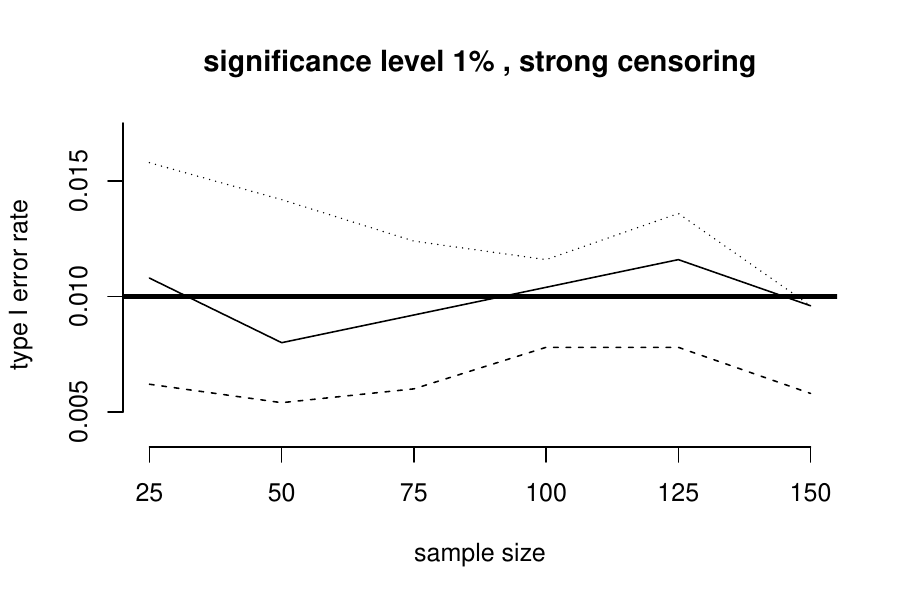}
 \includegraphics[width=0.33\textwidth]{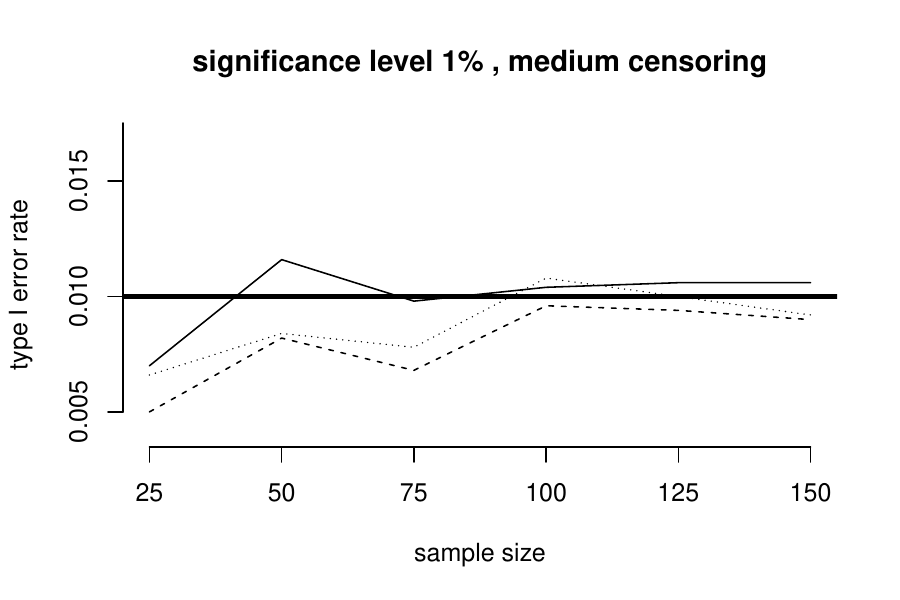}
 \includegraphics[width=0.33\textwidth]{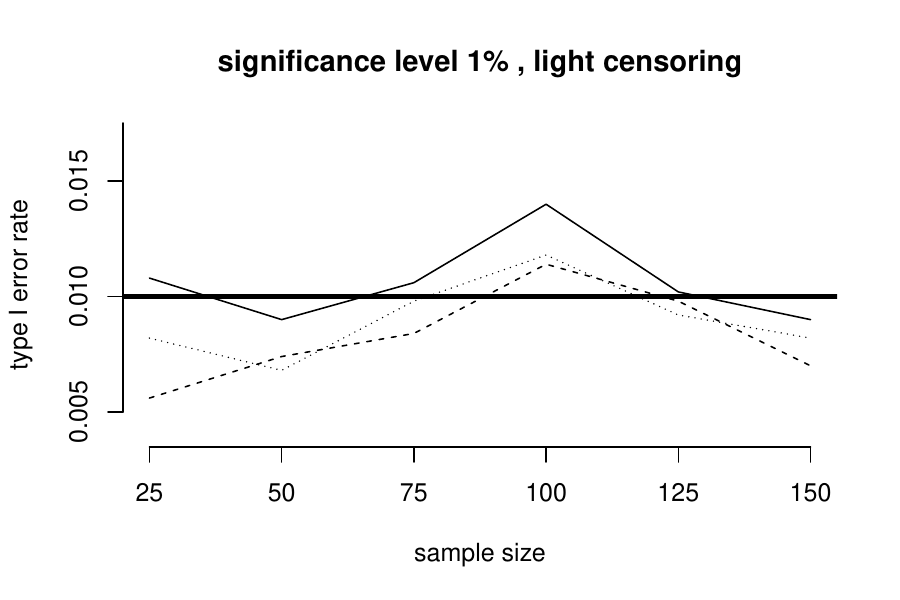} \\[-0.4cm]
 \includegraphics[width=0.33\textwidth]{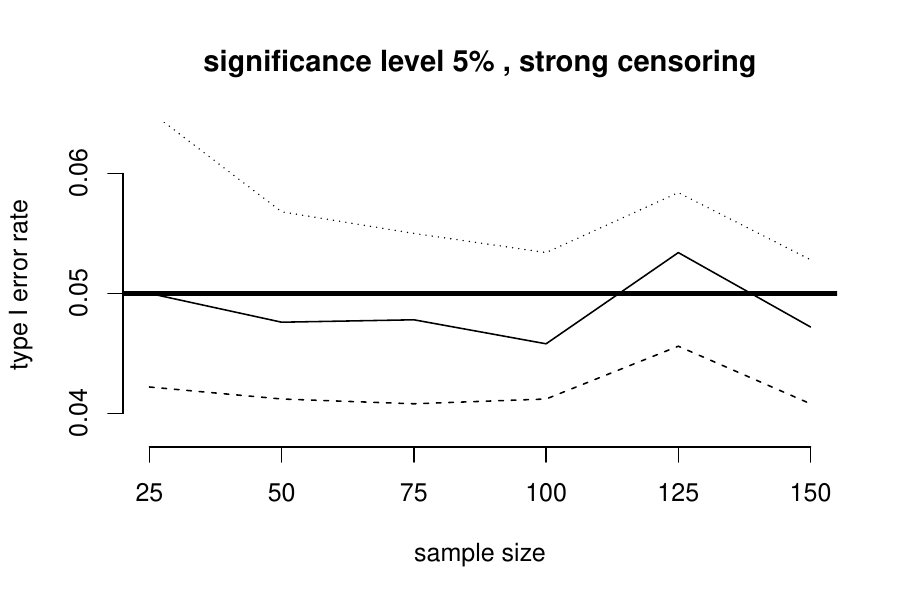}
 \includegraphics[width=0.33\textwidth]{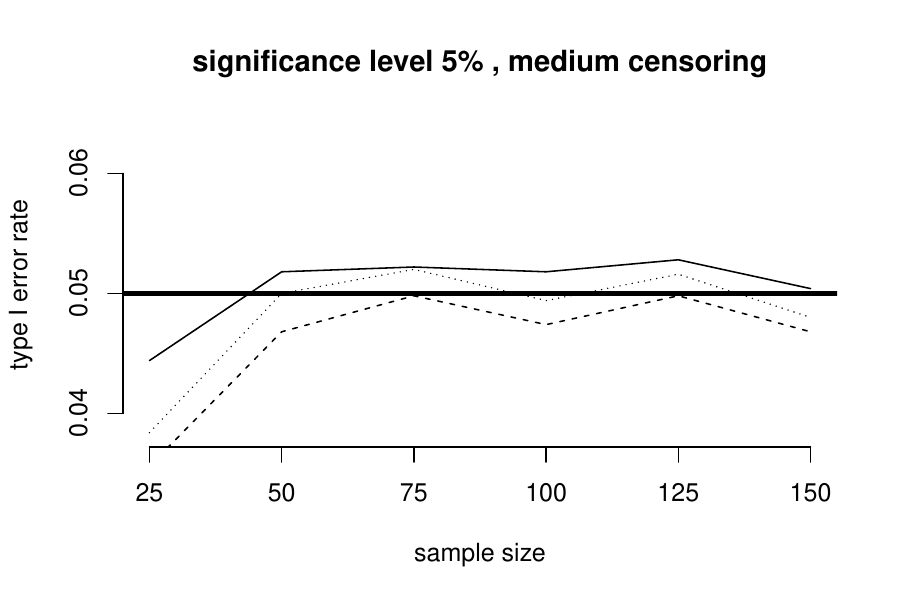}
 \includegraphics[width=0.33\textwidth]{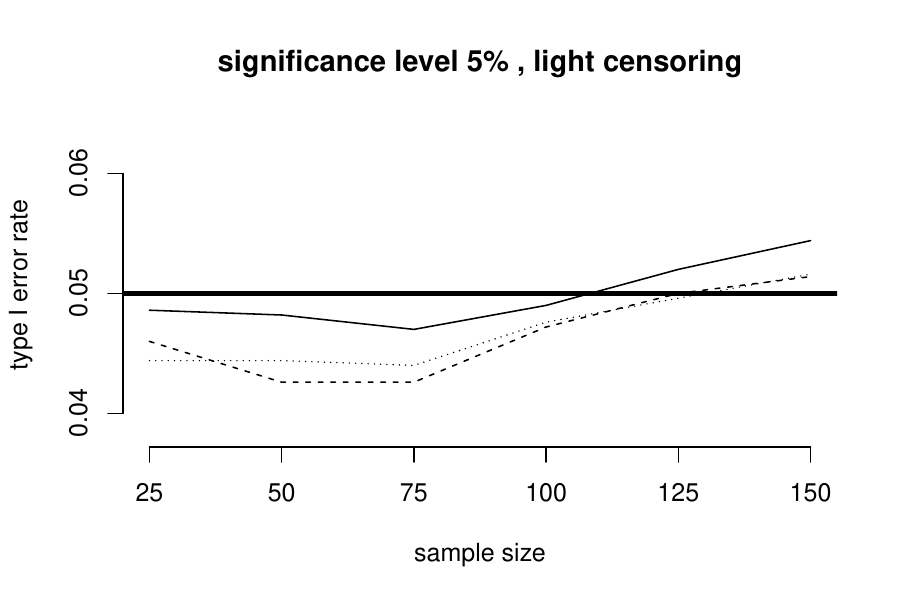} \\[-0.4cm]
 \includegraphics[width=0.33\textwidth]{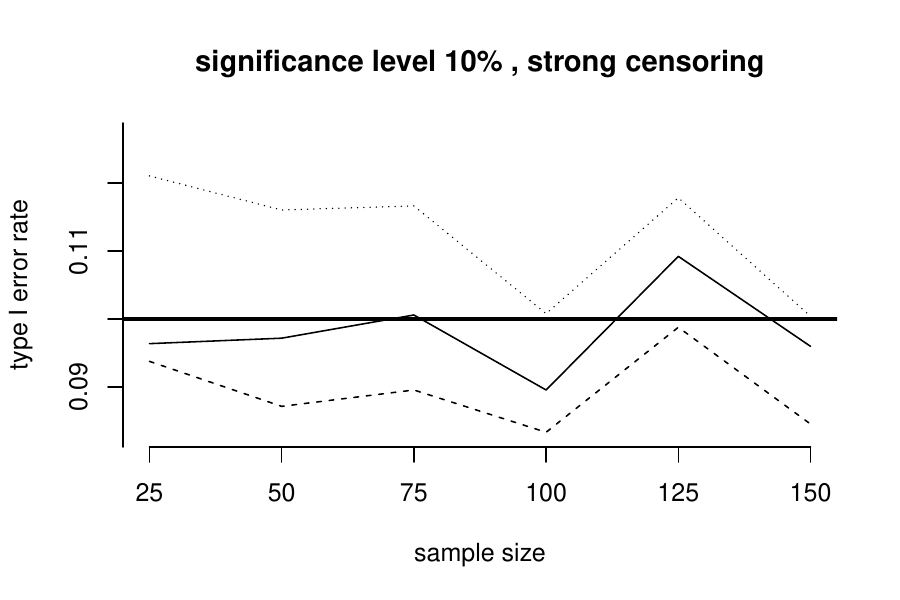}
 \includegraphics[width=0.33\textwidth]{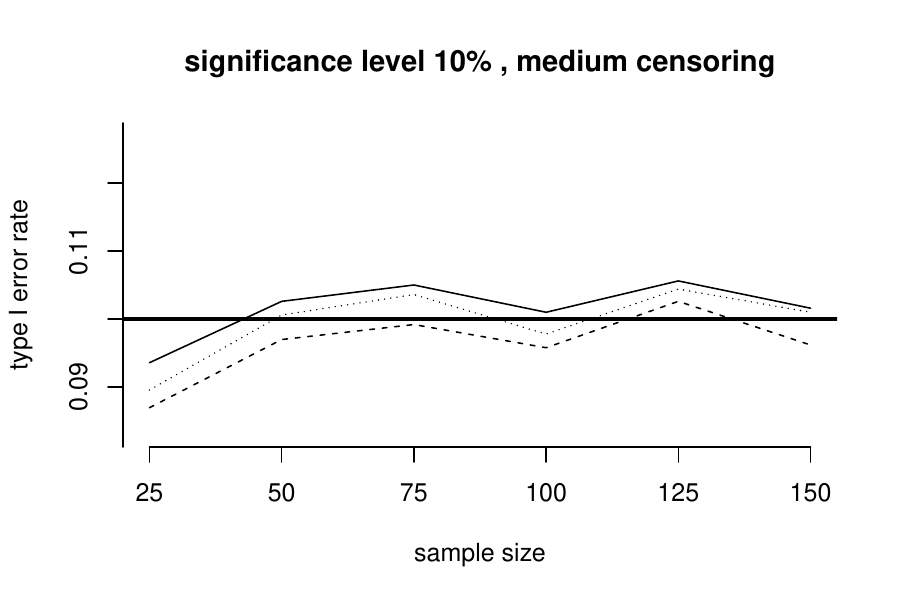}
 \includegraphics[width=0.33\textwidth]{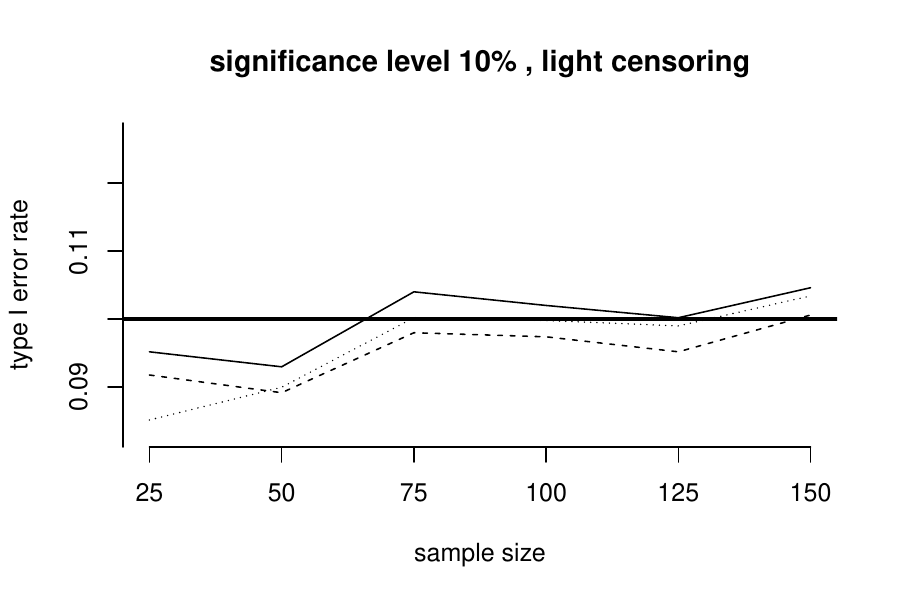} \\[-0.4cm]
 \includegraphics[width=0.33\textwidth]{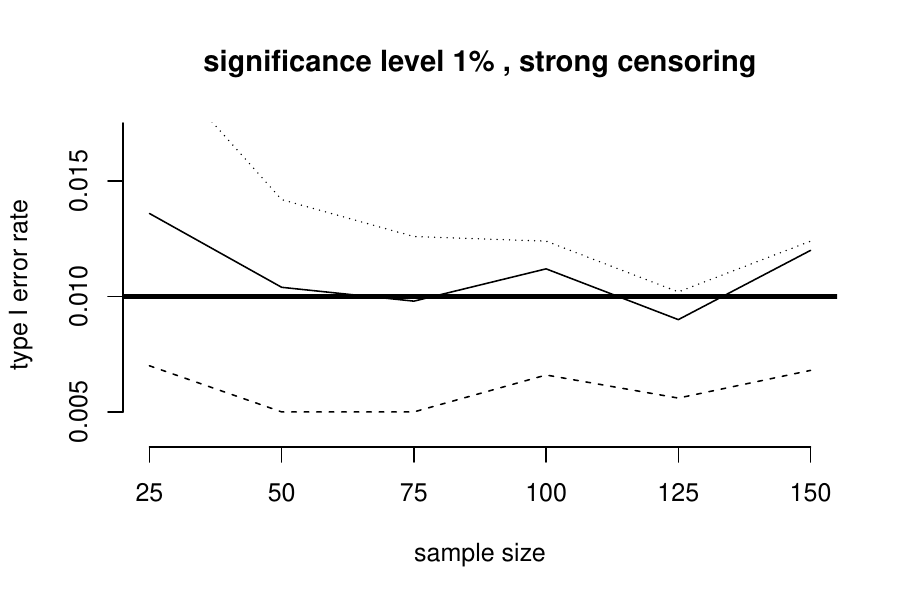}
 \includegraphics[width=0.33\textwidth]{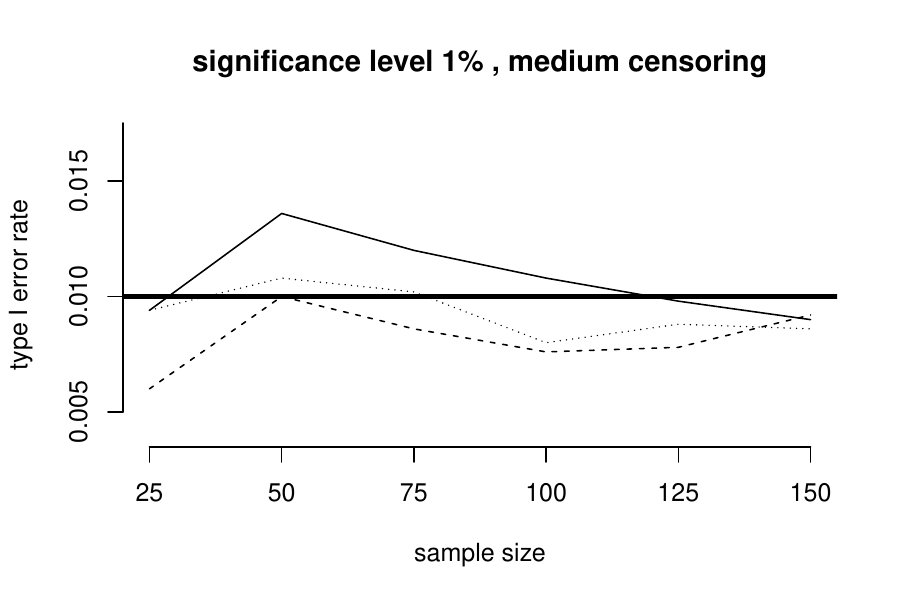}
 \includegraphics[width=0.33\textwidth]{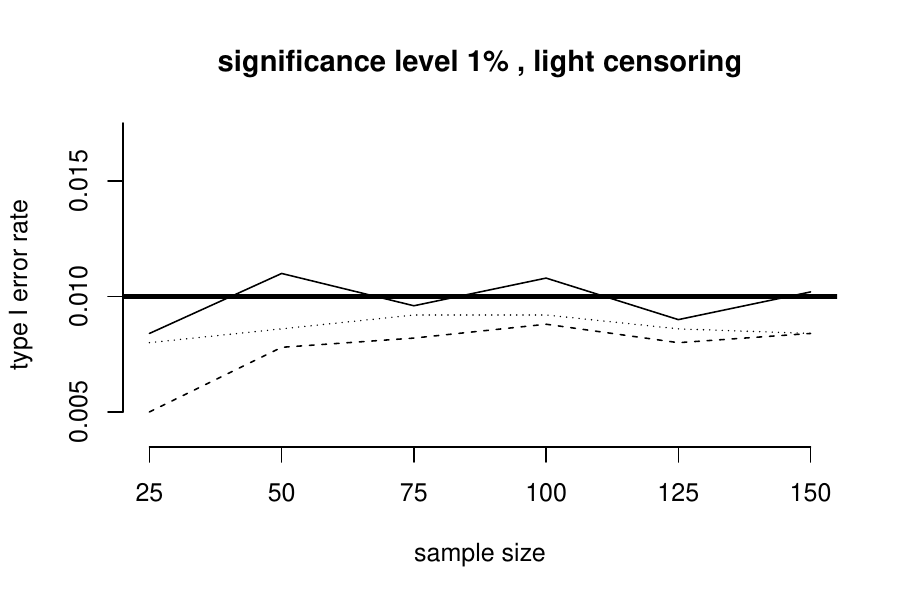} \\[-0.4cm]
 \includegraphics[width=0.33\textwidth]{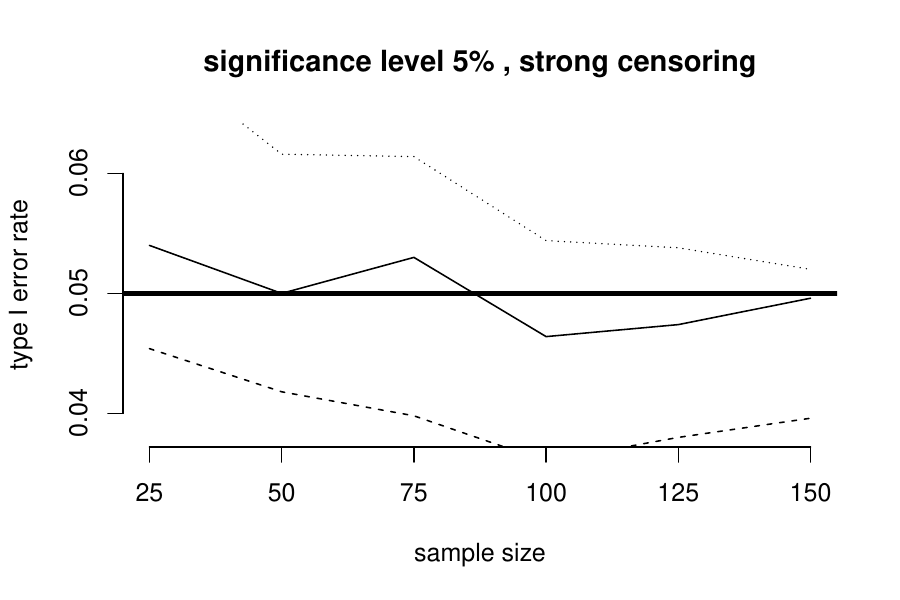}
 \includegraphics[width=0.33\textwidth]{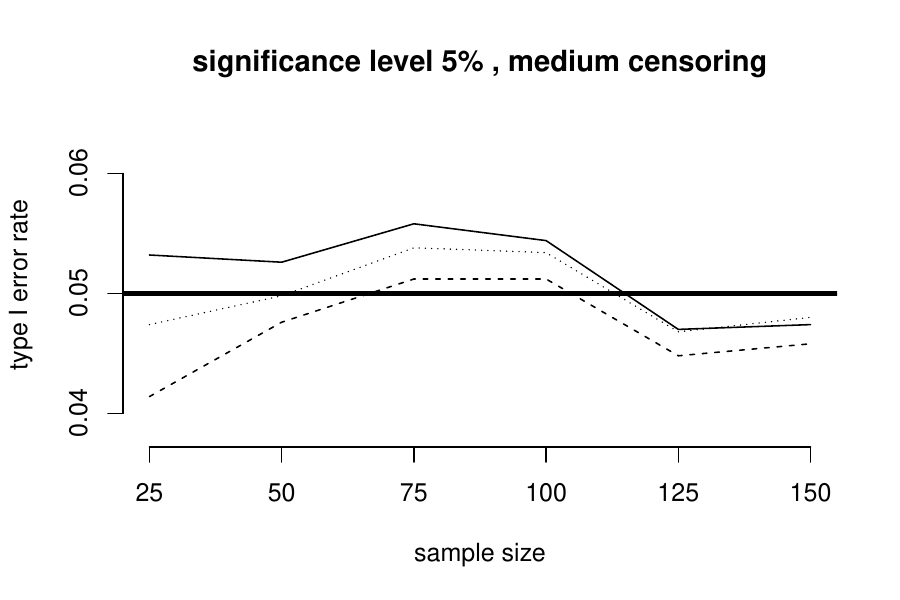}
 \includegraphics[width=0.33\textwidth]{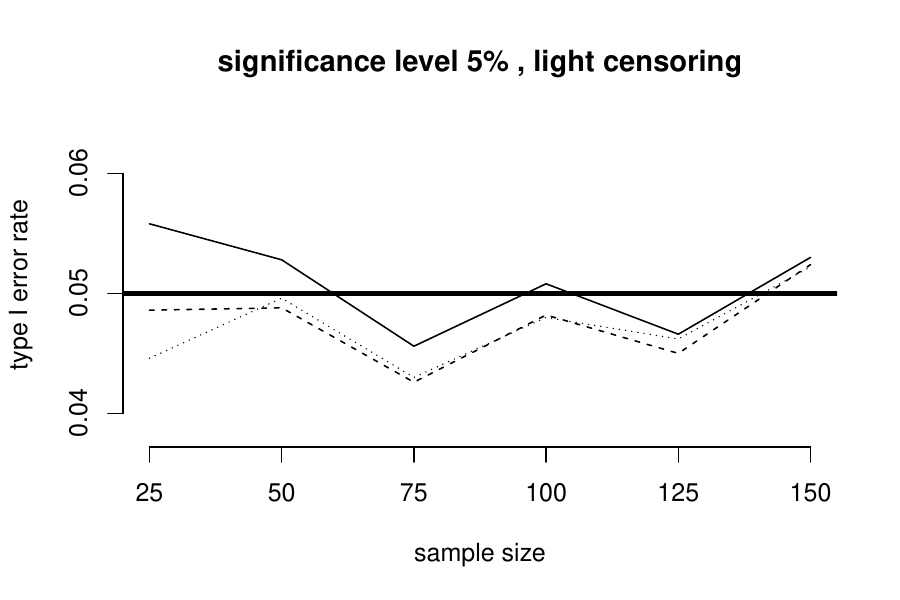} \\[-0.4cm]
 \includegraphics[width=0.33\textwidth]{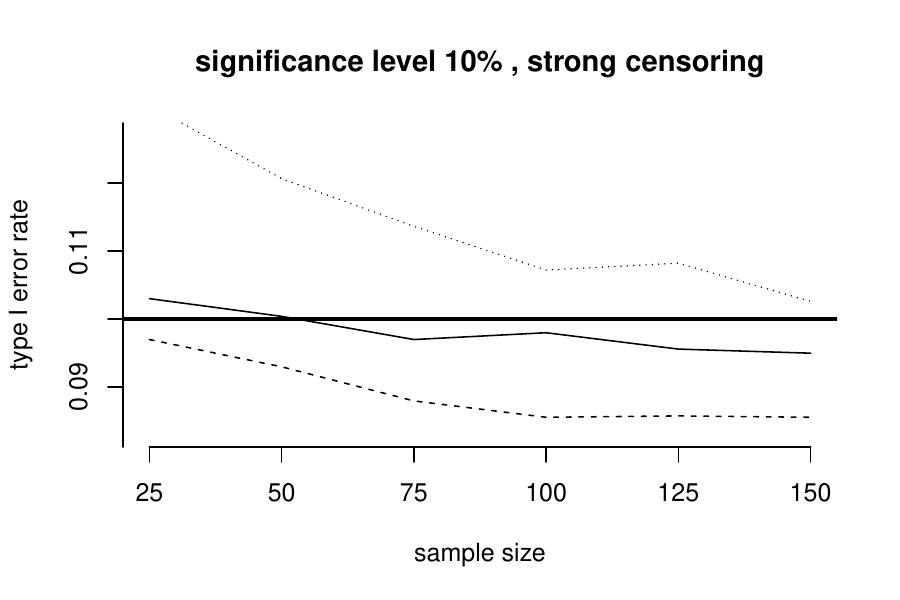}
 \includegraphics[width=0.33\textwidth]{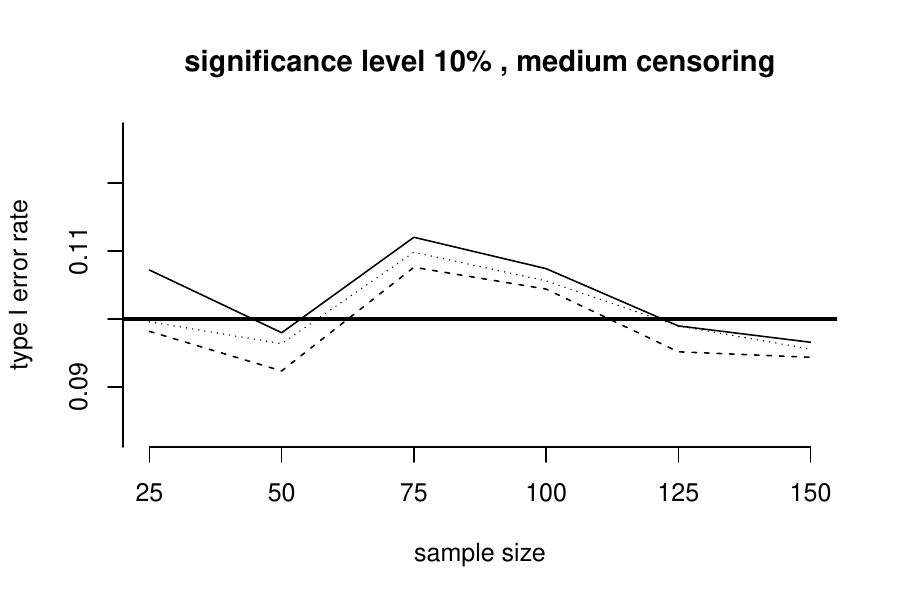}
 \includegraphics[width=0.33\textwidth]{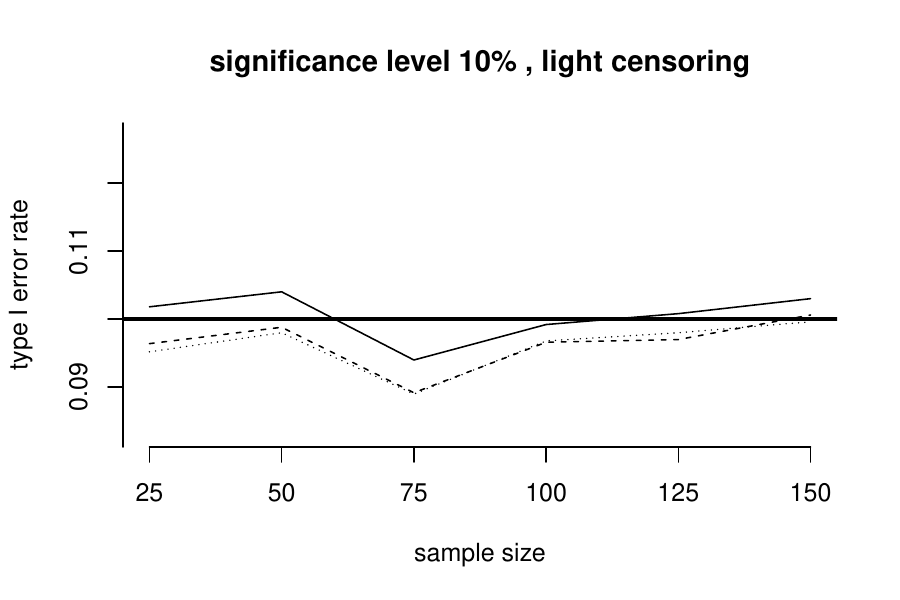} 
 \caption{Simulated type I error rates of the Mann-Whitney-type tests with the Gumbel-Hougaard copula underlying the data under strong (left), medium (middle), and light censoring (right); equal (upper half) and unequal marginal survival distributions (lower half); based on randomization (---), bootstrap (- -), normal quantiles ($\cdots$). 
 The horizonal line is the nominal significance level.}
 \label{fig:mw1}
\end{figure}

\begin{figure}[H]
 \includegraphics[width=0.33\textwidth]{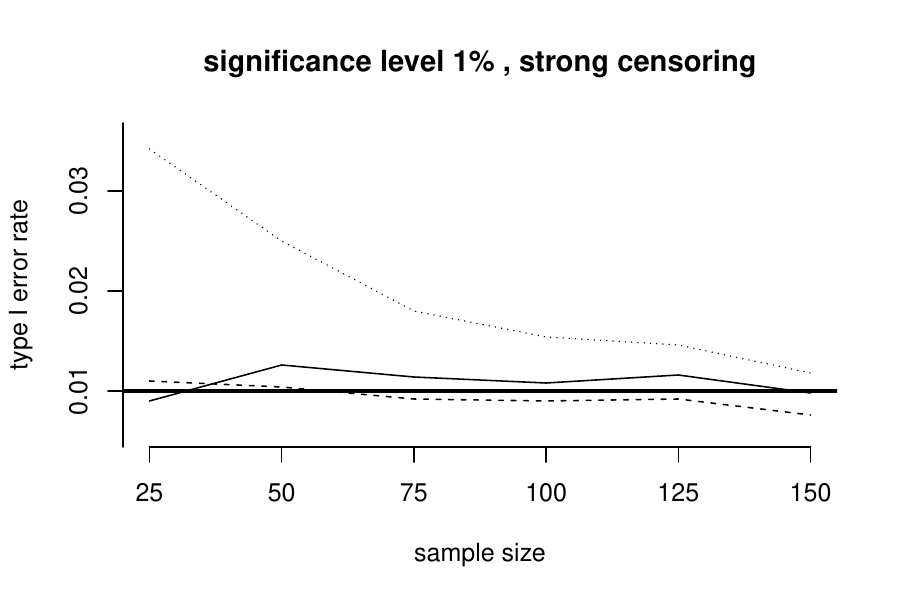}
 \includegraphics[width=0.33\textwidth]{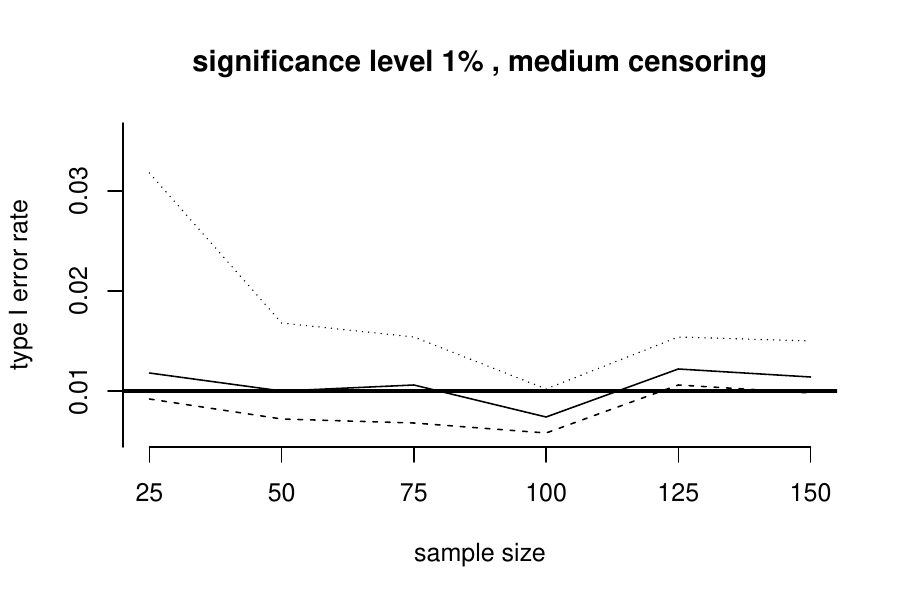}
 \includegraphics[width=0.33\textwidth]{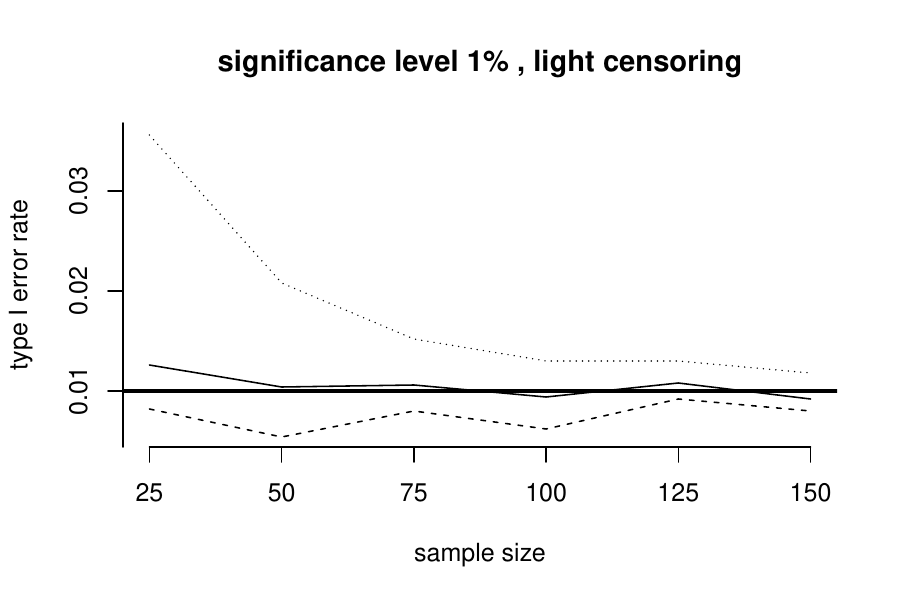} \\[-0.4cm]
 \includegraphics[width=0.33\textwidth]{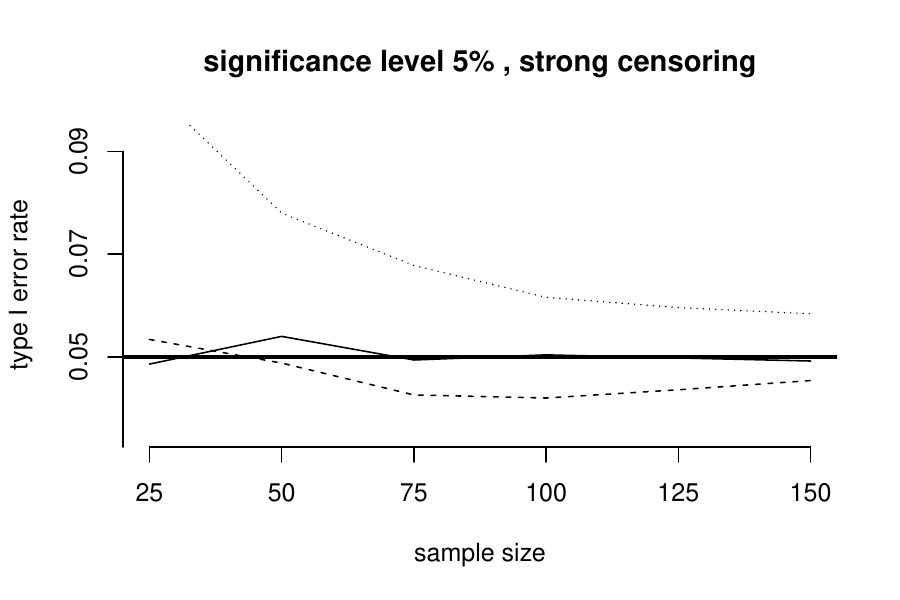}
 \includegraphics[width=0.33\textwidth]{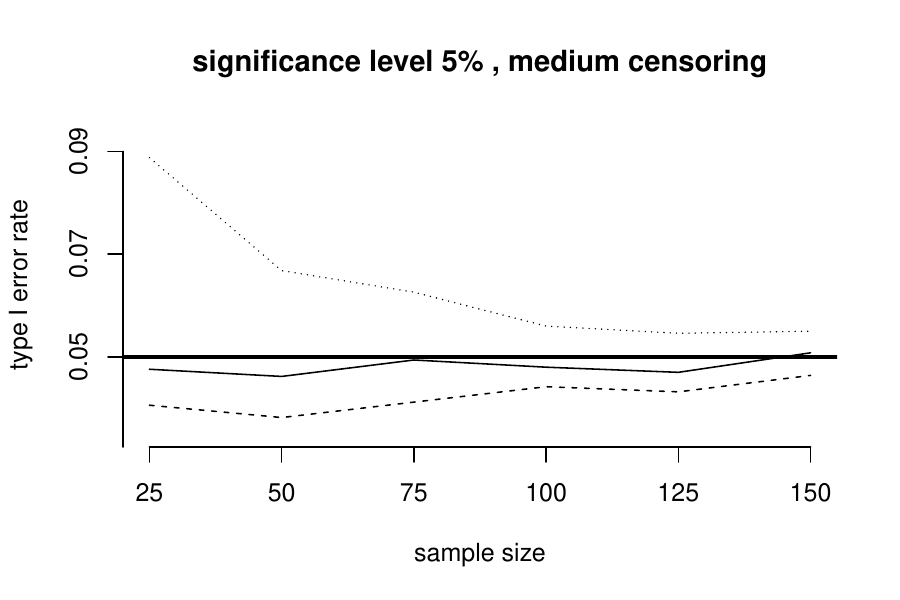}
 \includegraphics[width=0.33\textwidth]{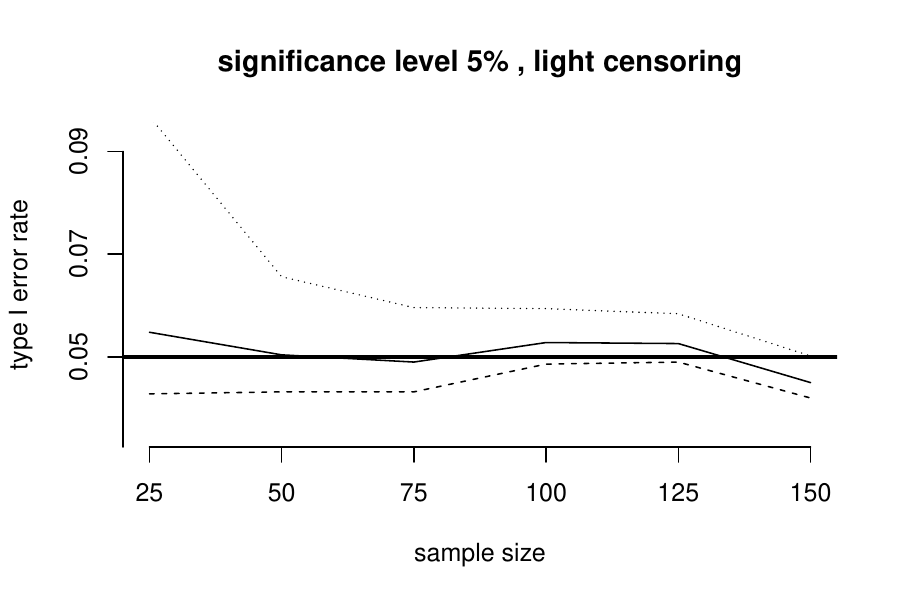} \\[-0.4cm]
 \includegraphics[width=0.33\textwidth]{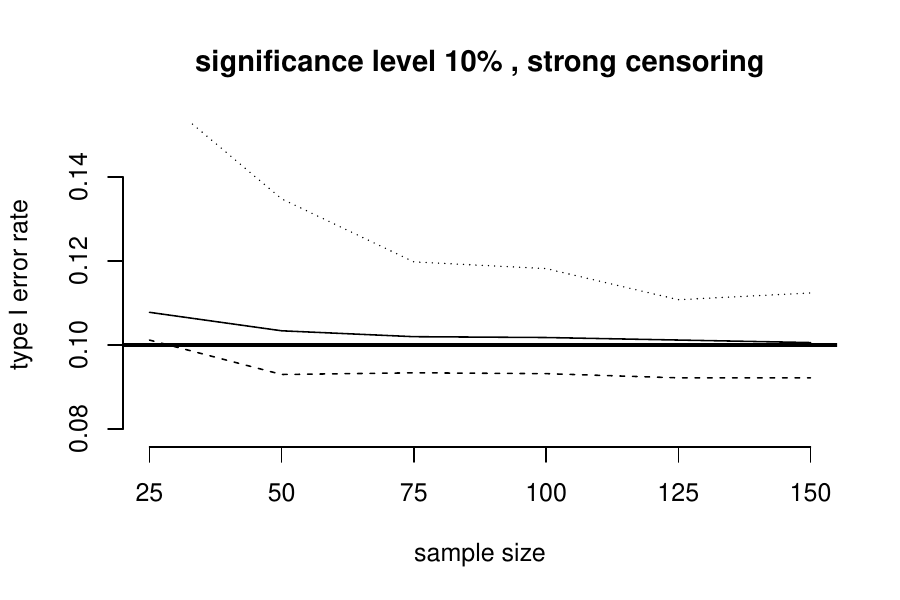}
 \includegraphics[width=0.33\textwidth]{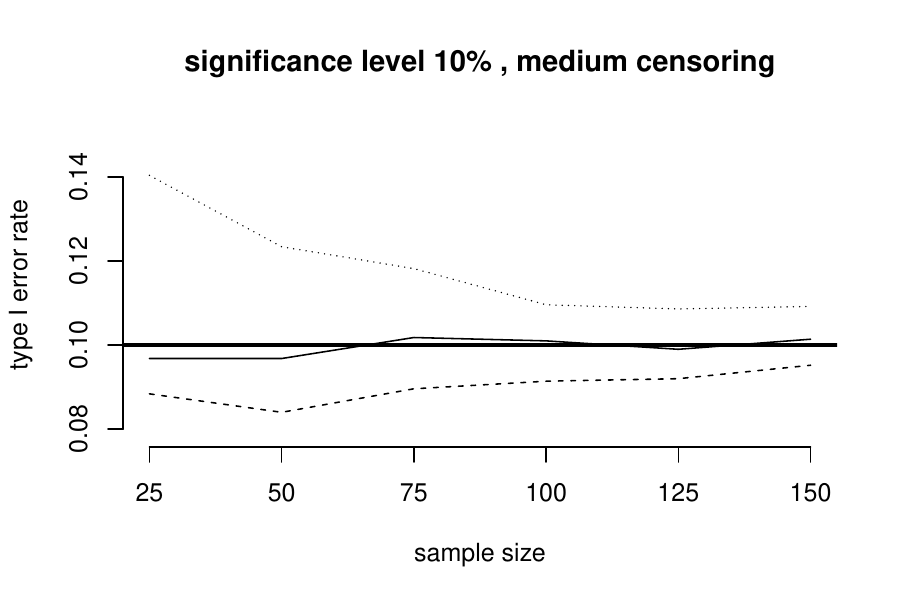}
 \includegraphics[width=0.33\textwidth]{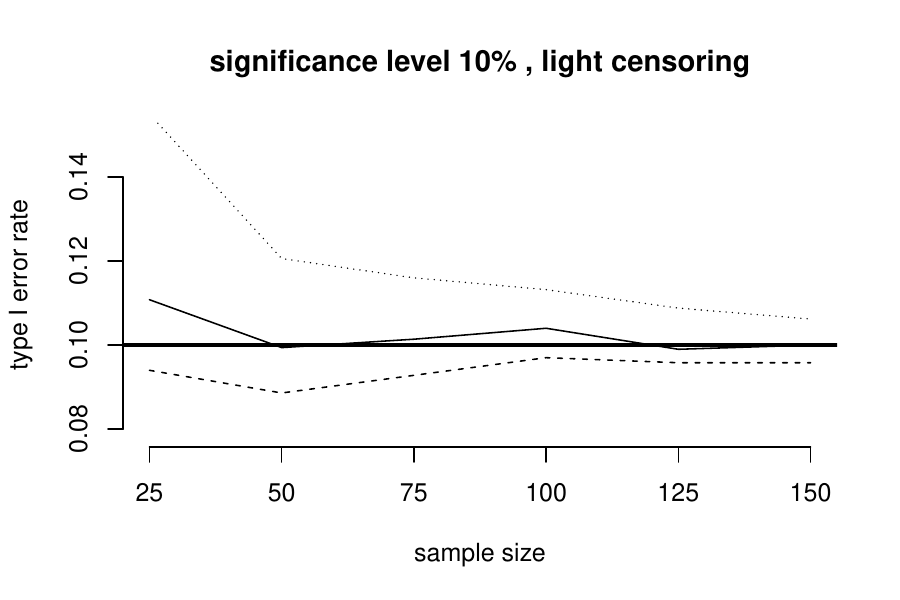} \\[-0.4cm]
 \includegraphics[width=0.33\textwidth]{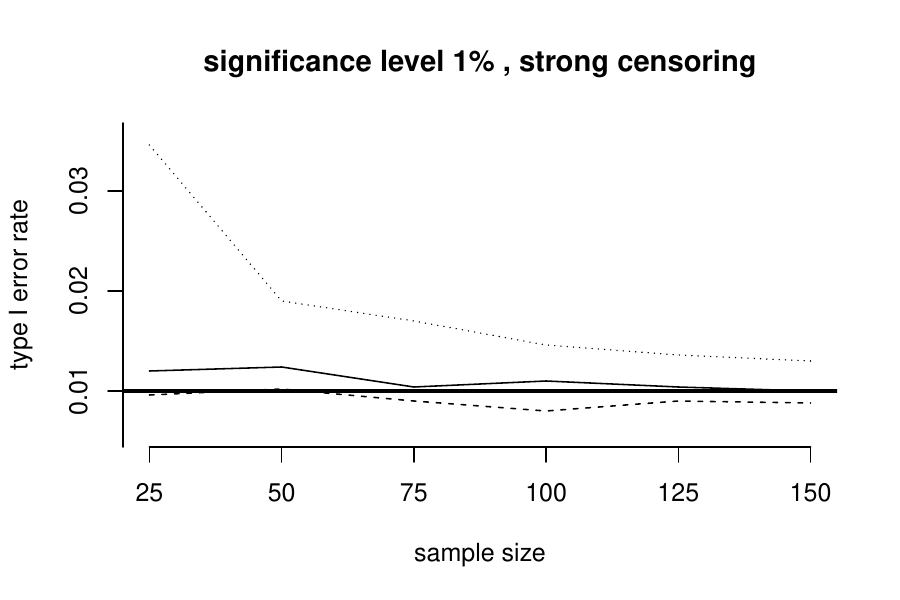}
 \includegraphics[width=0.33\textwidth]{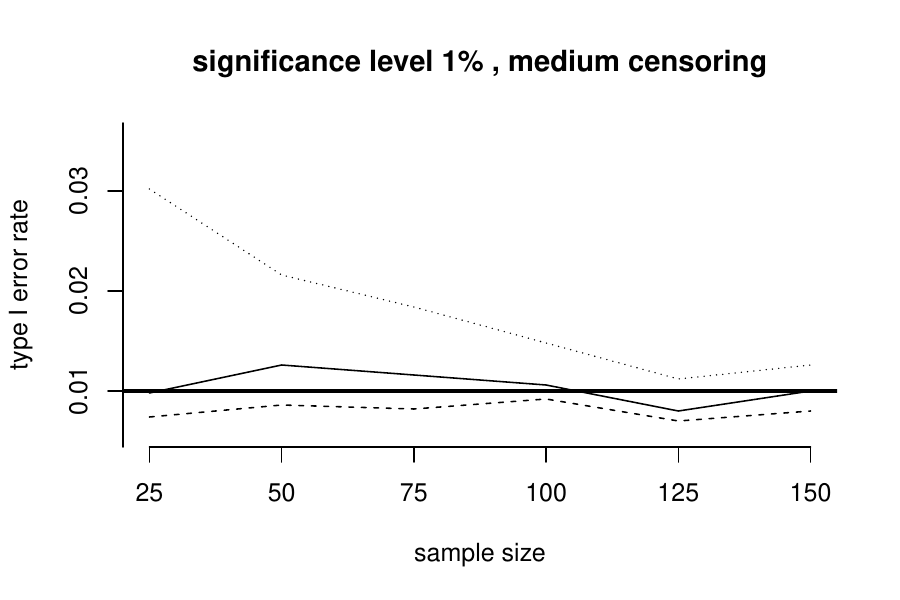}
 \includegraphics[width=0.33\textwidth]{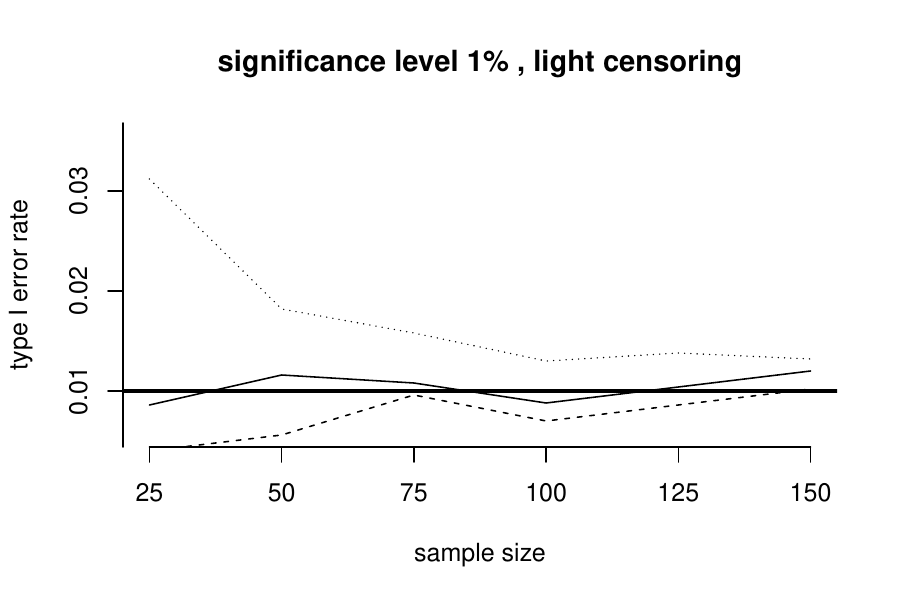} \\[-0.4cm]
 \includegraphics[width=0.33\textwidth]{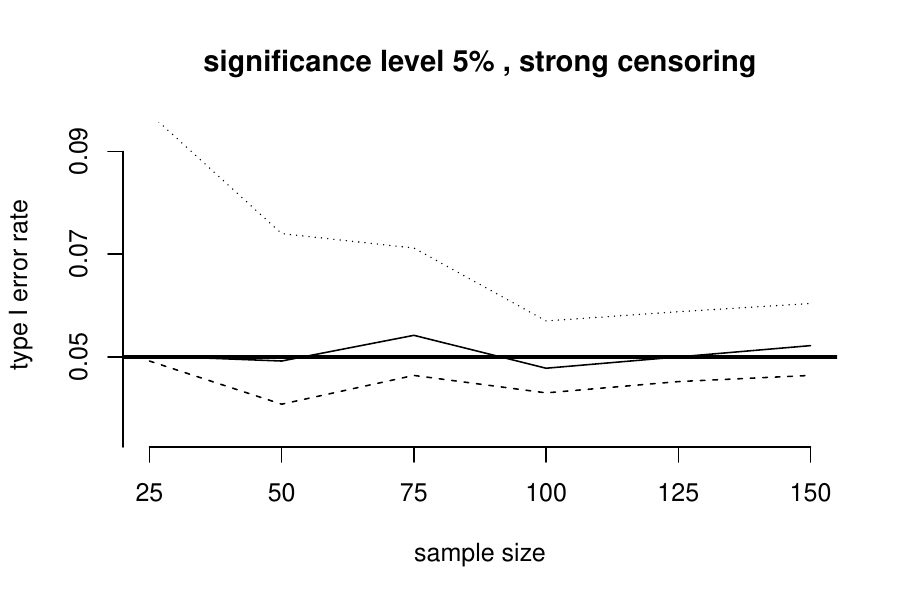}
 \includegraphics[width=0.33\textwidth]{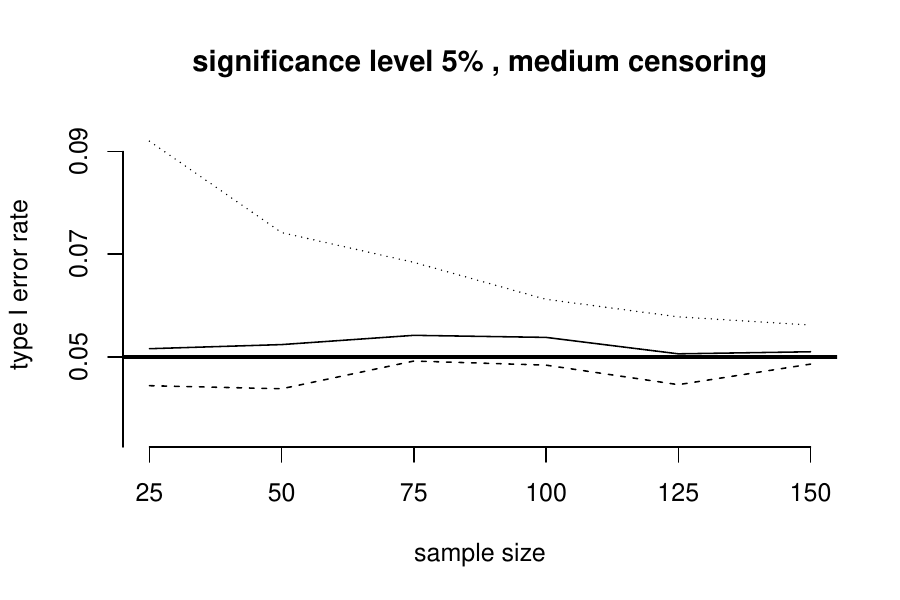}
 \includegraphics[width=0.33\textwidth]{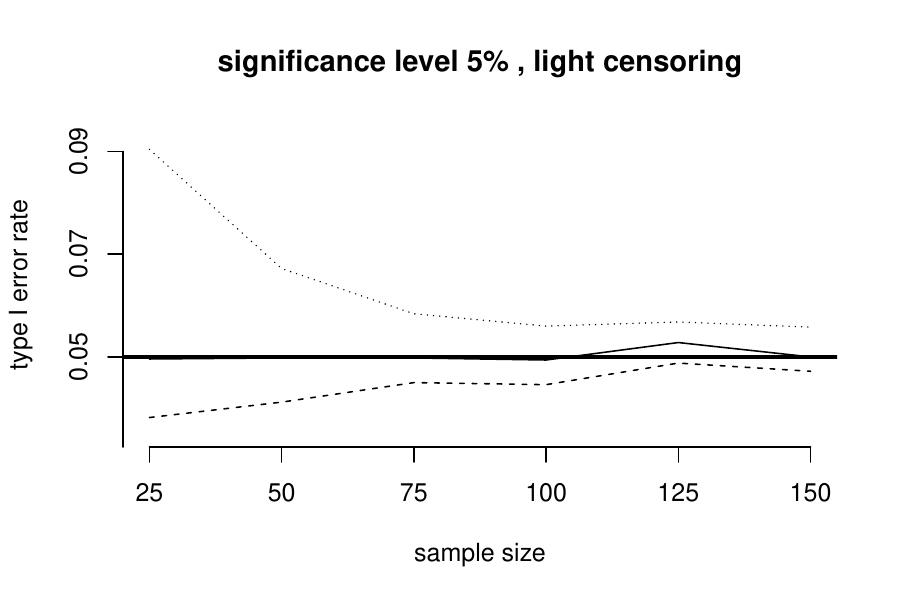} \\[-0.4cm]
 \includegraphics[width=0.33\textwidth]{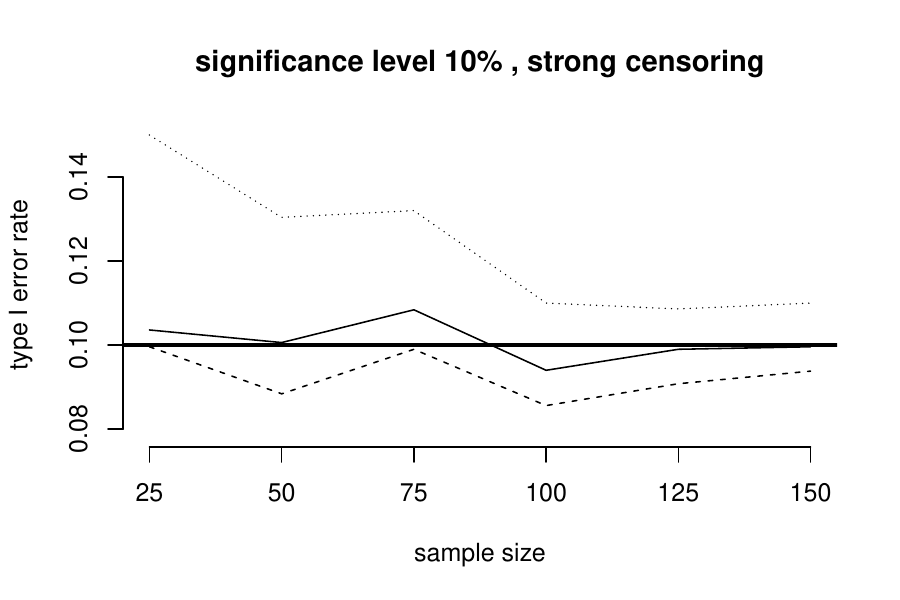}
 \includegraphics[width=0.33\textwidth]{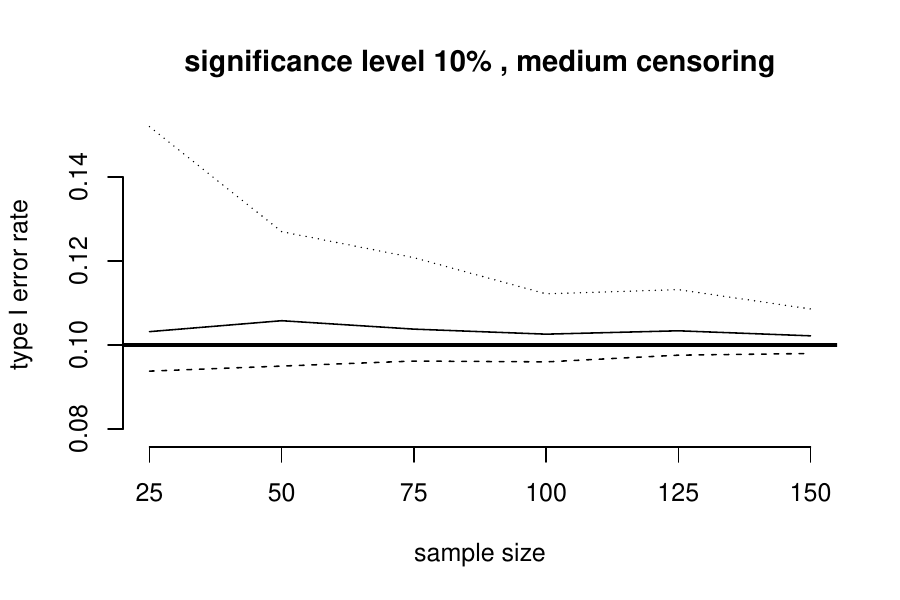}
 \includegraphics[width=0.33\textwidth]{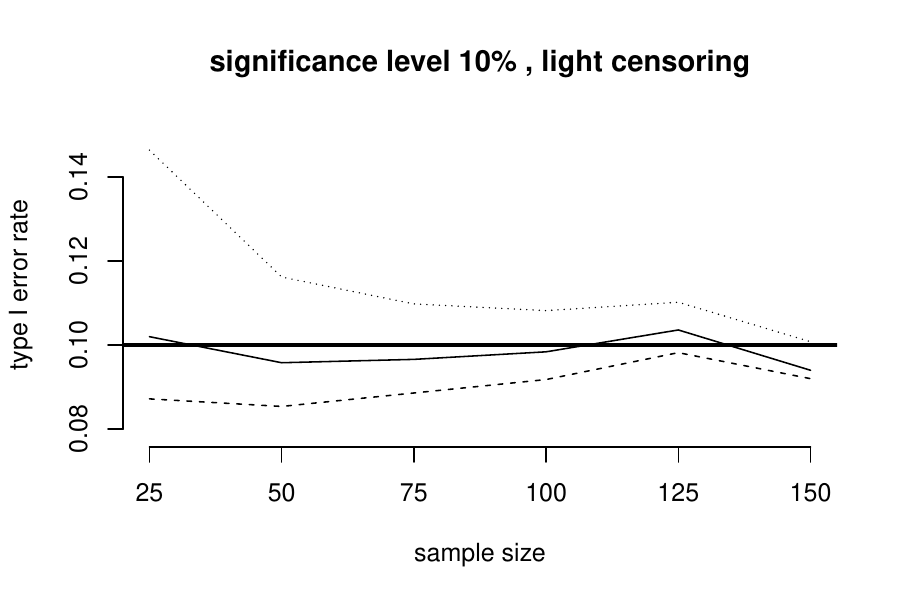} 
 \caption{Simulated type I error rates of the Mann-Whitney-type tests with the independence copula underlying the data under strong (left), medium (middle), and light censoring (right); equal (upper half) and unequal marginal survival distributions (lower half); based on randomization (---), bootstrap (- -), normal quantiles ($\cdots$). 
 The horizonal line is the nominal significance level.}
 \label{fig:mw3}
\end{figure}
{\phantom{X}}

\subsection{Mann-Whitney effect test: simulation results for unequal sample sizes and Exponential versus Gompertz distribution}
\label{sec:MW_uneq_simu}

{\blue In practice, it might not always be possible to include both members of a pair into a study.
For instance, if a serious surgical mistake happened to the first member whereas the second surgery was a success, there is usually no reason to throw away the data on the second member.
The only required assumption is that the surgical mistake had no impact on the second member as well.

If only one member of a pair can be included in a study, unequal sample sizes are the result.
Denote by $n_1$ and $n_2$ the amount of eligible first and second members of the pairs, respectively.
As motivated above, these numbers might in practice be random while the planned sample size $n \geq n_1, n_2$, i.e.\ the originally intended number of pairs, is fixed.
For ease of reference, we will refer to the $n_1$ and $n_2$ data points as two samples, even though some of them indeed form pairs and are thus dependent.

The next set of simulations for the case of unequal sample sizes was based on different marginals distribuions:
a Gompertz distribution with shape parameter $\eta=.6$ and scale parameter $b \approx 3.05605$ and an exponential distribution with rate parameter $\lambda=3$, cf.\ Figure~\ref{fig:go_exp} for an illustration.
The terminal time point is $\tau=.6$.
This configuration implies that the Mann-Whitney effect is $p=.5.$

\begin{figure}[!h]
\centering
 \includegraphics[width=0.8\textwidth]{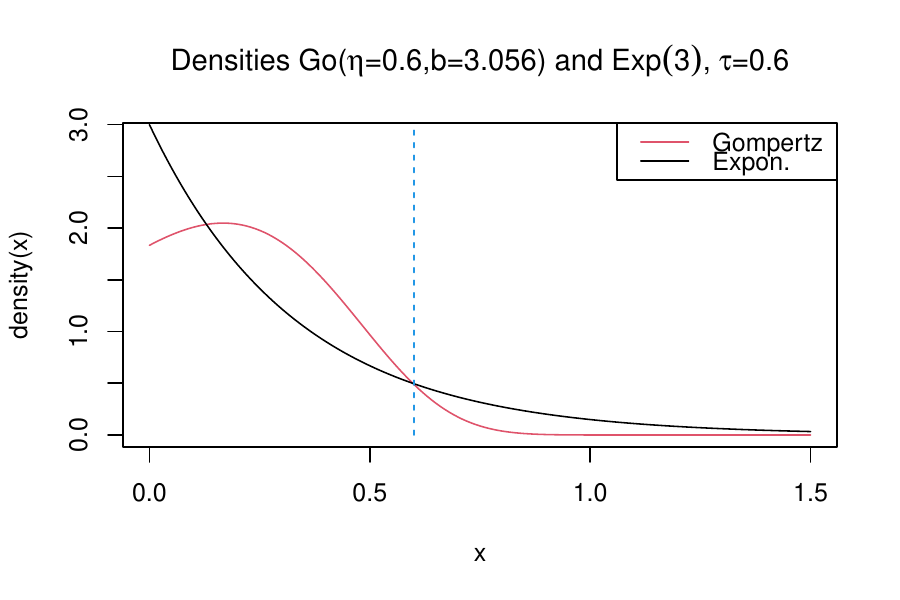}
 \caption{Densities of marginal distributions: Gompertz versus exponential.}
 \label{fig:go_exp}
\end{figure}

Denote by $n_p = n - n_1 - n_2$ the number of paired observations.
We have simulated two cases in which the sample sizes have the following multinomial distributions:\\
\phantom{I}(I) $(n_p, n_1, n_2) \sim mult(n, .7, .1, .2)$,
\\
(II) $(n_p, n_1, n_2) \sim mult(n, .7, .125, .125)$.
\\
The considered right-censoring rates are about 18.2\%/26.8\% (light), 28.5\%/34.4\% (medium), and 44.8\%\%/46.5\% (strong) for the Gompertz/exponential marginal distributions.
Censoring was caused by taking the minimum of $\tau$ and generated uniformly distributed random variables with minimum parameter 0 and maximum parameters 1.75, 1, and .6, respectively.
The remaining simulation parameters, i.e.\ $n, \alpha$, and the number of Monte-Carlo iterations are as in Section~\ref{sec:simus} in the main manuscript.

The simulation results are not much different from the previous ones:
the asymptotic test is always too liberal but one can clearly see the convergence of the rejection rates to the significance level as $n$ approaches the value 150.
In most cases, the randomization test is closer to $\alpha$, compared to the bootstrap test.
One notable exception is the case in which the Gumbel-Hougaard copula meets the strong censoring rate. Here, the randomization test shows a tendency to be conservative.
All in all, however, the randomization-based test performs best.}

\begin{figure}
 \includegraphics[width=0.34\textwidth]{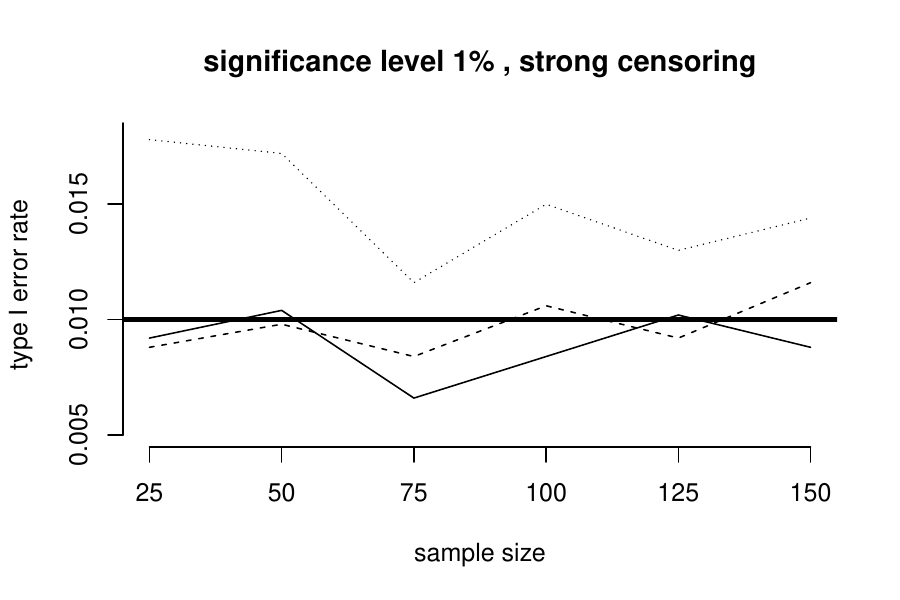}
 \hspace{-0.5cm}
 \includegraphics[width=0.34\textwidth]{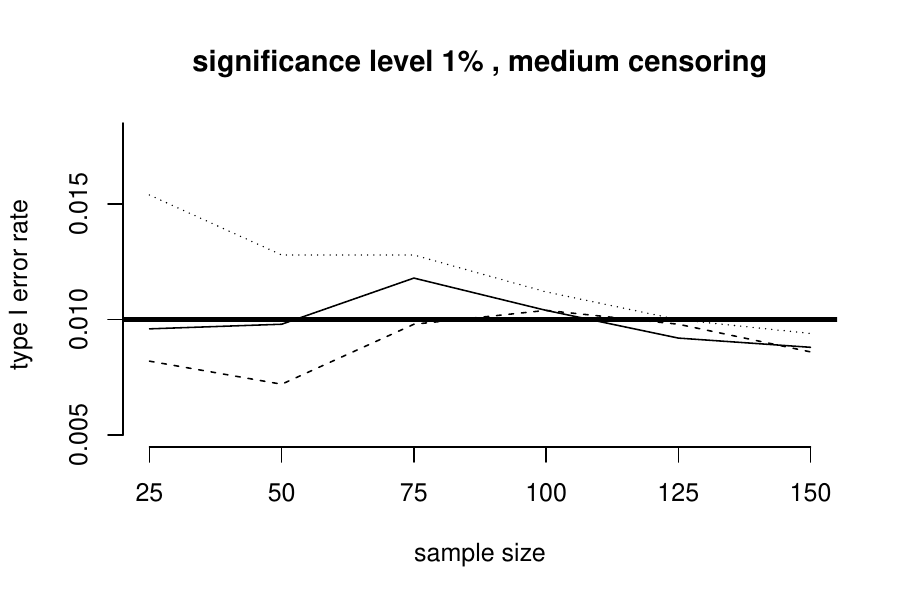}\hspace{-0.5cm}
 \includegraphics[width=0.34\textwidth]{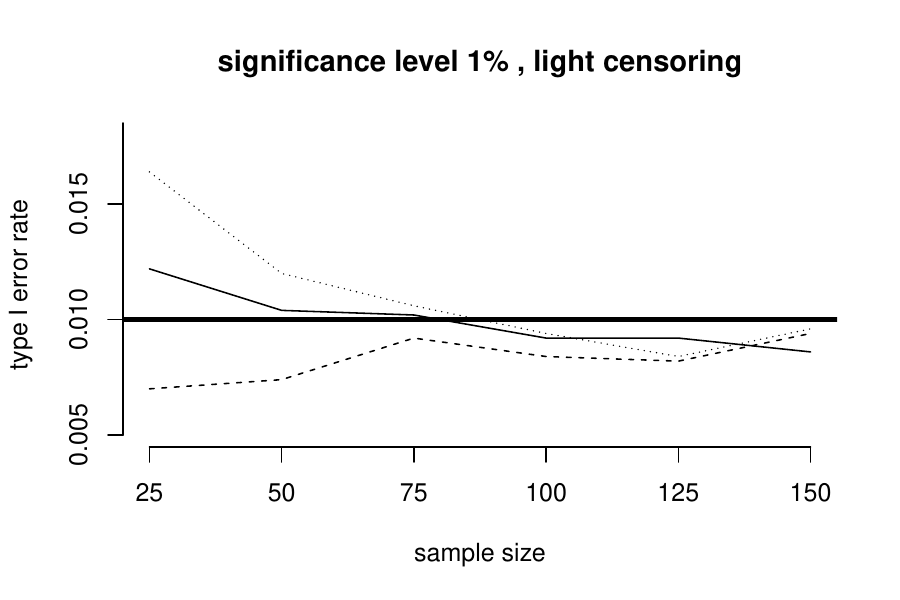} \\[-0.4cm]
 \includegraphics[width=0.34\textwidth]{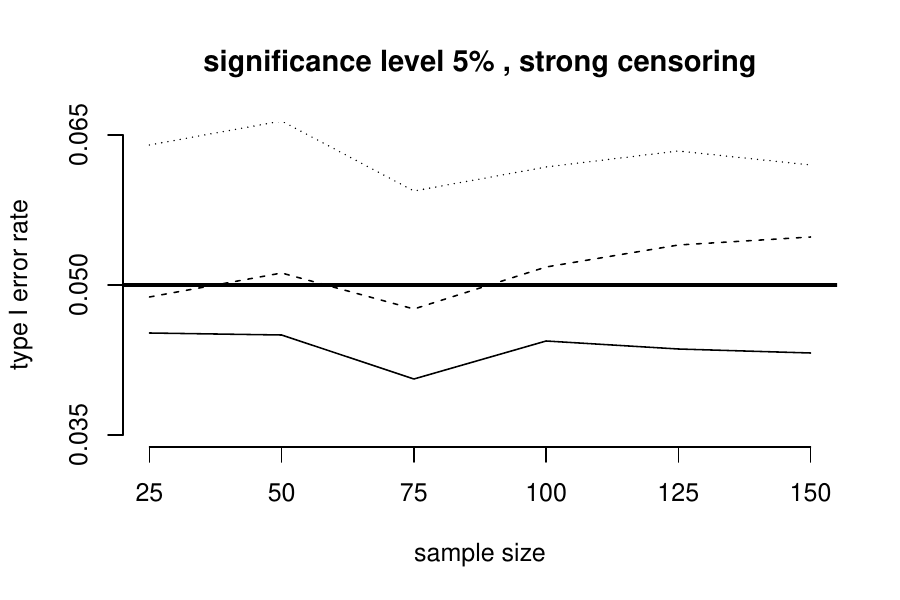}
 \hspace{-0.5cm}
 \includegraphics[width=0.34\textwidth]{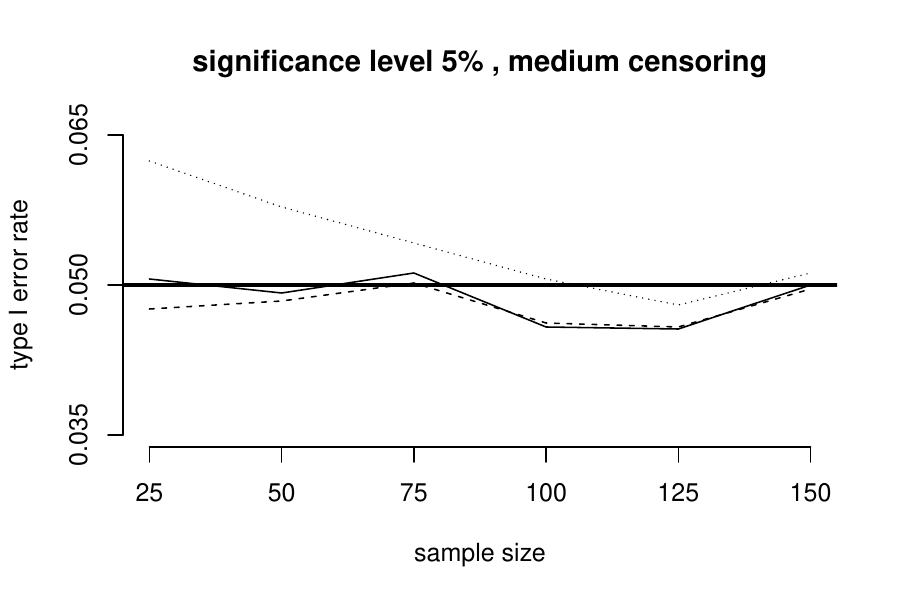}\hspace{-0.5cm}
 \includegraphics[width=0.34\textwidth]{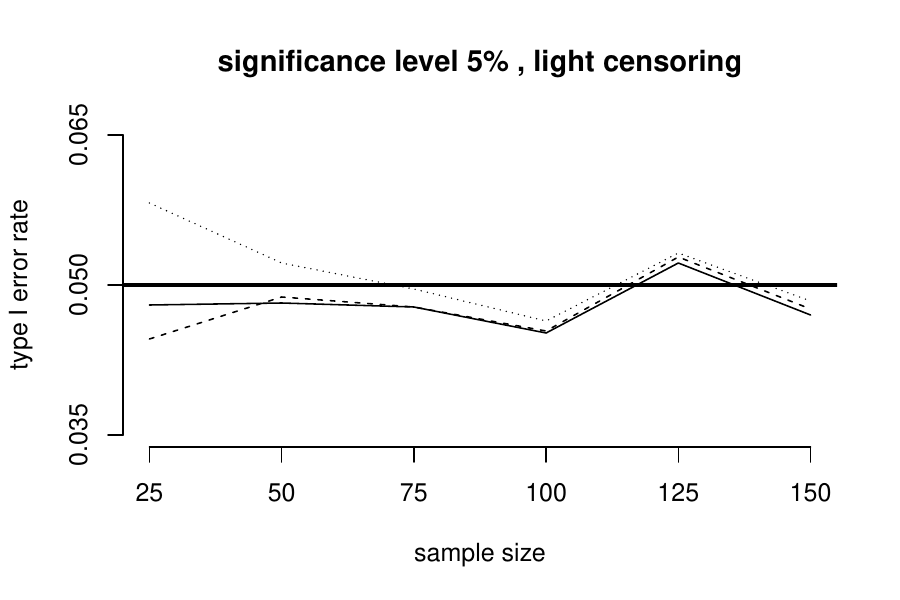} \\[-0.4cm]
 \includegraphics[width=0.34\textwidth]{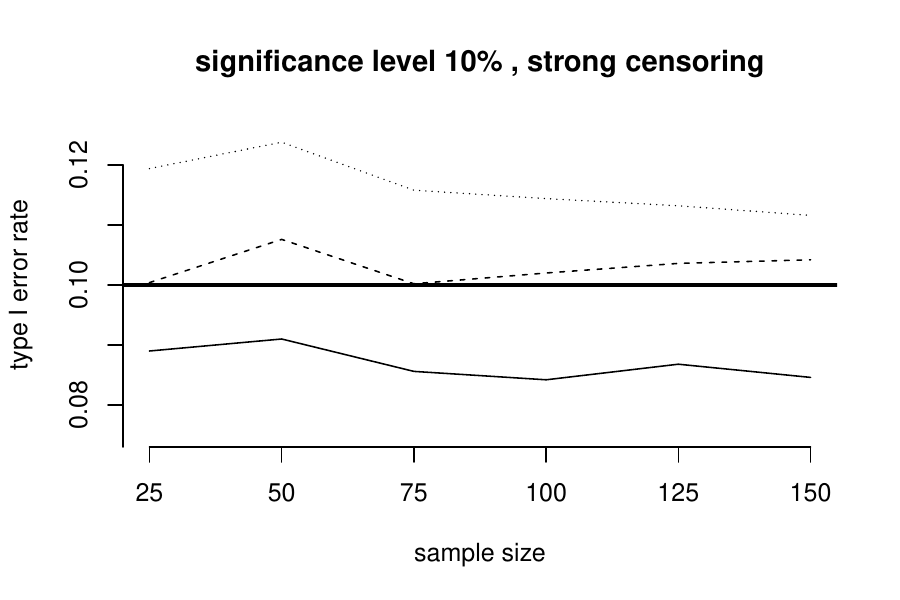}
 \hspace{-0.5cm}
 \includegraphics[width=0.34\textwidth]{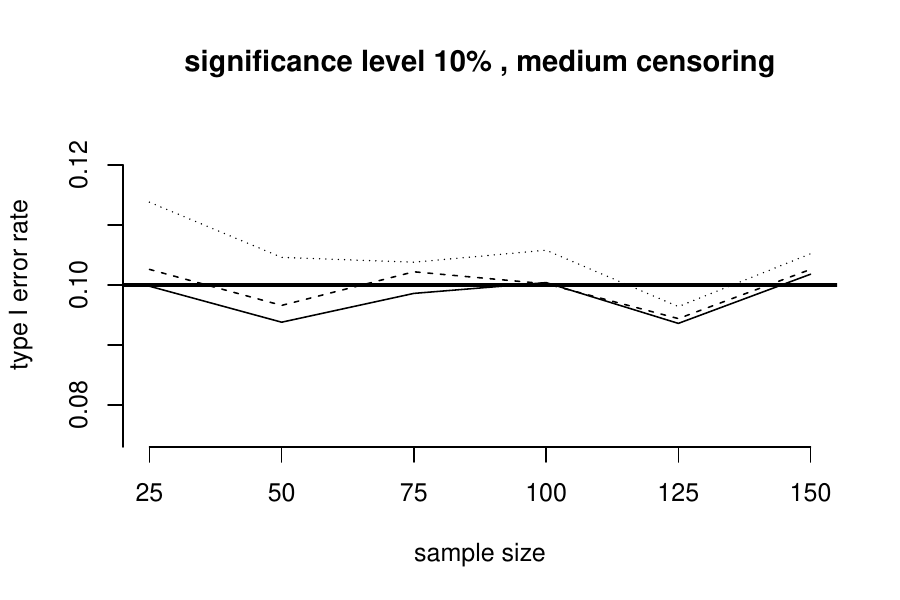}\hspace{-0.5cm}
 \includegraphics[width=0.34\textwidth]{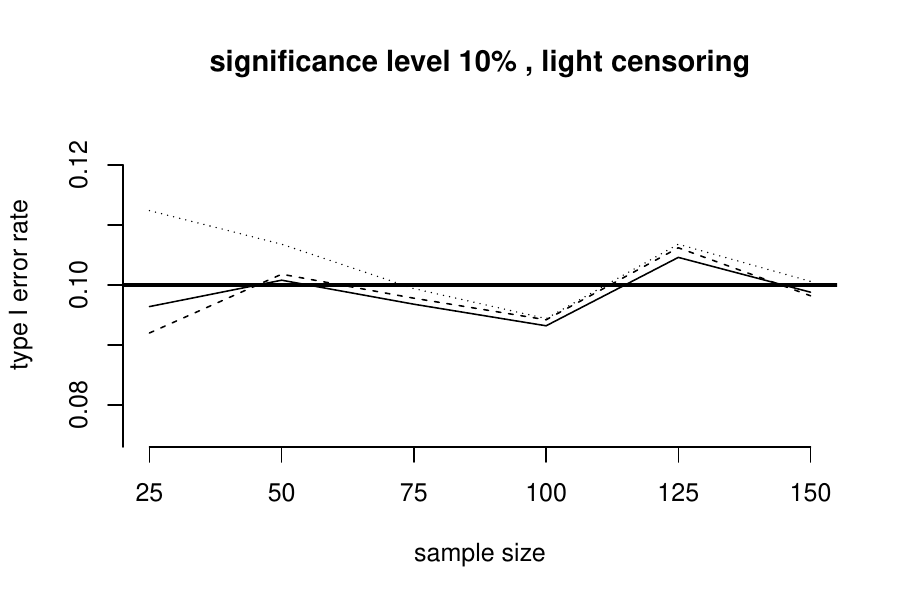} \\[-0.4cm]
 \includegraphics[width=0.34\textwidth]{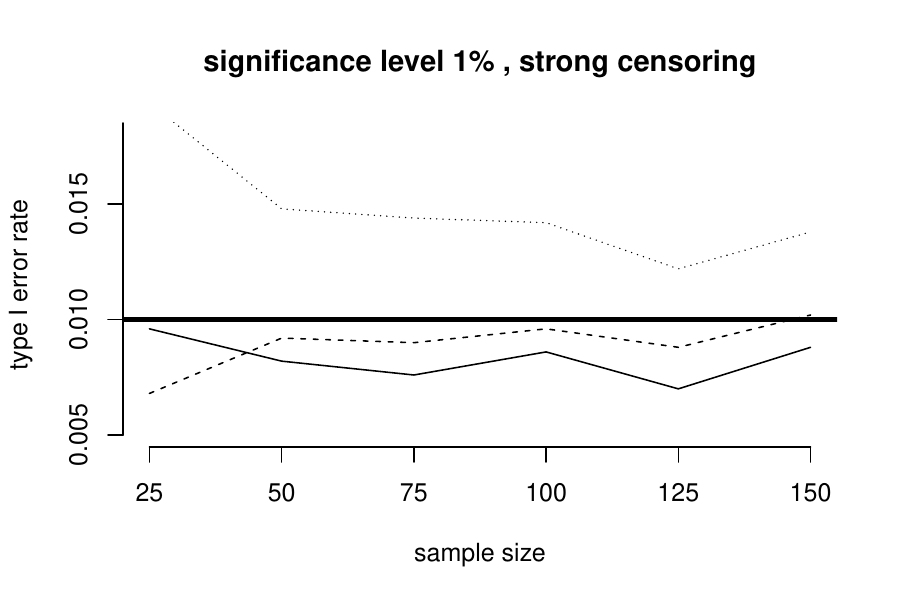}
 \hspace{-0.5cm}
 \includegraphics[width=0.34\textwidth]{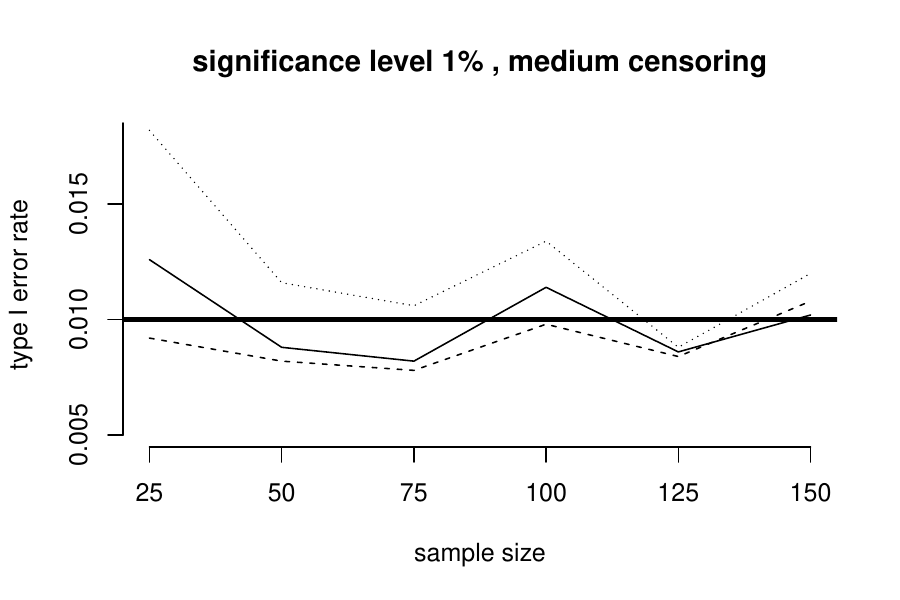}\hspace{-0.5cm}
 \includegraphics[width=0.34\textwidth]{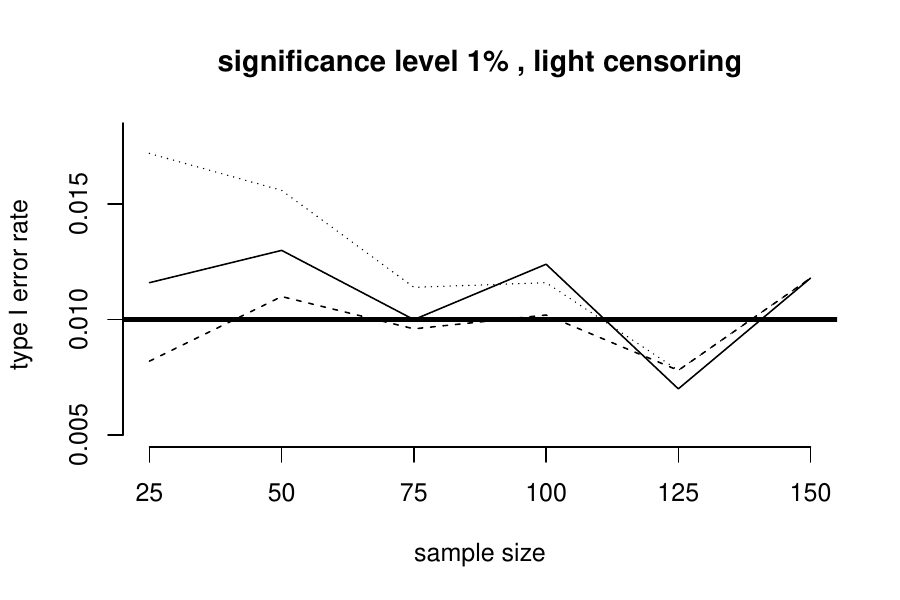} \\[-0.4cm]
 \includegraphics[width=0.34\textwidth]{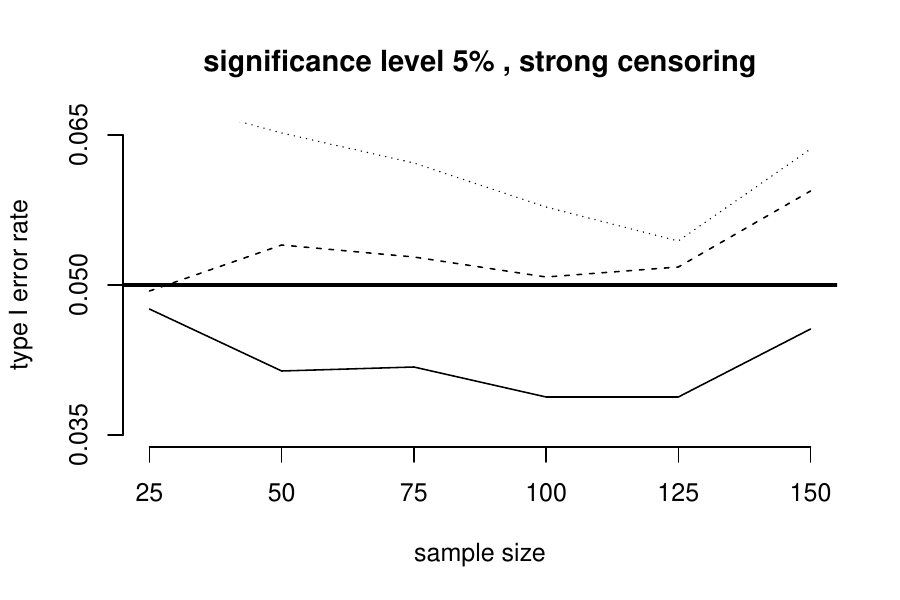}
 \hspace{-0.5cm}
 \includegraphics[width=0.34\textwidth]{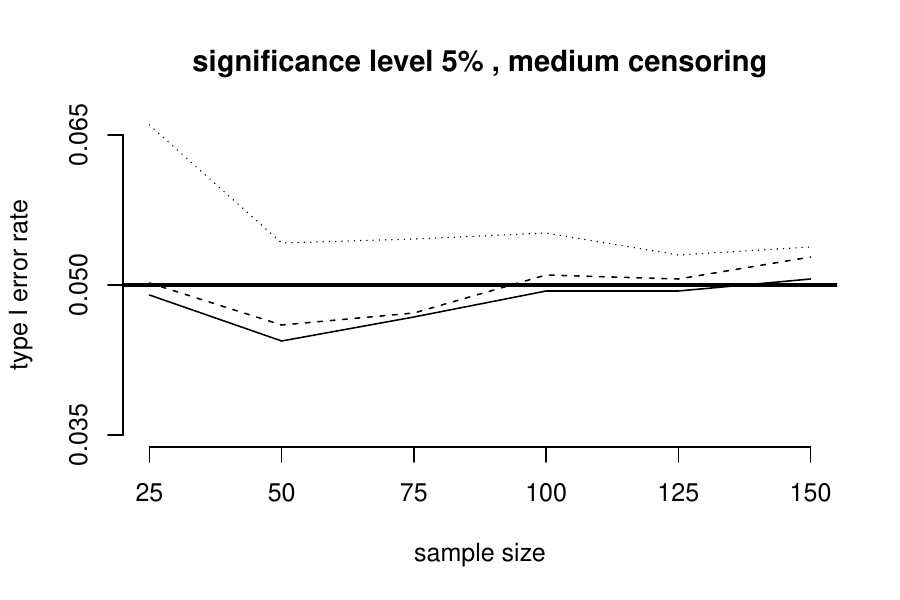}\hspace{-0.5cm}
 \includegraphics[width=0.34\textwidth]{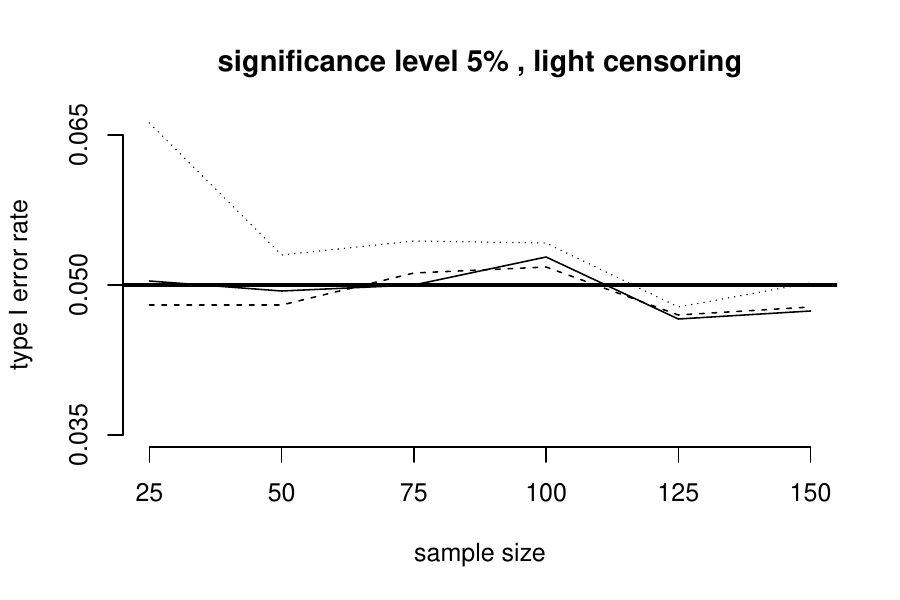} \\[-0.4cm]
 \includegraphics[width=0.34\textwidth]{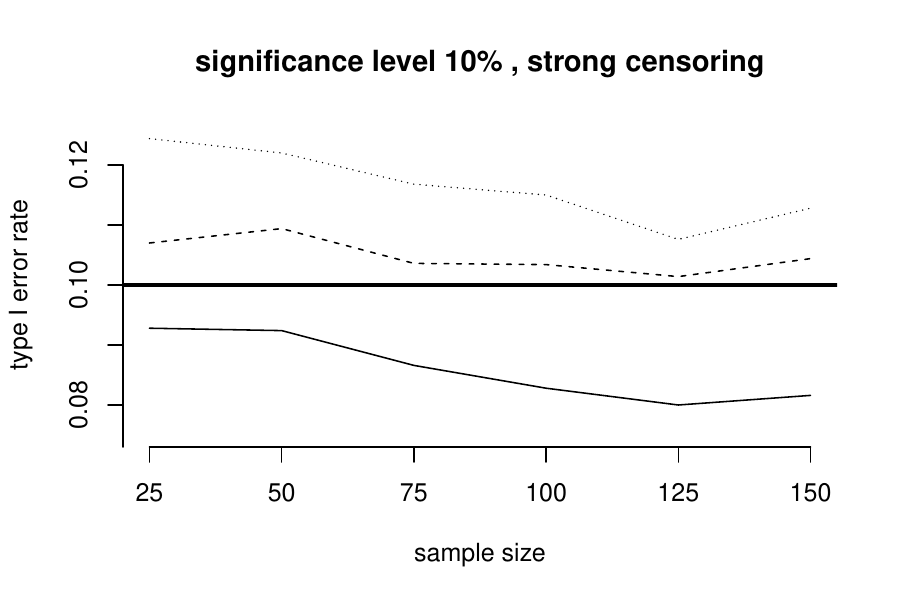}
 \hspace{-0.5cm}
 \includegraphics[width=0.34\textwidth]{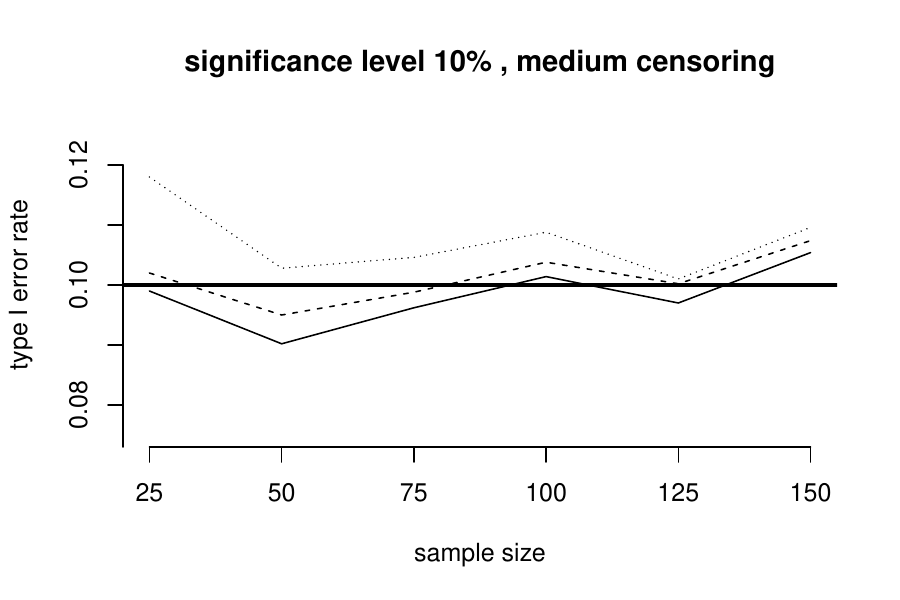}\hspace{-0.5cm}
 \includegraphics[width=0.34\textwidth]{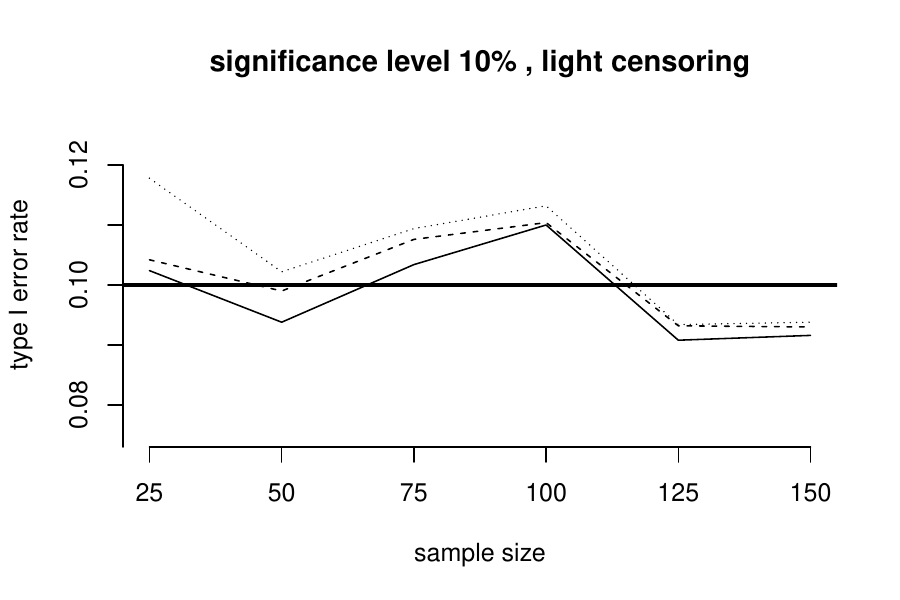}
 \vspace{-0.7cm}
 \caption{Simulated type I error rates of the Mann-Whitney-type tests with exponential versus Gompertz distributed data under strong (left), medium (middle), and light censoring (right); Gumbel-Hougaard copula; based on the randomization (---), bootstrap (- -), and asymptotic test ($\cdots$).
 The horizonal line is the nominal significance level.
 The missingness probabilities are $(\pi_1,\pi_2) = (0.125,0.125)$ (upper half) and $(\pi_1,\pi_2) = (0.1,0.2)$ (lower half).}
 \label{fig:mw_uneq}
\end{figure}

 \begin{figure}
 \includegraphics[width=0.34\textwidth]{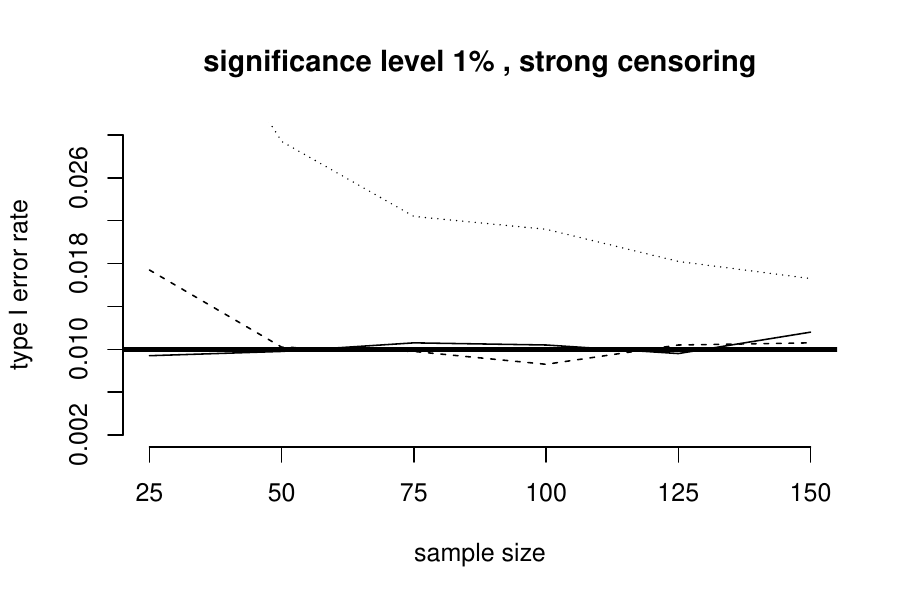}
 \hspace{-0.5cm}
 \includegraphics[width=0.34\textwidth]{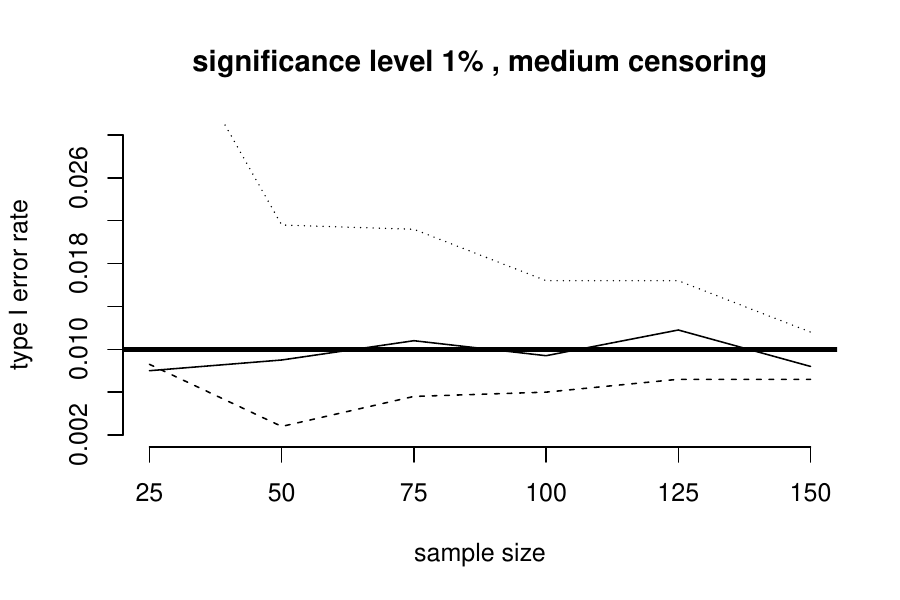}\hspace{-0.5cm}
 \includegraphics[width=0.34\textwidth]{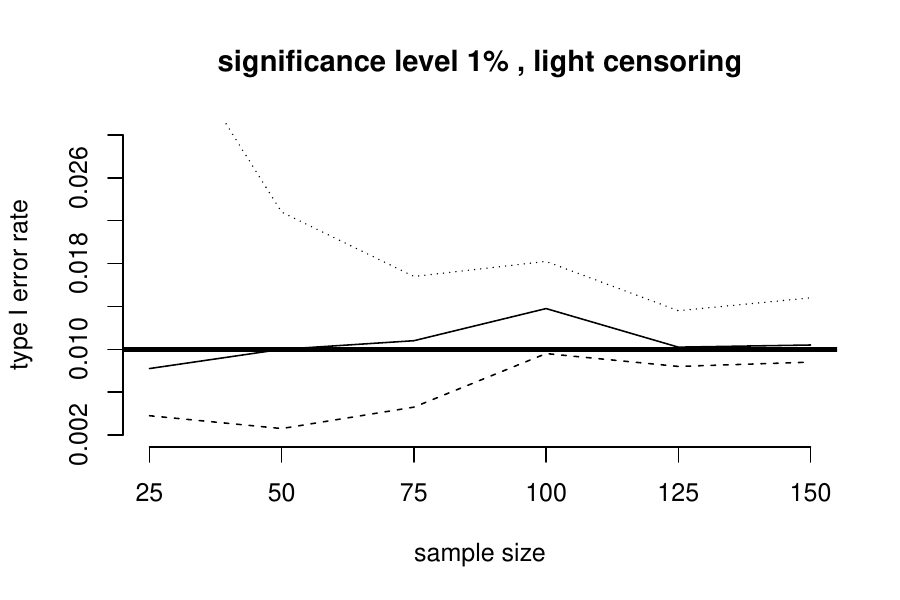} \\[-0.4cm]
 \includegraphics[width=0.34\textwidth]{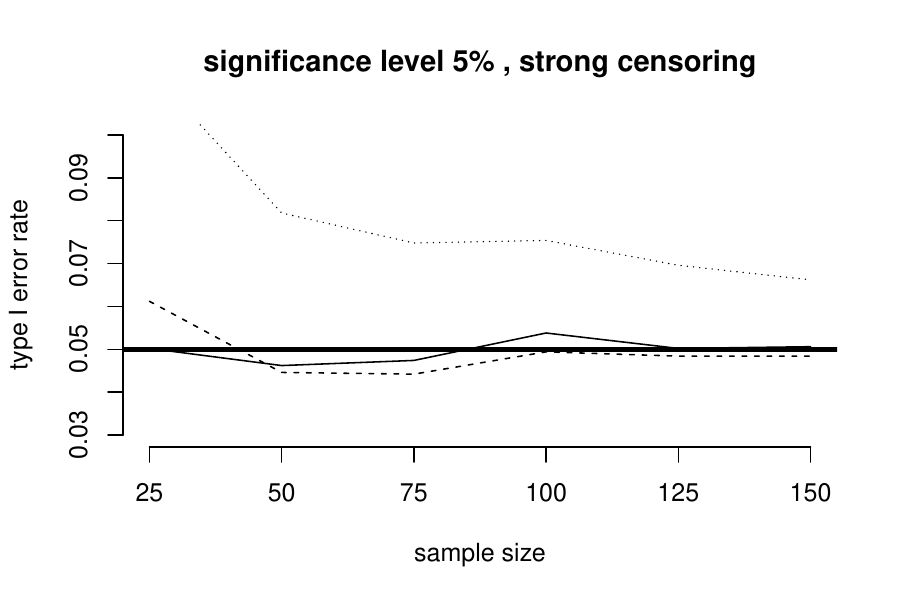}
 \hspace{-0.5cm}
 \includegraphics[width=0.34\textwidth]{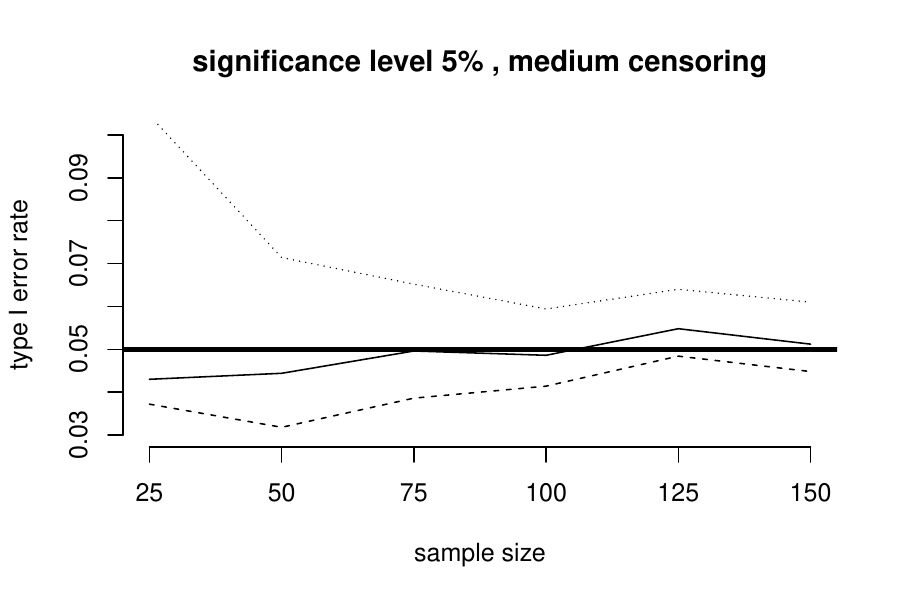}\hspace{-0.5cm}
 \includegraphics[width=0.34\textwidth]{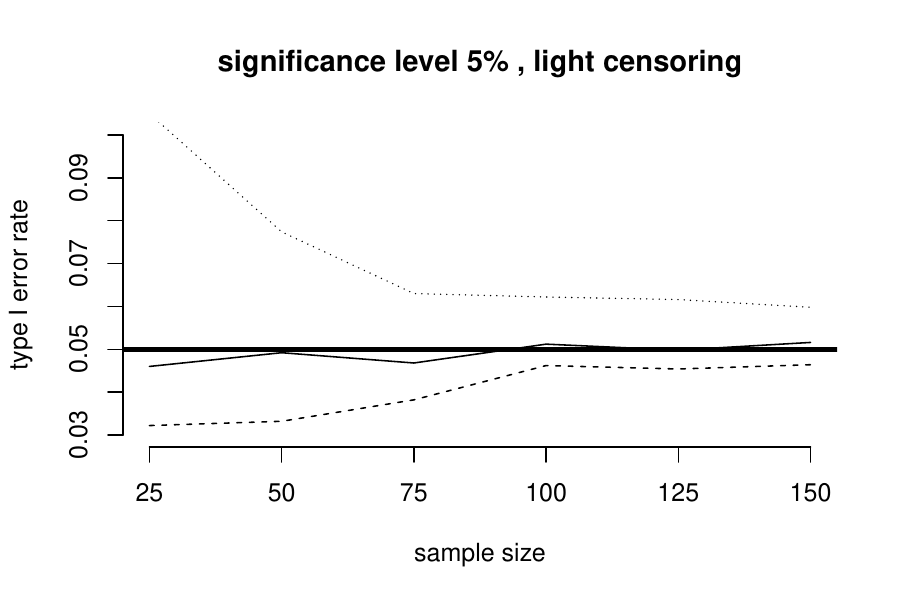} \\[-0.4cm]
 \includegraphics[width=0.34\textwidth]{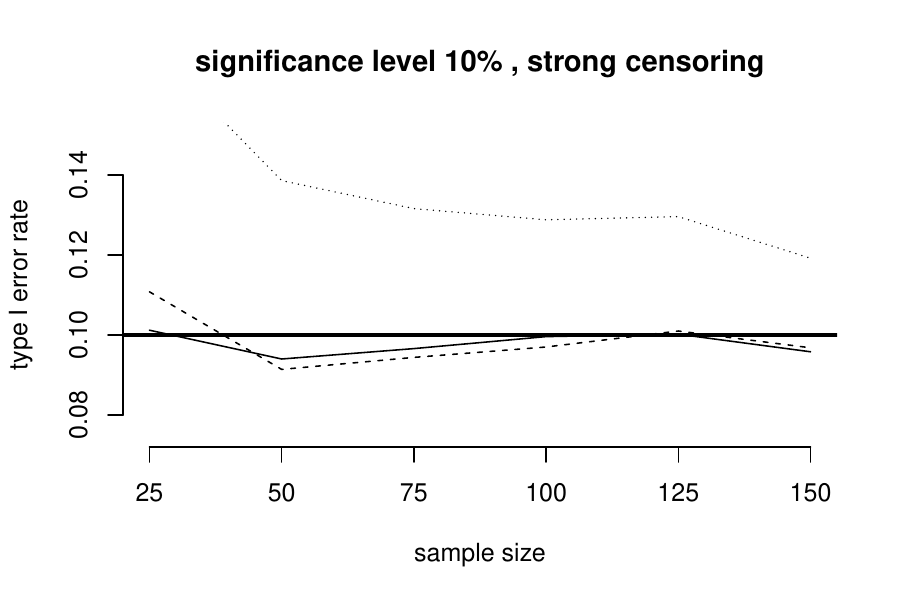}
 \hspace{-0.5cm}
 \includegraphics[width=0.34\textwidth]{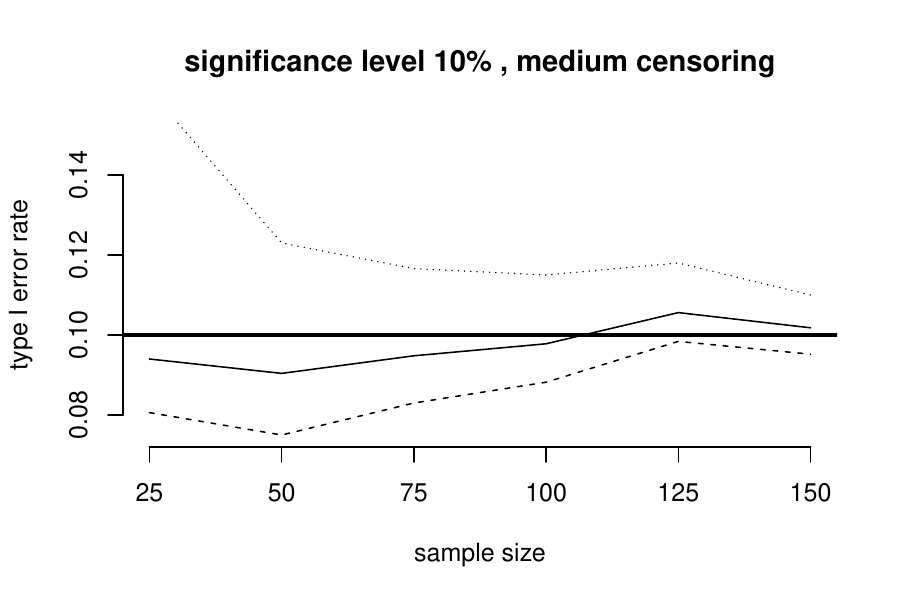}\hspace{-0.5cm}
 \includegraphics[width=0.34\textwidth]{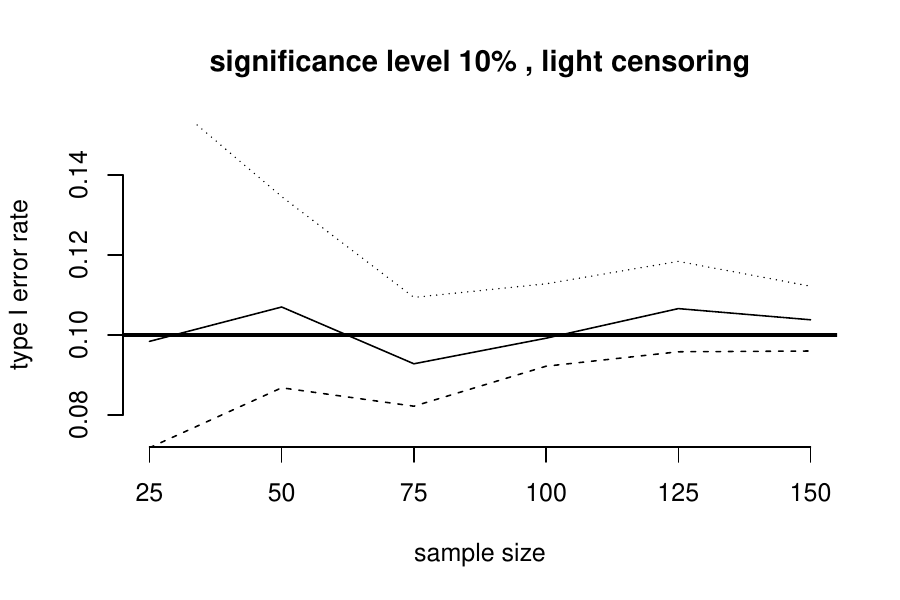}
  \\[-0.4cm]
 \includegraphics[width=0.34\textwidth]{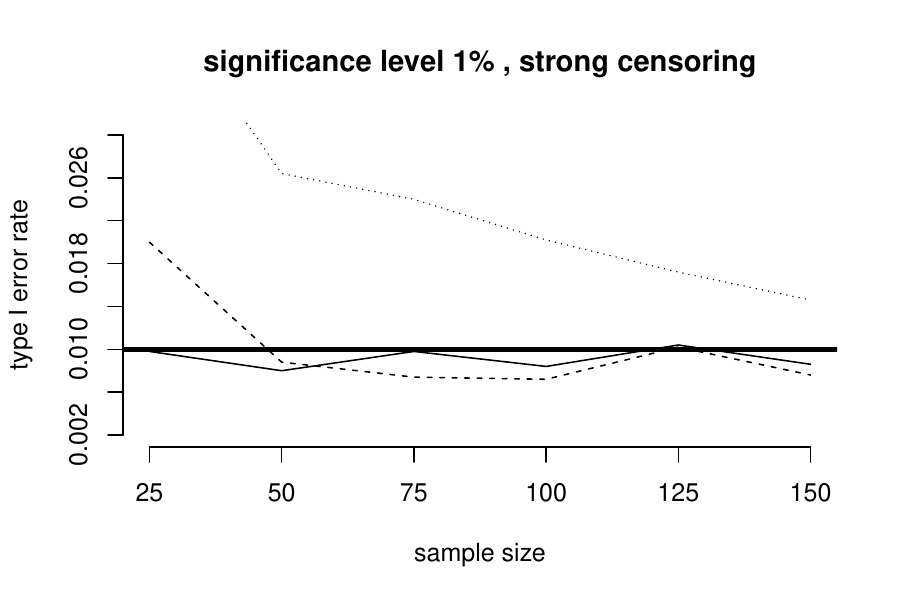}
 \hspace{-0.5cm}
 \includegraphics[width=0.34\textwidth]{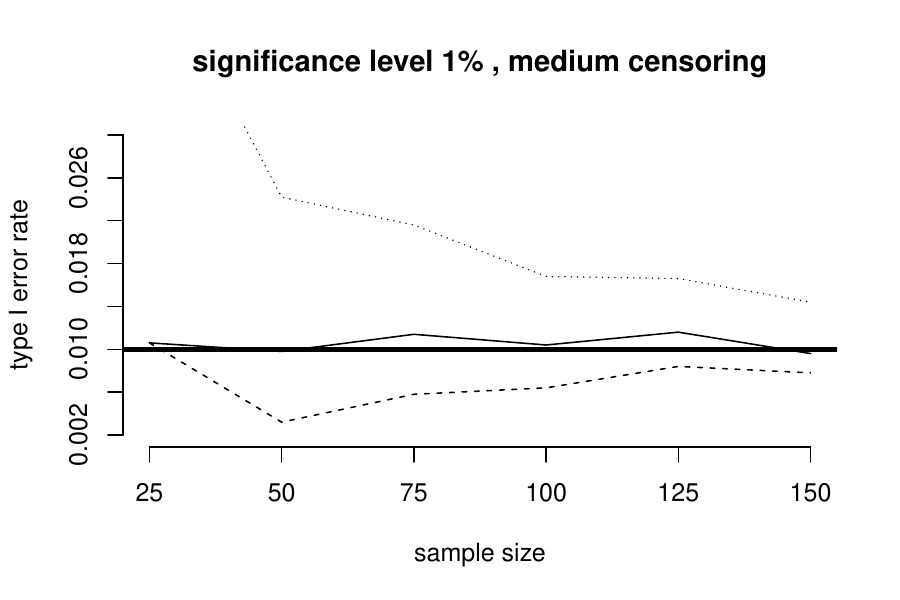}\hspace{-0.5cm}
 \includegraphics[width=0.34\textwidth]{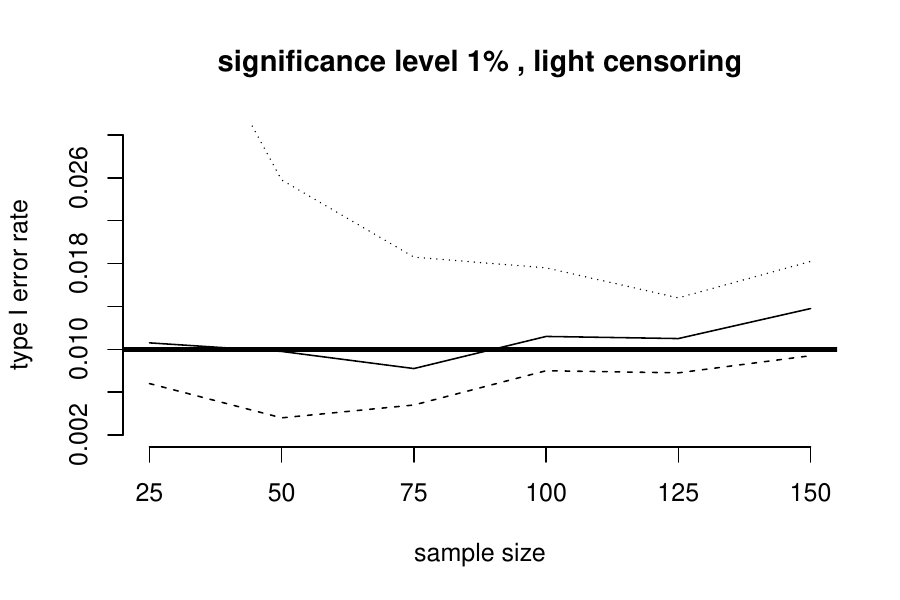} \\[-0.4cm]
 \includegraphics[width=0.34\textwidth]{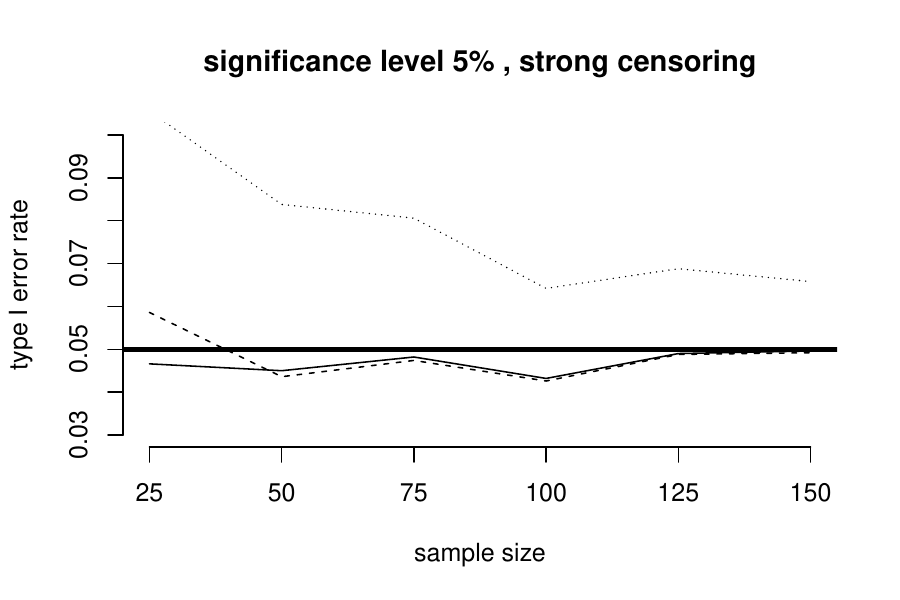}
 \hspace{-0.5cm}
 \includegraphics[width=0.34\textwidth]{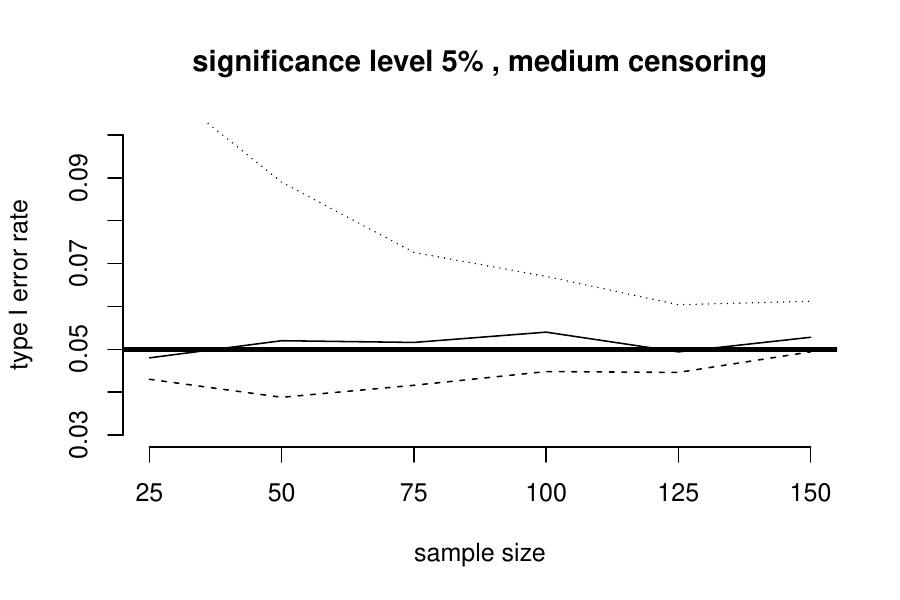}\hspace{-0.5cm}
 \includegraphics[width=0.34\textwidth]{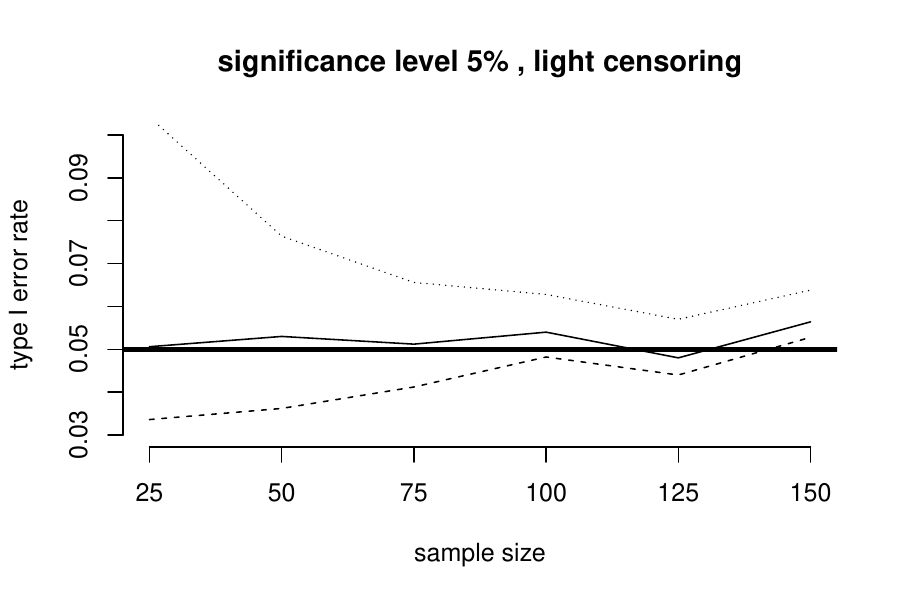} \\[-0.4cm]
 \includegraphics[width=0.34\textwidth]{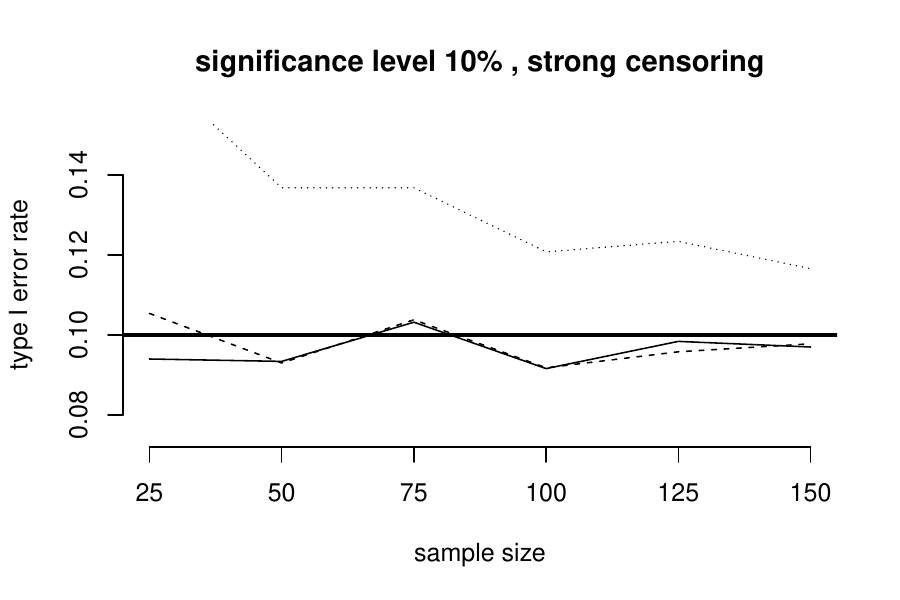}
 \hspace{-0.5cm}
 \includegraphics[width=0.34\textwidth]{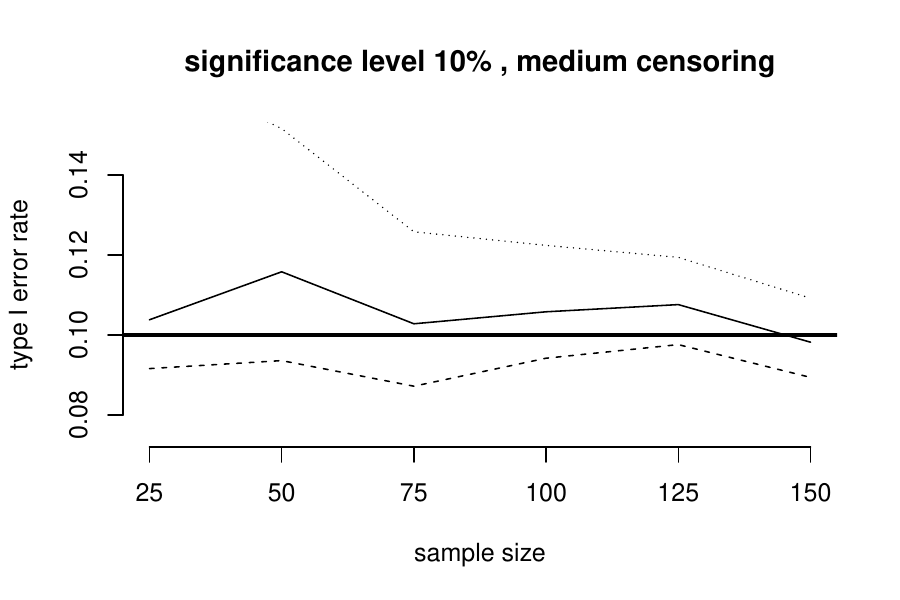}\hspace{-0.5cm}
 \includegraphics[width=0.34\textwidth]{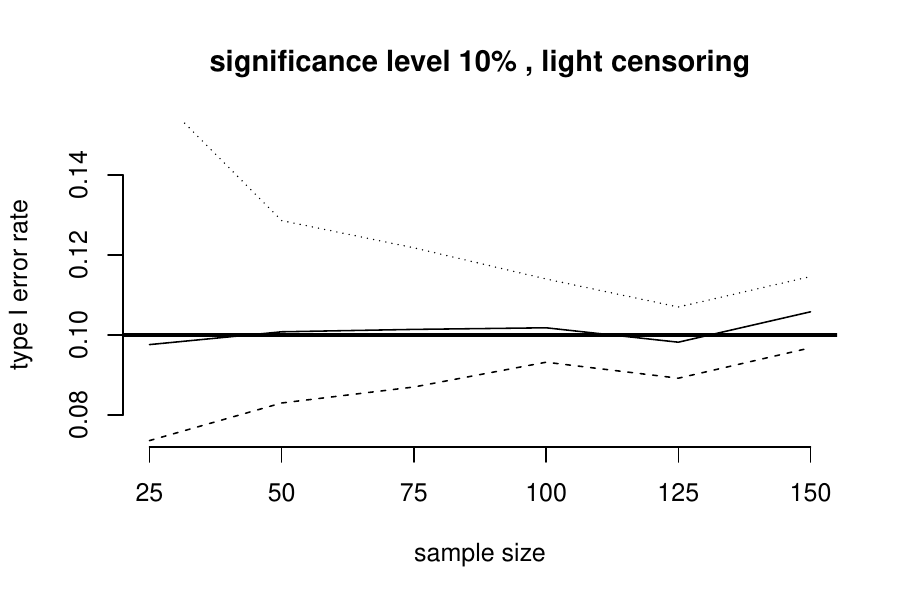}
 \vspace{-0.7cm}
 \caption{Simulated type I error rates of the Mann-Whitney-type tests with exponential versus Gompertz distributed data under strong (left), medium (middle), and light censoring (right); Clayton copula (lower half); based on the randomization (---), bootstrap (- -), and asymptotic test ($\cdots$).
 The horizonal line is the nominal significance level.
 The missingness probabilities are $(\pi_1,\pi_2) = (0.125,0.125)$ (upper half) and $(\pi_1,\pi_2) = (0.1,0.2)$ (lower half).}
 \label{fig:mw_uneq2}
\end{figure}

\begin{figure}
 \includegraphics[width=0.34\textwidth]{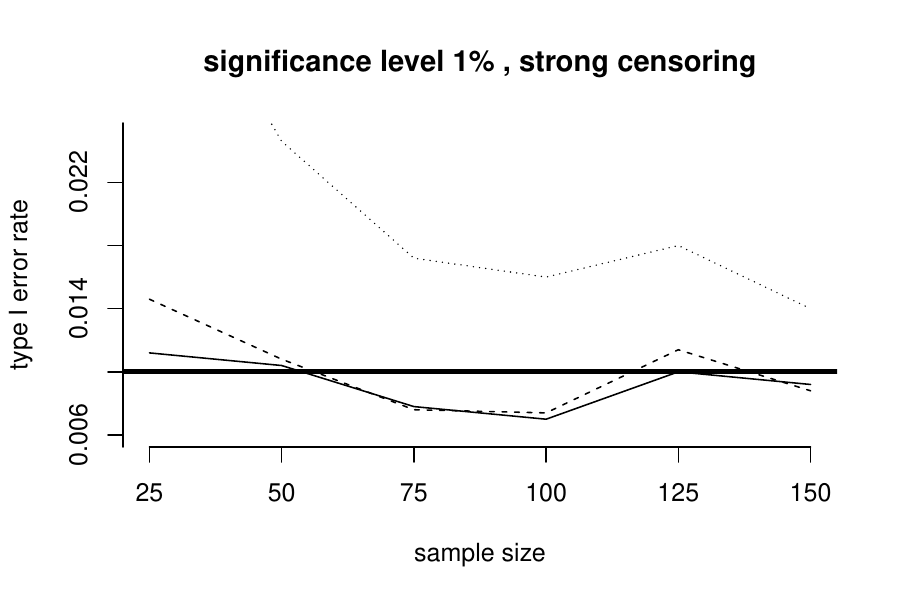}
 \hspace{-0.5cm}
 \includegraphics[width=0.34\textwidth]{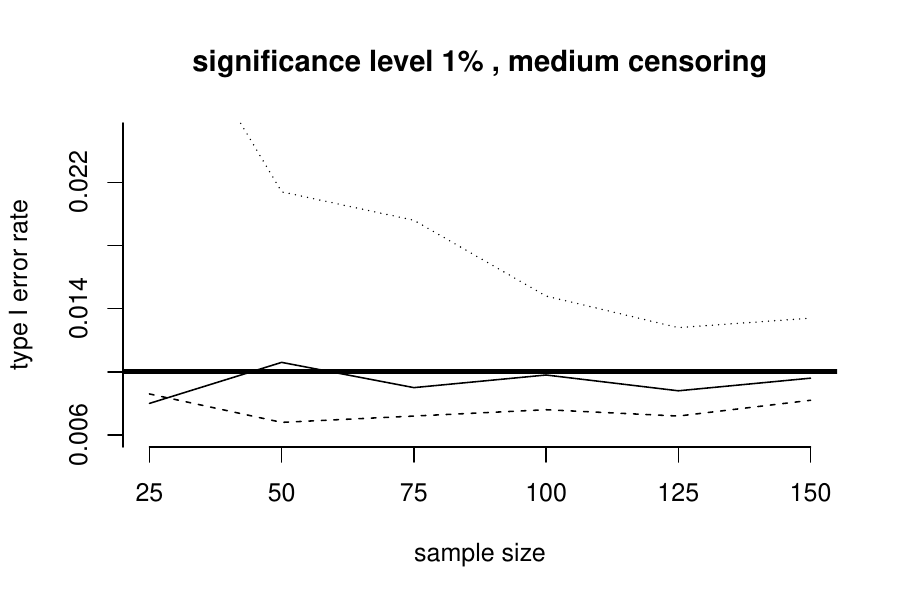}\hspace{-0.5cm}
 \includegraphics[width=0.34\textwidth]{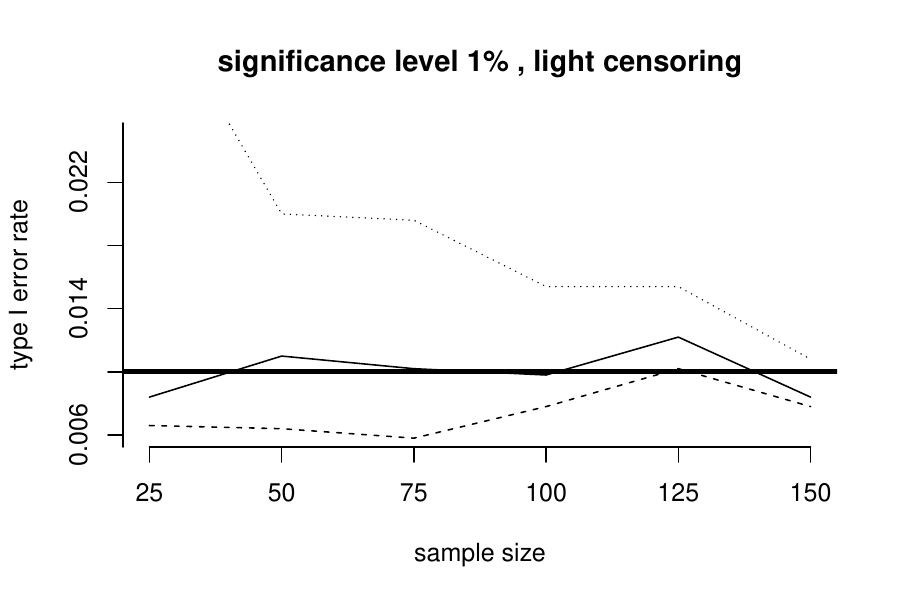} \\[-0.4cm]
 \includegraphics[width=0.34\textwidth]{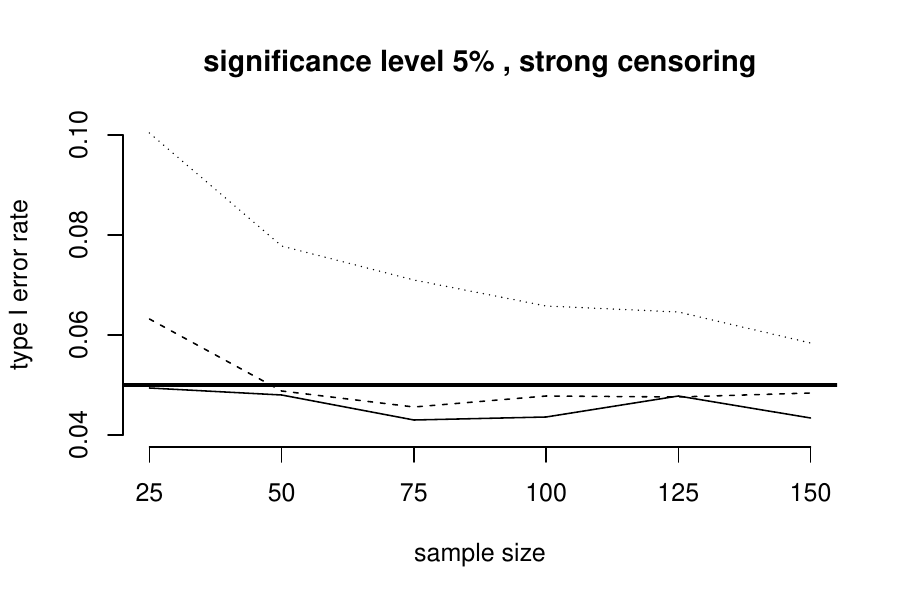}
 \hspace{-0.5cm}
 \includegraphics[width=0.34\textwidth]{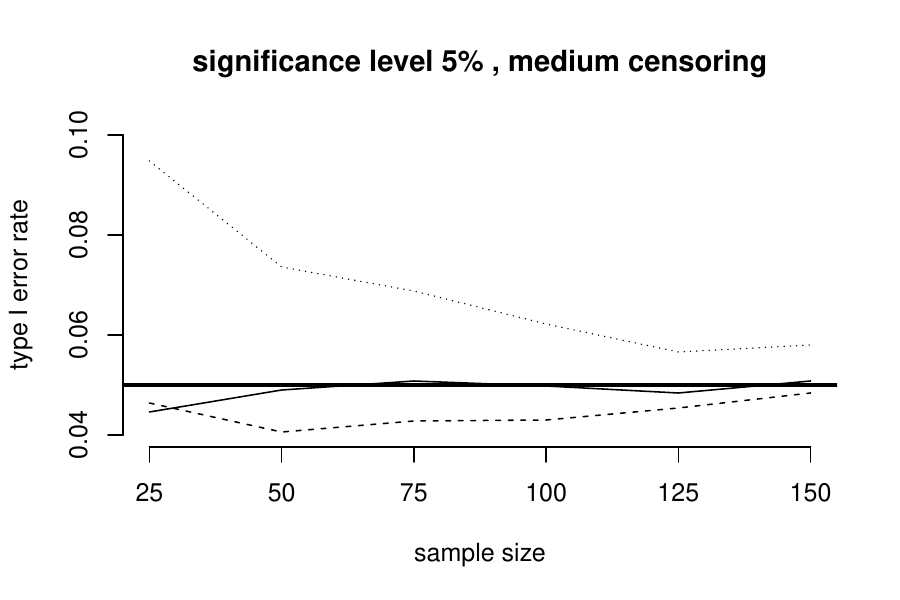}\hspace{-0.5cm}
 \includegraphics[width=0.34\textwidth]{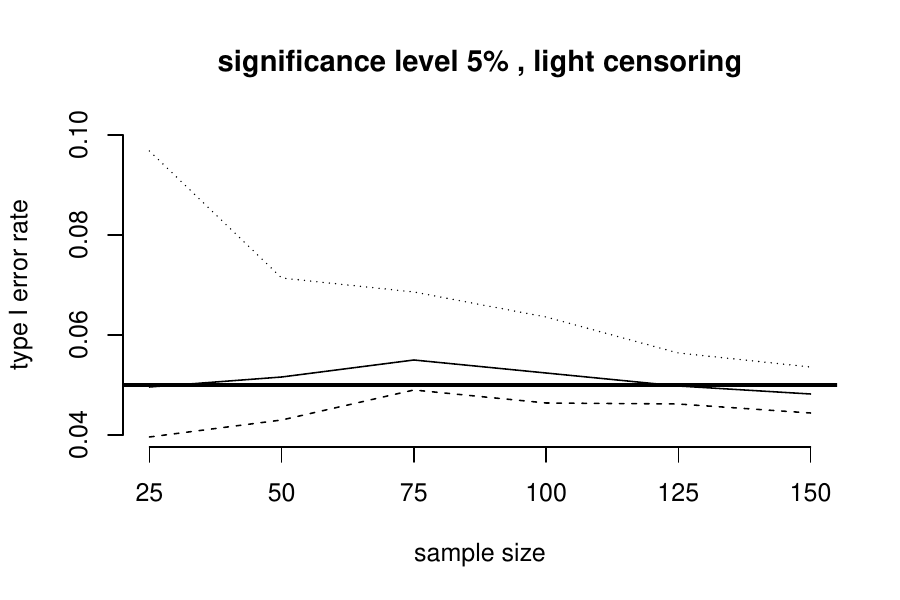} \\[-0.4cm]
 \includegraphics[width=0.34\textwidth]{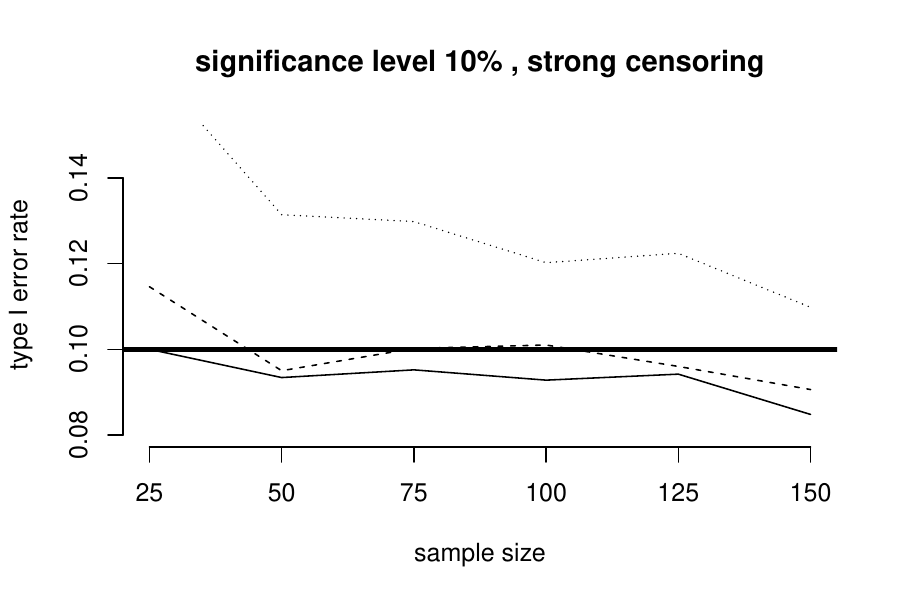}
 \hspace{-0.5cm}
 \includegraphics[width=0.34\textwidth]{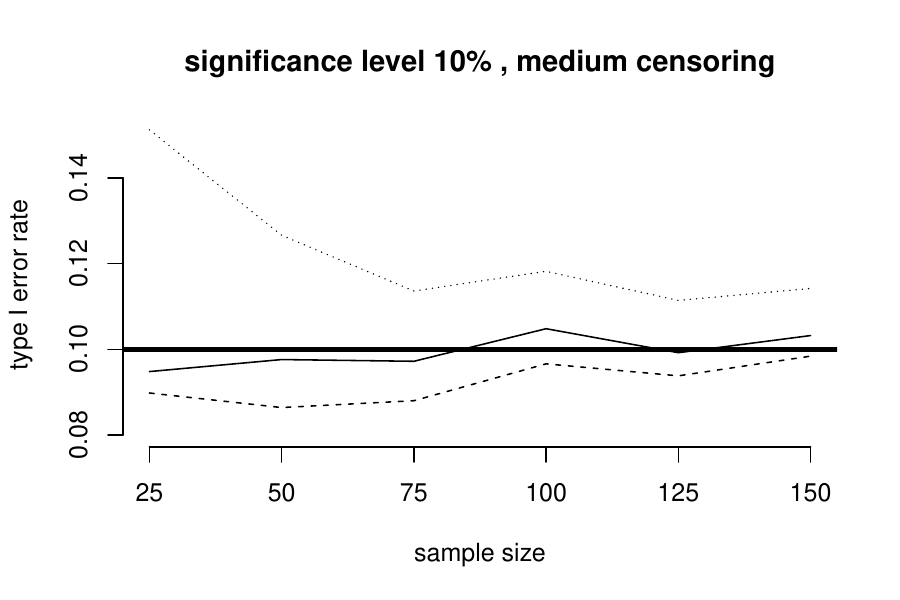}\hspace{-0.5cm}
 \includegraphics[width=0.34\textwidth]{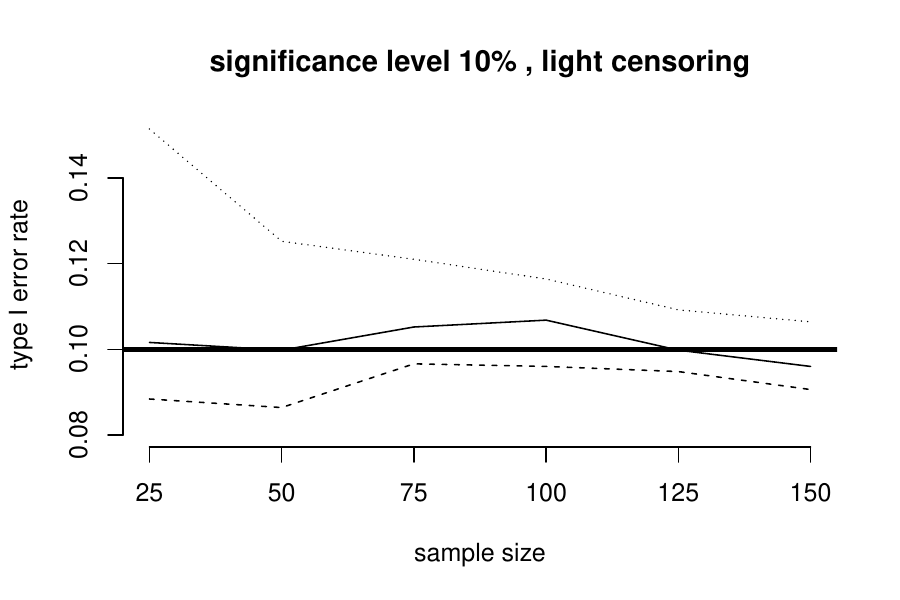}
 \\[-0.4cm]
 \includegraphics[width=0.34\textwidth]{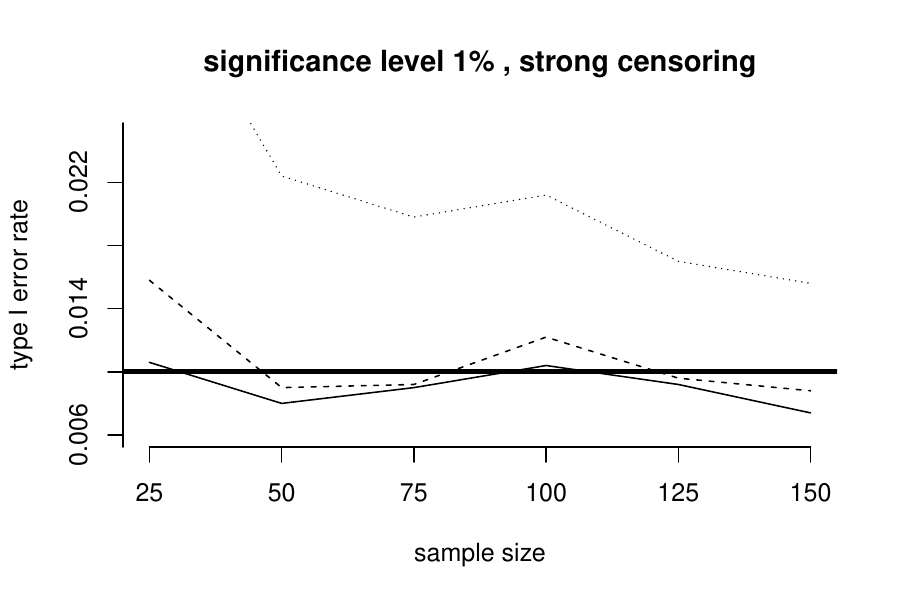}
 \hspace{-0.5cm}
 \includegraphics[width=0.34\textwidth]{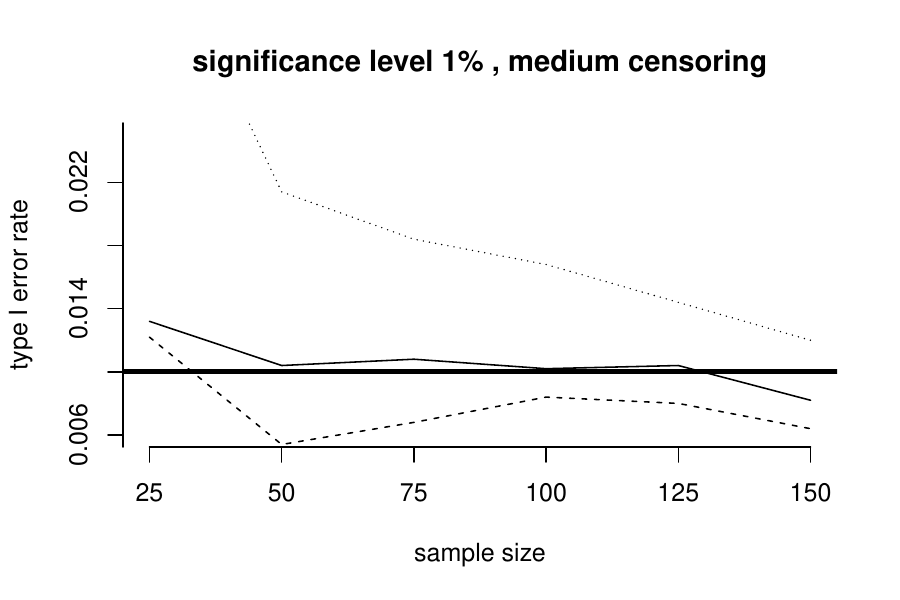}\hspace{-0.5cm}
 \includegraphics[width=0.34\textwidth]{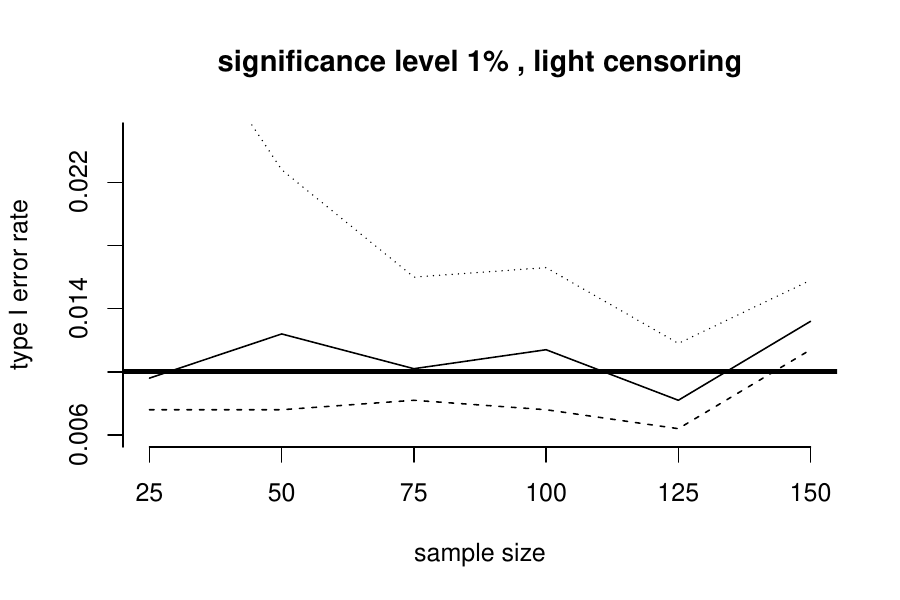} \\[-0.4cm]
 \includegraphics[width=0.34\textwidth]{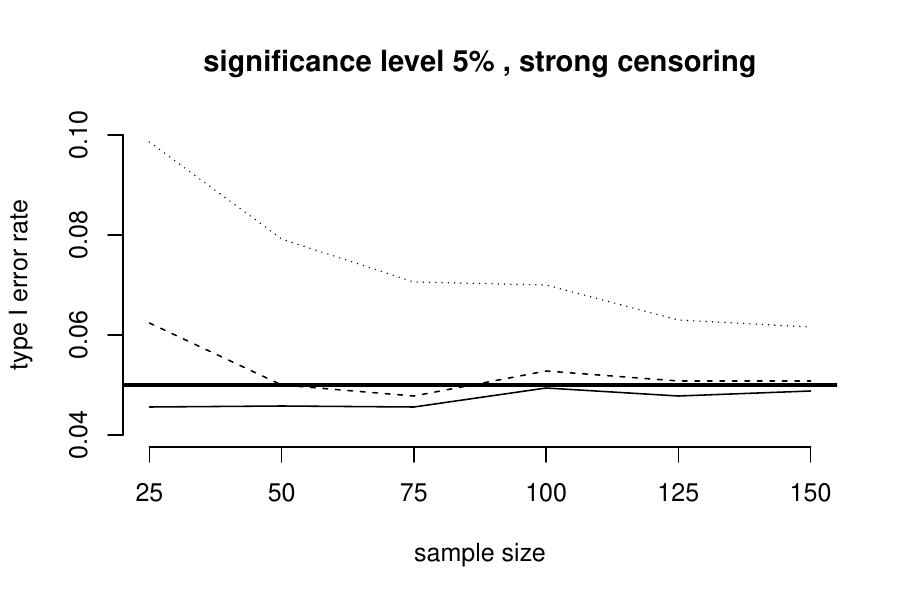}
 \hspace{-0.5cm}
 \includegraphics[width=0.34\textwidth]{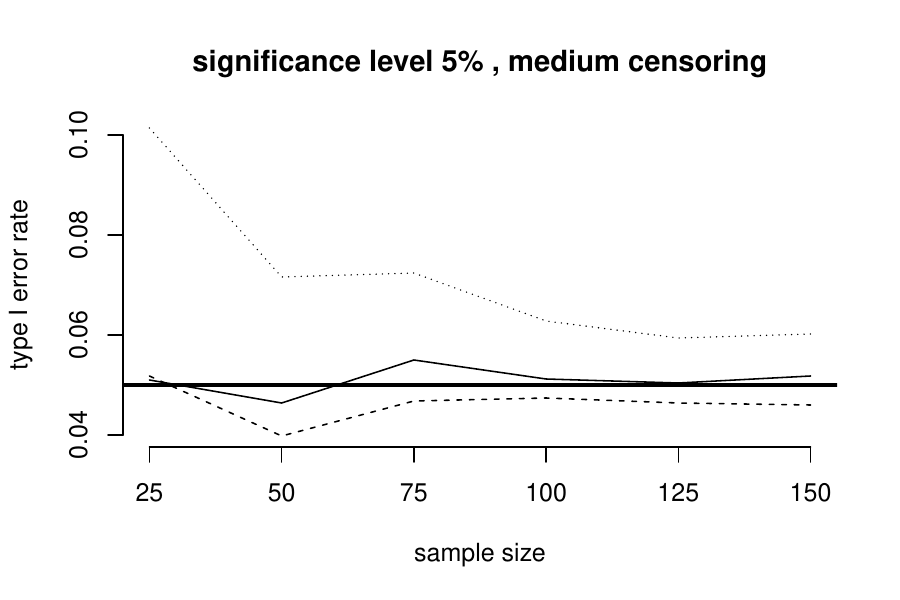}\hspace{-0.5cm}
 \includegraphics[width=0.34\textwidth]{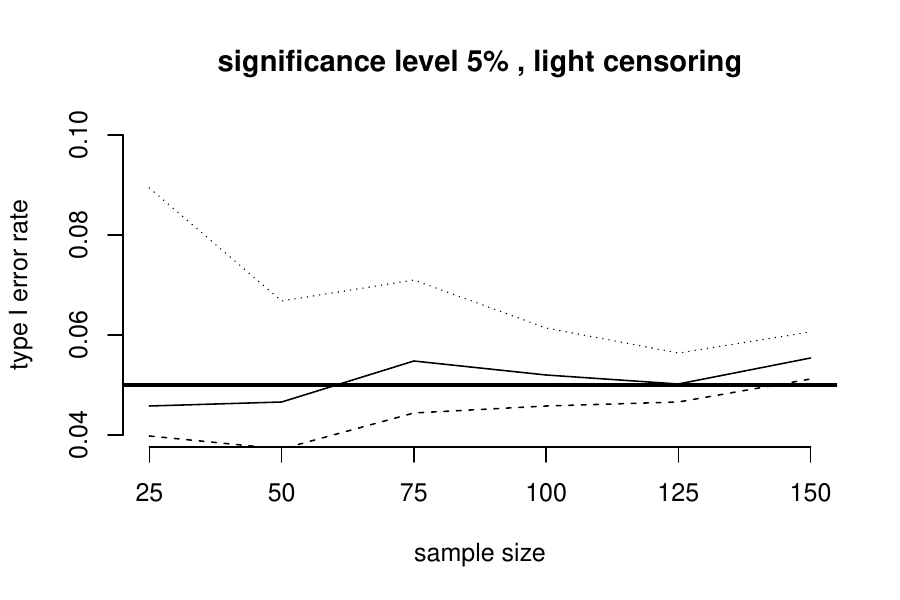} \\[-0.4cm]
 \includegraphics[width=0.34\textwidth]{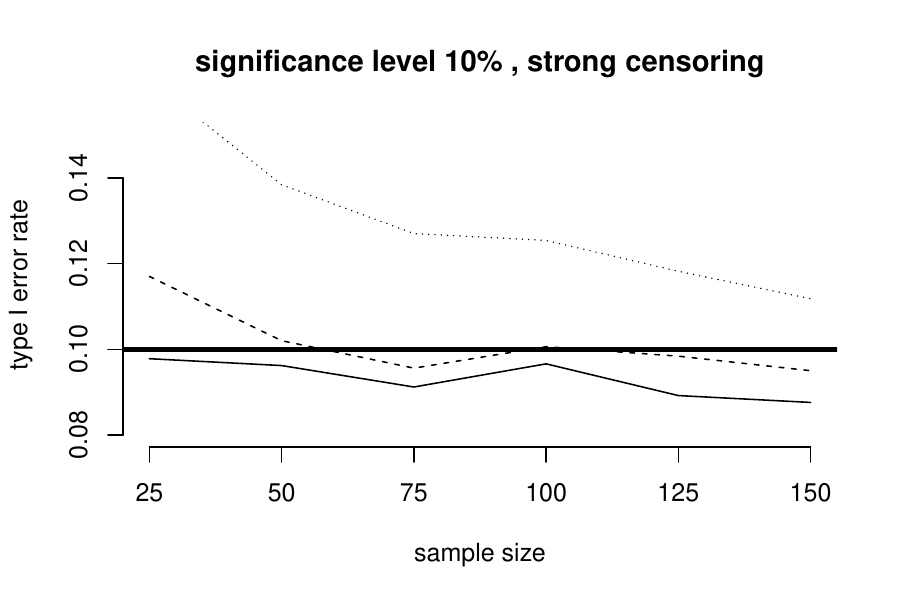}
 \hspace{-0.5cm}
 \includegraphics[width=0.34\textwidth]{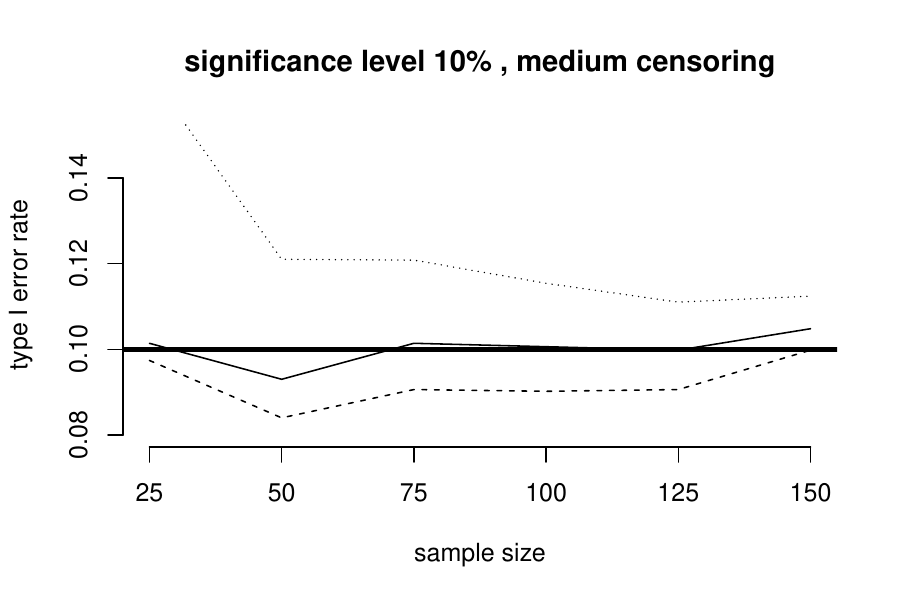}\hspace{-0.5cm}
 \includegraphics[width=0.34\textwidth]{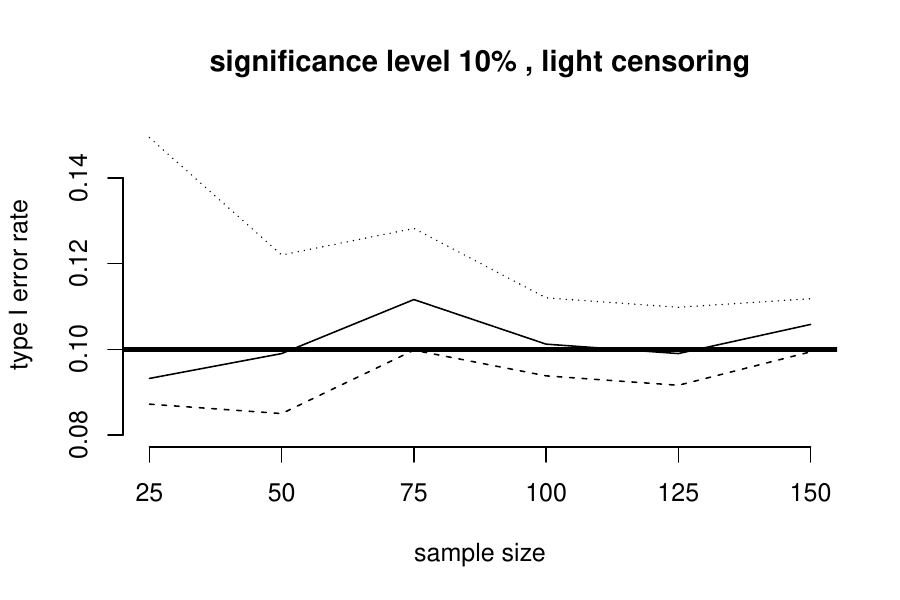}
 \vspace{-0.7cm}
 \caption{Simulated type I error rates of the Mann-Whitney-type tests with exponential versus Gompertz distributed data under strong (left), medium (middle), and light censoring (right); independence copula; based on the randomization (---), bootstrap (- -), and asymptotic test ($\cdots$). 
 The horizonal line is the nominal significance level.
 The missingness probabilities are $(\pi_1,\pi_2) = (0.125,0.125)$ (upper half) and $(\pi_1,\pi_2) = (0.1,0.2)$ (lower half).}
 \label{fig:mw_uneq3}
\end{figure}

\newpage

\section*{Acknowledgements}
The author wishes to thank Hein Putter, Liesbeth de Wreede, and Nan van Geloven for discussions on $p$ and $\check p$ and anonymous referees who  helped to improve the quality of this paper.

\bibliographystyle{plainnat}
\bibliography{literature}
\end{document}